\documentclass[10pt]{amsart}
\usepackage{amsfonts}
\usepackage{amsmath}
\usepackage{amssymb}
\usepackage{comment}
\usepackage[breaklinks,hyperindex=false]{hyperref}
\usepackage{mathrsfs}
\usepackage[noadjust]{cite}
\usepackage{enumerate}
\usepackage{latexsym}
\usepackage{xcolor}

\hypersetup{
  colorlinks = true, %Colours links instead of boxes
  urlcolor   = blue, %Colour for external hyperlinks
  linkcolor  = blue, %Colour of internal links
  citecolor  = red %Colour of citations
}

\newtheorem{theorem}{Theorem}[section]
\newtheorem{corollary}[theorem]{Corollary}
\newtheorem{lemma}[theorem]{Lemma}
\newtheorem{proposition}[theorem]{Proposition}

\theoremstyle{definition}
\newtheorem{example}[theorem]{Example}
\newtheorem{question}[theorem]{Question}
\newtheorem{remark}[theorem]{Remark}

\theoremstyle{remark}
\numberwithin{equation}{section}
\newtheorem{result}{Result}

\newcommand{\uast}{{}^*}
\newcommand{\Var}{\mathsf{Var}}
\newcommand{\fb}{finitely based}
\newcommand{\nfb}{non-finitely based}
\newcommand{\insem}{involution semigroup}
\newcommand{\Insem}{Involution semigroup}
\newcommand{\inv}{involution}
\newcommand{\ssw}{subsequence}
\newcommand{\Ablock}{ordered $A_0$-block}
\newcommand{\Bblock}{ordered $B_0$-block}
\newcommand{\sk}{\mathsf{S}}
\newcommand{\tr}{\mathsf{T}}
\newcommand{\order}{\prec}
\newcommand{\norder}{\nprec}
\newcommand{\rorder}{\succ}
\newcommand{\ordereq}{\preceq}
\newcommand{\Asf}{A}
\newcommand{\Bsf}{B}

\newcommand{\bfa}{\mathbf{a}}
\newcommand{\bfb}{\mathbf{b}}
\newcommand{\bfc}{\mathbf{c}}
\newcommand{\bfd}{\mathbf{d}}
\newcommand{\bfe}{\mathbf{e}}
\newcommand{\bff}{\mathbf{f}}
\newcommand{\bfg}{\mathbf{g}}
\newcommand{\bfh}{\mathbf{h}}
\newcommand{\bfs}{\mathbf{s}}
\newcommand{\bft}{\mathbf{t}}
\newcommand{\bfu}{\mathbf{u}}
\newcommand{\bfv}{\mathbf{v}}
\newcommand{\bfw}{\mathbf{w}}
\newcommand{\frakM}{\mathfrak{M}}

\newcommand{\fta}{\textsc}
\newcommand{\rma}{\fta{a}}

\newcommand{\rme}{\fta{e}}
\newcommand{\rmef}{\fta{ef}}
\newcommand{\rmf}{\fta{f}}

\newcommand{\cX}{\mathcal{X}}
\newcommand{\cXX}{\cX\cup\cX^*}
\newcommand{\Fsem}{F_\mathsf{inv}(\cX)}
\newcommand{\Fmon}{F_\mathsf{inv}^1(\cX)}

\newcommand{\up}{\textup}
\newcommand{\sfcon}{\mathsf{con}}
\newcommand{\sfint}{\mathsf{int}}
\newcommand{\sfocc}{\mathsf{occ}}

\newcommand{\sfh}{\mathsf{h}}
\newcommand{\sft}{\mathsf{t}}
\newcommand{\sfH}{\mathsf{H}}
\newcommand{\sfT}{\mathsf{T}}

\newcommand{\headu}{h} \newcommand{\tailu}{t} \newcommand{\headv}{H} \newcommand{\tailv}{T}

\begin{document}

\allowdisplaybreaks

\title[Finite basis problem for involution semigroups of order four]{Finite basis problem for involution semigroups of order four}

%\thanks{This research was partially supported by the National Natural Science Foundation of China (Nos. 12271224, 12171213) and the Fundamental Research Funds for the Central Universities (No. lzujbky-2023-ey16).}

\thanks{This research was partially supported by the National Natural Science Foundation of China (nos. 12271224, 12571018, and 12171213) and the Fundamental Research Funds for the Central Universities (no. lzujbky-2023-ey06).}

\author[M. Gao]{Meng Gao}
\address{School of Mathematics and Statistics, Lanzhou University, Lanzhou, Gansu 730000, P. R. China}
\email{gaom2015@lzu.edu.cn}

\author[E. W. H. Lee]{Edmond W. H. Lee}
\address{Department of Mathematics, Nova Southeastern University, Fort Lauderdale, FL 33328, USA}
\email{edmond.lee@nova.edu}

\author[Y. F. Luo]{Yan Feng Luo}
\address{School of Mathematics and Statistics, Lanzhou University, Lanzhou, Gansu 730000, P. R. China}
\email{luoyf@lzu.edu.cn}

\author[W. T. Zhang]{Wen Ting Zhang$^\star$}\thanks{$^\star$Corresponding author}
\address{School of Mathematics and Statistics, Lanzhou University, Lanzhou, Gansu 730000, P. R. China}
\email{zhangwt@lzu.edu.cn}

\subjclass[2020]{20M05}
\keywords{semigroup, involution semigroup, identity basis, finitely based, finite basis problem}

\begin{abstract}
Since the 1980s, it has been known that the smallest non-finitely based semigroups are of order six.
Surprisingly, for involution semigroups, a non-finitely based example of order five was recently discovered.
In this article, it is confirmed that every involution semigroup of order four is finitely based.
Since every involution semigroup of order three or less is already known to be finitely based, it follows that the smallest non-finitely based involution semigroups are of order five.
\end{abstract}

\maketitle

\section{Introduction} \label{sec: intro}

\subsection{Minimal {\nfb} {\insem}s}

An \textit{identity basis} for an algebra~$A$ is a set of identities of~$A$ that axiomatizes all the identities of~$A$.
An algebra is \textit{\fb} if it has some finite identity basis; otherwise, it is \textit{\nfb}.
A prominent research problem in universal algebra is the \textit{finite basis problem}: determine which finite algebras are {\fb}.
Finite groups~\cite{OP64}, finite associative rings \cite{BO75}, finite Lie rings \cite{Kru73,Lvo73}, and finite lattices~\cite{McK70} are {\fb}, but in general, not all finite algebras are {\fb}.
For instance, the multiplicative matrix semigroup \[ B_2^1 = \left\{ \begin{bmatrix}\,0\,\,\,0\,\\\,0\,\,\,0\,\end{bmatrix}, \begin{bmatrix}\,1\,\,\,0\,\\\,0\,\,\,0\,\end{bmatrix}, \begin{bmatrix}\,0\,\,\,1\,\\\,0\,\,\,0\,\end{bmatrix}, \begin{bmatrix}\,0\,\,\,0\,\\\,1\,\,\,0\,\end{bmatrix}, \begin{bmatrix}\,0\,\,\,0\,\\\,0\,\,\,1\,\end{bmatrix}, \begin{bmatrix}\,1\,\,\,0\,\\\,0\,\,\,1\,\end{bmatrix} \right\}, \] published by Perkins in 1969, is {\nfb}~\cite{Per69}.
The discovery of this example focused much attention upon the finite basis problem for small semigroups.
Decades of cumulative work that followed has shown that every semigroup of order five or less is {\fb} \cite{Lee13,Tra83,Tra91}, and among all semigroups of order six---28,634 of them up to isomorphism~\cite{DK14}---only four are {\nfb} \cite{LL11,LLZ12,LZ15}.
The four {\nfb} semigroups of order six, which include~$B_2^1$, are thus \textit{minimal \nfb}.

The present article is concerned with \textit{{\insem}s} $(S,\uast)$, that is, semigroups~$S$ with a unary operation~$\uast$ that satisfy the identities \begin{equation} (x^*)^* \approx  x, \quad (xy)^* \approx y^*x^*; \label{id: inv} \end{equation} the unary operation~$\uast$ is called an \textit{\inv} of~$S$.
An \textit{inverse semigroup} is an {\insem} $(S,\uast)$ that satisfies the additional identities \[ xx^*x \approx x, \quad xx^*yy^* \approx yy^*xx^*. \]
Examples of inverse semigroups include any group $(G,{}^{-1})$ with inversion~$^{-1}$ and the Perkins semigroup $(B_2^1,{}^\tr)$ under the usual matrix transposition~$^\tr$.
Examples of {\insem}s that are not inverse semigroups include the multiplicative ${n \times n}$ matrix semigroup $(M_n(\mathbb{F}),{}^\tr)$ over any field~$\mathbb{F}$ with the usual transposition~$^\tr$ and the Perkins semigroup $(B_2^1,{}^\sk)$ under the \textit{skew transposition}~$^\sk$ across the secondary diagonal, that is, \[ \begin{bmatrix}\,a\,\,\,b\,\\\,c\,\,\,d\,\end{bmatrix}^\sk = \begin{bmatrix}\,d\,\,\,b\,\\\,c\,\,\,a\,\end{bmatrix}. \]

Given how close {\insem}s are to semigroups, it seems reasonable to conjecture that a finite {\insem} $(S,\uast)$ and its semigroup \textit{reduct}~$S$ are always simultaneously {\fb}.
For instance, the {\insem}s $(B_2^1,{}^\tr)$ and $(B_2^1,{}^\sk)$ are both {\nfb}~\cite{ADV12,Kle79}, while their reduct~$B_2^1$ is also {\nfb}~\cite{Per69}.
However, this conjecture has been refuted by several counterexamples \cite{GZL20,JV10,Lee16,Lee19}, the smallest of which is the multiplicative matrix semigroup \[ A_0^1 = \left\{\begin{bmatrix}\,0\,\,0\,\,0\,\\\,0\,\,0\,\,0\,\\\,0\,\,0\,\,0\,\end{bmatrix}, \begin{bmatrix}\,1\,\,1\,\,0\,\\\,0\,\,0\,\,0\,\\\,0\,\,0\,\,0\,\end{bmatrix}, \begin{bmatrix}\,0\,\,0\,\,1\,\\\,0\,\,0\,\,0\,\\\,0\,\,0\,\,0\,\end{bmatrix}, \begin{bmatrix}\,0\,\,0\,\,0\,\\\,0\,\,0\,\,1\,\\\,0\,\,0\,\,1\,\end{bmatrix}, \begin{bmatrix}\,1\,\,0\,\,0\,\\\,0\,\,1\,\,0\,\\\,0\,\,0\,\,1\,\end{bmatrix}\right\} \] under the skew transposition~$^\sk$.
As noted above, all semigroups of order five or less---which include~$A_0^1$---are {\fb}, but the {\insem} $(A_0^1,{}^\sk)$ is {\nfb}~\cite{GZL20}.
It follows that minimal {\nfb} {\insem}s are of order at most five; this result is quite unexpected given that minimal {\nfb} semigroups are  of order six~\cite{LLZ12}.

It is of fundamental importance to examine if there exists a {\nfb} {\insem} that is smaller than $(A_0^1,{}^\sk)$.
Since every {\insem} of order three or less is {\fb}~\cite{Lee22}, the answer would require addressing the finite basis problem for those of order four.
Solving this problem is the objective of the present article.

\begin{theorem} \label{Thm: main}
Every {\insem} of order four is {\fb}.
\end{theorem}

Consequently, minimal {\nfb} {\insem}s are of order five, and $(A_0^1,{}^\sk)$ is one such example.
An obvious next step in the investigation is to question the uniqueness of this example.

\begin{question} \label{Q: example}
Is there a {\nfb} {\insem} of order five that is not isomorphic to $(A_0^1,{}^\sk)$?
\end{question}

It has recently been confirmed that up to isomorphism, $(A_0^1,{}^\sk)$ is the unique smallest {\nfb} {\insem} within the class of all {\insem}s with a unit element~\cite{HZL}.
Therefore, to answer Question~\ref{Q: example}, it suffices to only examine {\insem}s without a unit element.

\subsection{Finite basis problem for finite (\inv) semigroups}

Let $\frakM_n$ denote the set of all subsemigroups of $M_n(\mathbb{R})$ consisting of binary matrices and let $\frakM_\infty = \bigcup_{n \geq 1} \frakM_n$.
It is not a coincidence that all explicit examples of finite {\insem}s given so far are semigroups from $\frakM_\infty$ with transpositions~$^\tr$ or~$^\sk$, given that every finite semigroup is isomorphic to some semigroup in $\frakM_\infty$ and every finite inverse semigroup is isomorphic to some semigroup in $\frakM_\infty$ with the usual transposition~$^\tr$; see, for instance, Howie \cite[Theorems~1.1.2 and~5.1.7]{How95}.

Regarding finite {\insem}s in general, it turns out that every one of them is isomorphic to some semigroup in $\frakM_\infty$ with the skew transposition~$^\sk$ but not necessarily the usual transposition~$^\tr$ \cite{Lee23b}.
Therefore, when addressing the finite basis problem for finite semigroups (with \inv)---which is currently open---it is equivalent to focus on finite semigroups in $\frakM_\infty$ (with the skew transposition~$^\sk$).
Refer to the survey by Volkov~\cite{Vol01} for more information on the finite basis problem for finite semigroups.

On the other hand, the finite basis problem for finite algebras is undecidable in general~\cite{McK96}.

\subsection{Organization}

Notation and background information are first given in Section \ref{sec: prelim}.
An outline of the proof of Theorem~\ref{Thm: main} is then given in Section~\ref{sec: proof}, while the finer details of the proof are deferred to Sections~\ref{sec: S2}--\ref{sec: B0}.
Multiplication tables of all {\insem}s of order up to four are listed in Section~\ref{sec: tables}.

\section{Preliminaries} \label{sec: prelim}

Most of the notation and background results of this article are given in this section.
Refer to the monograph of Burris and Sankappanavar~\cite{BS81} for any undefined notation and terminology of universal algebra in general.

\subsection{Words}

Let~$\ordereq$ be a total order on a countably infinite alphabet~$\cX$ that excludes the symbol~$0$; write $x \order y$ to indicate that $x \ordereq y$ and $x \neq y$.

Let $\cX^* =\{x^* \,|\, x \in \cX \}$ be a disjoint copy of~$\cX$.
Elements of $\cXX$ are called \textit{variables}.
The \textit{free \insem} over~$\cX$ is the free semigroup $\Fsem = (\cXX)^+$ with unary operation given by $(x^*)^* = x$ for all $x \in \cX$ and \[ (x_1x_2 \cdots x_{n-1} x_n)^* = x_n^* x_{n-1}^* \cdots x_2^*x_1^* \] for all $x_1,x_2,\ldots,x_{n-1},x_n \in \cXX$.
The \textit{free {\inv} monoid} over~$\cX$ is $\Fmon = \Fsem \cup \{1\}$, where~$1$ is the empty word with $1^* = 1$.
Elements of $\Fmon$ are called \textit{words} and words in $\cX^+ \cup \{1\}$ are said to be \textit{plain}.

Any word $\bfw \in \Fmon$ can be written in the form \[ \bfw = x_1^{\circledast_1} x_2^{\circledast_2} \cdots x_n^{\circledast_n} \] for some $x_1, x_2,\ldots, x_n \in \cX$ and $\circledast_1,\circledast_2,\ldots,\circledast_n \in \{ 1,*\}$ with $n \geq 0$; the \textit{plain projection} of such a word is the plain word \[ \overline{\bfw}=x_1x_2\cdots x_n. \]

For any words $\bfu,\bfv \in \Fmon$, write $\bfu \hookrightarrow \bfv$ to indicate that $\bfu$ is a \textit{\ssw} of~$\bfv$, that is, $\bfu = x_1 x_2 \cdots x_n$ for some $x_1, x_2, \ldots, x_n \in \cXX$ and \[ \bfv = \bfv_0 x_1 \bfv_1 x_2 \bfv_2 \cdots x_n \bfv_n \] for some $\bfv_0, \bfv_1, \ldots, \bfv_n \in \Fmon$.
Specifically, a {\ssw} $\bfu$ of $\bfv$ such that $\sfcon(\overline{\bfu}) = \{x,y\}$ for some $x,y \in \cX$ is called an \textit{$\{x,y\}$-\ssw} of $\bfv$.

For any word $\bfw \in \Fsem$, the \textit{content} of~$\bfw$, denoted by $\sfcon (\bfw)$, is the set of variables occurring in~$\bfw$; the number of times that a variable $x \in \cXX$ occurs in~$\bfw$ is denoted by $\sfocc(x,\bfw)$; the \textit{head} of~$\bfw$, denoted by $\sfh(\bfw)$, is the first variable occurring in $\bfw$; the \textit{tail} of~$\bfw$, denoted by $\sft(\bfw)$, is the last variable occurring in $\bfw$; and the length of $\bfw$ is denoted by $|\bfw|$.

\begin{example}
If $\bfw =x^*xy^*x^2yz^*x^*y$ for some $x,y,z\in \cX$, then
\begin{itemize}
  \item $\overline{\bfw}=x^2yx^2yzxy$;
  \item $\sfcon(\bfw)=\{x, x^*, y, y^*, z^*\}$;
  \item $\sfocc(x, \bfw) = 3$,\, $\sfocc(x^*,\bfw) =2$,\, $\sfocc(y,\bfw) = 2$,\, $\sfocc(y^*,\bfw) = \sfocc(z^*,\bfw)=1$;
  \item $\sfh(\bfw)=x^*$,\, $\sft(\bfw)=y$; and
  \item $|\bfw|=9$.
\end{itemize}
\end{example}

For any word $\bfw \in \Fsem$, a variable $x \in \sfcon(\bfw)$ is \textit{simple} if $\sfocc(\overline{x},\overline{\bfw}) =1$.
A word $\bfw$ is \textit{simple} if every variable in $\bfw$ is simple.
If $x,x^* \in \sfcon(\bfw)$ for some $x \in \cX$, then $\{x, x^*\}$ is a \textit{mixed pair} of $\bfw$.
A word $\bfw$ is \textit{mixed} if it has some mixed pair; otherwise, $\bfw$ is \textit{bipartite}.

Two words $\bfw_1, \bfw_2 \in \Fmon$ are \textit{disjoint} if $\sfcon(\overline{\bfw_1}) \cap \sfcon(\overline{\bfw_2}) = \emptyset$.
A non-simple word~$\bfw$ is \textit{connected} if it cannot be decomposed into a product of two disjoint nonempty words.

\subsection{Identities}

An \textit{identity} is an expression $\bfu \approx \bfv$ formed by words $\bfu, \bfv \in \Fsem$; it is \textit{nontrivial} if $\bfu \neq \bfv$.
An identity $\bfu \approx \bfv$ is \textit{bipartite} if both $\bfu$ and $\bfv$ are bipartite words.
A bipartite identity $\bfu \approx \bfv$ is \textit{plain} if $\bfu,\bfv \in \cX^+$.

An {\insem} $(S,\uast)$ \textit{satisfies} an identity $\bfs \approx \bft$, or $\bfs \approx \bft$ is \textit{satisfied by} $(S,\uast)$, if for any substitution $\varphi: \cX \to S$, the elements $\varphi(\bfs)$ and $\varphi(\bft)$ of $S$ coincide; in this case, $\bfs \approx \bft$ is also said to be an \textit{identity of} $(S,\uast)$.

\begin{lemma}[Lee {\cite[Lemma~9]{Lee17a}}] \label{Lem: A0 mixed bipartite}
An {\insem} satisfies a bipartite identity $\bfu \approx \bfv$ with $\sfcon({\bfu})=\sfcon({\bfv})$ if and only if it satisfies $\overline{\bfu} \approx \overline{\bfv}$.
\end{lemma}

\begin{lemma}[Lee {\cite[Lemma~2.12]{Lee23a}}] \label{Lem: reversal}
An {\insem} satisfies an identity $\bfu \approx \bfv$ if and only if it satisfies the identity $\overleftarrow{\bfu} \approx \overleftarrow{\bfv}$\textup, where $\overleftarrow{\bfu}$ and $\overleftarrow{\bfv}$ are the words~$\bfu$ and~$\bfv$ written in reverse order.
\end{lemma}

Recall that a \textit{semilattice} is a semigroup that is commutative and idempotent.
Up to isomorphism, the smallest semilattice with nontrivial {\inv} is the multiplicative matrix semigroup \[ S\ell_3 = \left\{ \begin{bmatrix}\,0\,\,\,0\,\\\,0\,\,\,0\,\end{bmatrix}, \begin{bmatrix}\,1\,\,\,0\,\\\,0\,\,\,0\,\end{bmatrix}, \begin{bmatrix}\,0\,\,\,0\,\\\,0\,\,\,1\,\end{bmatrix} \right\} \] under the skew transposition~$^\sk$.

\begin{lemma}[Lee {\cite[Lemma~8]{Lee17a}}] \label{Lem: A0 twisted}
Let $\bfu \approx \bfv$ be any identity of $(S\ell_3,{}^\sk)$.
Then $\bfu$ is bipartite if and only if $\bfv$ is bipartite\up; in this case\up, $\sfcon({\bfu})=\sfcon({\bfv})$.
\end{lemma}

A set $\Sigma$ of identities of $(S,\uast)$ is an \textit{identity basis} for $(S,\uast)$ if every identity of $(S,\uast)$ is deducible from $\Sigma$.
An {\insem} is \textit{\fb} if it possesses a finite identity basis.
It is unambiguous and sometimes convenient to take the {\inv} axioms~\eqref{id: inv} for granted and omit them from an identity basis for an {\insem}.

\section{Proof of Theorem~\ref{Thm: main}} \label{sec: proof}

Every finite commutative {\insem} is {\fb} \cite[Proposition~2.2]{HZL}.
Since every {\insem} of order three or less is commutative~\cite{Lee22}, they are all {\fb}.
Therefore, it remains to consider only noncommutative {\insem}s of order four.
With the help of a computer, it is routine to check that up to isomorphism, there exist 83 {\insem}s of order four; see Section~\ref{sec: tables}.
Only six of these 83 {\insem}s are noncommutative:
\begin{align*}
&
\begin{tabular}{c | c c c c}
$S_1$ & 1 & 2 & 3 & 4 \\ \hline
1     & 1 & 1 & 1 & 1 \\
2     & 1 & 1 & 1 & 1 \\
3     & 1 & 1 & 1 & 1 \\
4     & 1 & 1 & 2 & 1 \\ \hline\hline
$x$   & 1 & 2 & 3 & 4 \\ \hline
$x^*$ & 1 & 2 & 4 & 3
\end{tabular}
&&
\begin{tabular}{c | c c c c}
$S_2$ & 1 & 2 & 3 & 4 \\ \hline
1     & 1 & 1 & 1 & 1 \\
2     & 1 & 1 & 1 & 1 \\
3     & 1 & 1 & 1 & 3 \\
4     & 1 & 2 & 1 & 4 \\ \hline\hline
$x$   & 1 & 2 & 3 & 4 \\ \hline
$x^*$ & 1 & 3 & 2 & 4
\end{tabular}
&&
\begin{tabular}{c | c c c c}
$S_3$ & 1 & 2 & 3 & 4 \\ \hline
1     & 1 & 1 & 1 & 1 \\
2     & 1 & 1 & 1 & 1 \\
3     & 1 & 1 & 2 & 1 \\
4     & 1 & 1 & 2 & 2 \\ \hline\hline
$x$   & 1 & 2 & 3 & 4 \\ \hline
$x^*$ & 1 & 2 & 4 & 3
\end{tabular}
\\[0.05in]
&
\begin{tabular}{c | c c c c}
$S_4$ & 1 & 2 & 3 & 4 \\ \hline
1     & 1 & 1 & 1 & 1 \\
2     & 1 & 1 & 1 & 2 \\
3     & 1 & 2 & 3 & 1 \\
4     & 1 & 1 & 1 & 4 \\ \hline\hline
$x$   & 1 & 2 & 3 & 4 \\ \hline
$x^*$ & 1 & 2 & 4 & 3
\end{tabular}
&&
\begin{tabular}{c | c c c c}
$S_5$ & 1 & 2 & 3 & 4 \\ \hline
1     & 1 & 1 & 1 & 1 \\
2     & 1 & 1 & 1 & 2 \\
3     & 1 & 2 & 3 & 2 \\
4     & 1 & 1 & 1 & 4 \\ \hline\hline
$x$   & 1 & 2 & 3 & 4 \\ \hline
$x^*$ & 1 & 2 & 4 & 3
\end{tabular}
&&
\begin{tabular}{c | c c c c}
$S_6$ & 1 & 2 & 3 & 4 \\ \hline
1     & 1 & 1 & 3 & 3 \\
2     & 2 & 2 & 4 & 4 \\
3     & 1 & 1 & 3 & 3 \\
4     & 2 & 2 & 4 & 4 \\ \hline\hline
$x$   & 1 & 2 & 3 & 4 \\ \hline
$x^*$ & 1 & 3 & 2 & 4
\end{tabular}
\end{align*}

%The {\insem}s $(S_1,\uast)$ and $(S_3,\uast)$ are nilpotent and so are easily shown to be \fb, while the identities of the {\insem} $(S_6,\uast)$ are axiomatized by $\{ x^2 \approx x,\, xyz \approx xz \}$ \cite[Lemma~2]{Faj71}.
Since $(S_1,\uast)$ and $(S_3,\uast)$ satisfy the identity $x_1x_2x_3 \approx y_1y_2y_3$, their identities can be axiomatized by those formed by words of length at most four, whence they are finitely based.
The {\insem}s $(S_2,\uast)$, $(S_5,\uast)$, and $(S_4,\uast)$ are shown to be {\fb} in Sections~\ref{sec: S2}, \ref{sec: A0}, and~\ref{sec: B0}, respectively.
%The {\insem} $(S_6,\uast)$---a rectangular band with \inv---is long known to be finitely based; its identities are axiomatized by $\{ x^2 \approx x,\, xyz \approx xz \}$ \cite[Lemma~2]{Faj71}.
The identities of $(S_6,\uast)$---a rectangular band with \inv---is long known to be axiomatized by $\{ x^2 \approx x,\, xyz \approx xz \}$ \cite[Lemma~2]{Faj71}.
Consequently, every {\insem} of order four is {\fb}.

\section{The {\insem} $(S_2,\uast)$} \label{sec: S2}

For any word $\bfw \in \Fsem$, the \textit{interior} of~$\bfw$, denoted by $\sfint(\bfw)$, is the word obtained from~$\bfw$ by removing its first and last variables.
Specifically, if $\bfw = h \bfw_0 t$ for some $h,t \in \cXX$ and $\bfw_0 \in \Fmon$, then $\sfint(\bfw) = \bfw_0$.
Note that if $|\bfw| \leq 2$, then $\sfint(\bfw) = 1$.

\begin{proposition}
The identities 
\begin{subequations} \label{id: S2 basis}
\begin{gather}
x^3 \approx x^2, \quad xyx \approx x^2y^2, \quad xyx \approx y^2x^2, \quad xy^2z \approx xyz, \label{id: S2 plain} \\
x^*x \approx x^2, \quad xx^* \approx x^2, \quad x^*yx \approx xyx, \quad xyx^* \approx xyx, \quad xy^*z \approx xyz \label{id: S2 nonplain}
\end{gather}
\end{subequations}
constitute an identity basis for the {\insem} $(S_2,\uast)$.
\end{proposition}

\begin{proof}
In this proof, a word~$\bfw$ is in \textit{canonical form} if every variable in $\sfint(\bfw)$ and every non-simple variable in~$\bfw$ are plain.
It is easy to see that the identities~\eqref{id: S2 nonplain} can be used to convert any word into one in canonical form.

It is routine to check that $(S_2,\uast)$ satisfies the identities~\eqref{id: S2 basis}.
Hence, there exists some set~$\Sigma$ of identities of $(S_2,\uast)$, formed by words in canonical form, such that $\{ \eqref{id: S2 basis} \} \cup \Sigma$ is an identity basis for $(S_2,\uast)$.
Now the identities~\eqref{id: S2 plain} in fact constitute an identity basis for the semigroup~$S_2$ \cite[Variety~24]{AACLR23}, so every plain identity of $(S_2,\uast)$ is deducible from~\eqref{id: S2 plain}.
Therefore, the identities in~$\Sigma$ can be assumed non-plain.
Let $\bfu \approx \bfv$ be any identity in~$\Sigma$, say with~$\bfu$ non-plain.
If $|\bfu| = 1$, then $\bfu = x^*$ for some $x \in \cX$, whence it is easy to show that the identity $\bfu \approx \bfv$ is trivial and so is clearly deducible from~\eqref{id: S2 basis}.

It remains to consider the case $|\bfu| \geq 2$.
Since~$\bfu$ is in canonical form and contains a non-plain variable, say~$x^*$ with $x \in \cX$, the variable~$x^*$ is simple and so is either the head or the tail of~$\bfu$.
In view of Lemma~\ref{Lem: reversal}, it suffices to assume the former, so that $\bfu = x^* \sfint(\bfu)\sft(\bfu)$ with $\sfint(\bfu) \in \cX^+ \cup \{ 1 \}$ and $\sft(\bfu) \in \cXX$ such that $x \notin \sfcon(\sfint(\bfu)\overline{\sft(\bfu)})$.
If $\sfh(\bfv) \neq x^*$, then letting $\varphi_1: \cX \to S_2$ be the substitution that maps~$x$ to~$2$ and every other variable to~$4$, the contradiction $\varphi_1(\bfu) = 3 \neq \varphi_1(\bfv)$ is obtained.
Therefore, $\sfh(\bfv) = x^*$, so that $\bfv = x^* \sfint(\bfv)\sft(\bfv)$ with $\sfint(\bfv) \in \cX^+ \cup \{ 1 \}$ and $\sft(\bfv) \in \cXX$ such that $x \notin \sfcon(\sfint(\bfv)\overline{\sft(\bfv)})$.
There are two cases.

\medskip

\noindent{\sc Case~1:} $\sft(\bfu),\sft(\bfv) \in \cX$.
Then $\sfint(\bfu)\sft(\bfu)$ and $\sfint(\bfv)\sft(\bfv)$ are plain, so that $\overline{\bfu} = x \sfint(\bfu)\sft(\bfu)$ and $\overline{\bfv} = x \sfint(\bfv)\sft(\bfv)$.

\medskip

\noindent{\sc Case~2:} $\sft(\bfu) \notin \cX$ or $\sft(\bfv) \notin \cX$.
By symmetry, suppose that $\sft(\bfu) \in \cX^*$, say $\sft(\bfu) = y^*$ for some $y \in \cX \backslash \{x\}$.
Then $\bfu = x^* \sfint(\bfu) y^*$ with $\sfint(\bfu) \in \cX^+ \cup \{ 1 \}$ such that $x,y \notin \sfcon(\sfint(\bfu))$.
If $\sft(\bfv) \neq y^*$, then letting $\varphi_2: \cX \to S_2$ be the substitution that maps~$y$ to~$3$ and every other variable to~$4$, the contradiction $\varphi_2(\bfu) = 2 \neq \varphi_2(\bfv)$ is obtained.
Therefore, $\sft(\bfv) = y^*$, so that $\bfv = x^* \sfint(\bfv) y^*$ with $\sfint(\bfv) \in \cX^+ \cup \{ 1 \}$ such that $x,y \notin \sfcon(\sfint(\bfv))$.
It follows that $\overline{\bfu} = x \sfint(\bfu) y$ and $\overline{\bfv} = x \sfint(\bfv) y$.

\medskip

It is clear that in both cases, the identities $\{ \eqref{id: inv},\, \bfu \approx \bfv \}$ and $\{ \eqref{id: inv},\, \overline{\bfu} \approx \overline{\bfv} \}$ are deducible from one another.
Since $\overline{\bfu} \approx \overline{\bfv}$ is an identity of $S_2$, it is deducible from~\eqref{id: S2 plain}.
It follows that $\bfu \approx \bfv$ is deducible from $\{ \eqref{id: inv}, \eqref{id: S2 basis} \}$.
Consequently, every identity in~$\Sigma$ is deducible from $\{ \eqref{id: inv}, \eqref{id: S2 basis} \}$, so that~\eqref{id: S2 basis} is an identity basis for $(S_2,\uast)$.
\end{proof}

\section{The {\insem} $(S_5,\uast)$} \label{sec: A0}

The {\insem} $(S_5,\uast)$ is isomorphic to the semigroup \[ A_0 = \langle \rme,\rmf \,|\, \rme^2=\rme, \, \rmf^2 = \rmf, \, \rmf\rme = 0 \rangle = \{ 0, \rme,\rmf,\rmef \} \] with the operation~$\uast$ that interchanges $\rme$ and $\rmf$ and fixes every other element.
\[
\begin{tabular}{c | c c c c}
$A_0$   & 0 & $\rmef$ & $\rme$ & $\rmf$ \\ \hline
0       & 0 & 0       & 0      & 0 \\
$\rmef$ & 0 & 0       & 0      & $\rmef$ \\
$\rme$  & 0 & $\rmef$ & $\rme$ & $\rmef$ \\
$\rmf$  & 0 & 0       & 0      & $\rmf$ \\ \hline\hline
$x$     & 0 & $\rmef$ & $\rme$ & $\rmf$ \\ \hline
$x^*$   & 0 & $\rmef$ & $\rmf$ & $\rme$
\end{tabular}
\]
The {\insem} $(A_0,\uast)$ is isomorphic to the {\inv} subsemigroup of $(A_0^1,{}^\sk)$ that consists of its non-unit elements.

The semigroup $A_0$ has been known to be {\fb} for over 40 years~\cite{Edm80}.
The finite basis problem for the {\insem} $(A_0,\uast)$ has not been considered and has only recently been questioned \cite[Question~6.4]{Lee17b}.
The present section addresses this problem by showing that $(A_0,\uast)$ is \fb.

\begin{proposition} \label{Prop: A0 basis}
The identities
\begin{subequations} \label{id: A0 basis}
\begin{align}
xx^*xy \approx xx^*x, \quad yxx^*x & \approx xx^*x, \quad xx^*x \approx yy^*y, \label{id: A0 xx*xy=xx*x} \\
xyx^* & \approx xy^*x^*, \label{id: A0 xyx*=xy*x*} \\
x^2 \sfH x^* \approx x \sfH x^*, & \quad x \sfH (x^*)^2 \approx x \sfH x^*, \label{id: A0 xxHx*=xHx*} \\
xy \sfH y^* & \approx yx \sfH x^*, \label{id: A0 xyHy*=yxHx*} \\
x \sfH x^* y & \approx y^* x \sfH x^*, \label{id: A0 xHx*y=y*xHx*} \\
x^2 \sfH y\sfT y^* & \approx x \sfH y\sfT y^*, \label{id: A0 xxHyTy*=xHyTy*} \\
xy \sfH z\sfT z^* & \approx yx \sfH z\sfT z^*, \label{id: A0 xyHzTz*=yxHzTz*} \\
x^3 \approx x^2, \quad x^2yx & \approx xyx, \quad xyx^2 \approx xyx, \label{id: A0 xxx=xx} \\
xyx & \approx yxy, \label{id: A0 xyx=yxy} \\
x \sfH yz \sfT x & \approx x \sfH zy \sfT x, \label{id: A0 xHyzTx=xHzyTx}
\end{align}
\end{subequations}
where $\sfH \in \{ 1,h \}$ and $\sfT \in \{ 1,t \}$\up, constitute an identity basis for $(A_0,\uast)$.
\end{proposition}

It is easily checked that $(A_0,\uast)$ satisfies the identities~\eqref{id: A0 basis}.
In Section~\ref{subsec: A0 identities}, some information on identities of $(A_0,\uast)$ are given.
In Section~\ref{subsec: A0 forms}, it is shown that the identities of $(A_0,\uast)$ can be used to convert every mixed word into one of two specific forms.
Based on these results, it is shown in Section~\ref{subsec: A0 proof} that every identity of $(A_0,\uast)$ is deducible from $\{ \eqref{id: inv}, \eqref{id: A0 basis} \}$.
This completes the proof of Proposition~\ref{Prop: A0 basis}.

\begin{corollary} \label{Cor: A0 simplified}
The identities
\begin{equation} \label{id: A0 simplified}
\begin{split}
x^3 & \approx x^2, \quad xyxy \approx xyx, \quad x^2x^* \approx xx^*, \quad x^2yx^* \approx xyx^*, \\
& xy^*x^* \approx xyx^*, \quad xx^*x \approx yy^*y, \quad xyx^*z \approx z^*xyx^*
\end{split}
\end{equation}
constitute an identity basis for $(A_0, \uast)$.
\end{corollary}

\begin{proof}
It is routine to check, say with Prover9~\cite{McC05}, that the identities $\{ \eqref{id: inv}, \eqref{id: A0 basis} \}$ and $\{ \eqref{id: inv}, \eqref{id: A0 simplified} \}$ are deducible from one another.
\end{proof}

\begin{remark}
\begin{enumerate}[\ \rm(i)]
\item Not only is the semigroup~$A_0$ \fb, it is \textit{hereditarily \fb} in the sense that every semigroup in the variety $\Var\,A_0$ is {\fb} \cite[Corollary~4.3]{Lee04}.
\item In contrast, the {\fb} {\insem} $(A_0,\uast)$ is not hereditarily {\fb} because the variety $\Var(A_0,\uast)$ contains a {\nfb} {\insem} \cite[Proposition~3.8]{Lee22}.
\end{enumerate}
\end{remark}

\subsection{Some identities of $(A_0,\uast)$} \label{subsec: A0 identities}

\begin{lemma} \label{Lem: id: A0 plain}
The identities $\{ \eqref{id: A0 xxx=xx}, \eqref{id: A0 xyx=yxy} \}$ constitute an identity basis for the semigroup $A_0$.
\end{lemma}

\begin{proof}
The identities $\Sigma = \{ x^3 \approx x^2,\, xyxy \approx xyx,\, xyxy \approx yxy \}$ constitute an identity basis for~$A_0$ \cite[Theorem~4.1]{LV07}; it is routine to verify that $\Sigma$ and $\{ \eqref{id: A0 xxx=xx}, \eqref{id: A0 xyx=yxy} \}$ are deducible from one another.
\end{proof}

\begin{lemma} \label{Lem: A0 bipartite}
Let $\bfu \approx \bfv$ be any identity of $(A_0,\uast)$ such that either~$\bfu$ or~$\bfv$ is bipartite.
Then $\bfu \approx \bfv$ is deducible from $\{ \eqref{id: inv}, \eqref{id: A0 basis} \}$.
\end{lemma}

\begin{proof}
Since $(S\ell_3,{}^\sk)$ is isomorphic to $(A_0, \uast)$ modulo the ideal $\{ 0,\rmef \}$, the identity $\bfu \approx \bfv$ is satisfied by $(S\ell_3,{}^\sk)$.
Then since either~$\bfu$ or~$\bfv$ is bipartite, by Lemma~\ref{Lem: A0 twisted}, both~$\bfu$ and~$\bfv$ are bipartite with $\sfcon({\bfu})=\sfcon({\bfv})$.
It follows from Lemma~\ref{Lem: A0 mixed bipartite} that $(A_0, \uast)$ satisfies the plain identity $\overline{\bfu} \approx \overline{\bfv}$.
By Lemma~\ref{Lem: id: A0 plain}, the identities $\{ \eqref{id: A0 xxx=xx}, \eqref{id: A0 xyx=yxy} \}$ constitute an identity basis for~$A_0$, so that $\overline{\bfu} \approx \overline{\bfv}$ is deducible from $\{ \eqref{id: A0 xxx=xx}, \eqref{id: A0 xyx=yxy} \}$.
It then follows from Lemma~\ref{Lem: A0 mixed bipartite} that $\bfu \approx \bfv$ is deducible from $\{ \eqref{id: inv}, \eqref{id: A0 basis} \}$.
\end{proof}

An \textit{\Ablock} is a word of the form \[ \bfc = (y_1 y_2 \cdots y_k)^2, \] where $y_1,y_2,\ldots,y_k \in \cXX$ are such that $\overline{y_1} \order \overline{y_2} \order \cdots \order \overline{y_k}$ in~$\cX$ and $k \geq 1$.
Note that every {\Ablock} is bipartite and connected.

\begin{lemma}[Lee {\cite[Lemma~2.1]{Lee23a}}] \label{Lem: A0 plain id}
Let $\bfu,\bfv \in \cX^+$ be any plain connected words such that $\sfcon(\bfu) = \sfcon(\bfv)$.
Then $\bfu \approx \bfv$ is an identity of the semigroup~$A_0$.
\end{lemma}

\begin{lemma}
\label{Lem: A0 block}
Let $\bfw \in \Fsem$ be any bipartite connected word such that $\sfcon(\bfw) = \{ y_1,y_2,\ldots,y_k\}$ and $\overline{y_1} \order \overline{y_2} \order \cdots \order \overline{y_k}$ in~$\cX$.
Then the identities $\{\eqref{id: A0 xxx=xx},\eqref{id: A0 xyx=yxy} \}$ can be used to convert~$\bfw$ into the {\Ablock} $\bfc = (y_1 y_2 \cdots y_k)^2$.
\end{lemma}

\begin{proof}
Since~$\bfw$ is bipartite, there exists a partition $I \cup J = \{ 1,2, \ldots,k\}$ such that $y_i \in \cX$ for all $i \in I$ and $y_j \in \cX^*$ for all $j \in J$.
Let $\varphi : \cX \to \cXX$ denote the substitution \[ \varphi(t) = \begin{cases} t^* & \text{if $t \in \{ y_j,y_j^* \,|\, j \in J \}$}, \\ t & \text{otherwise}. \end{cases} \]
Then it is easy to check that for any word $\bfv \in \sfcon(\bfw)^+$, we have $\varphi(\bfv) = \overline{\bfv}$ and $\varphi(\overline{\bfv}) = \bfv$.
Since $\overline{\bfw}$ and $\overline{\bfc}$ are plain words such that $\sfcon(\overline{\bfw}) = \sfcon(\overline{\bfc})$, by Lemma~\ref{Lem: A0 plain id}, the identity $\overline{\bfw} \approx \overline{\bfc}$ is satisfied by~$A_0$.
Then by Lemma~\ref{Lem: id: A0 plain}, $\overline{\bfw} \approx \overline{\bfc}$ is deducible from $\{ \eqref{id: A0 xxx=xx},\eqref{id: A0 xyx=yxy} \}$.
Consequently, \[ \bfw = \varphi(\overline{\bfw}) \stackrel{\eqref{id: A0 xxx=xx},\eqref{id: A0 xyx=yxy}}{\approx} \varphi(\overline{\bfc}) = \bfc. \qedhere \]
\end{proof}

\subsection{Some special forms of words} \label{subsec: A0 forms}

It is easily checked that for any substitution $\varphi: \cX \to A_0$ and any variable $z \in \cX$, we have $\varphi(zz^*z) = 0$ in $A_0$.
Therefore, in the $\Var(A_0,\uast)$-free algebra over~$\cX$, the class $[zz^*z]$ containing $zz^*z$ is its zero element.
This phenomenon can also be seen from the identities~\eqref{id: A0 xx*xy=xx*x} of $(A_0,\uast)$.

Words of other possible forms in the class $[zz^*z]$ are given in the following result.

\begin{lemma} \label{Lem: A0 zz*z}
Let $\bfw \in \Fsem$.
Suppose that either $xx^*x \hookrightarrow \bfw$ for some $x \in \cXX$ or $xx^*yy^* \hookrightarrow \bfw$ for some $x,y \in \cXX$.
Then the identities~\eqref{id: A0 basis} can be used to convert~$\bfw$ into the word $zz^*z$ for any $z \in \cXX$.
\end{lemma}

\begin{proof}
If $\bfw = \bfw_0 x \bfw_1 x^* \bfw_2 x \bfw_3$ for some $\bfw_0, \bfw_1, \bfw_2, \bfw_3 \in \Fmon$, then \[ \bfw \stackrel{\eqref{id: A0 xHyzTx=xHzyTx}}{\approx} \bfw_0 x \bfw_1 \bfw_2 x^*x \bfw_3 \stackrel{\eqref{id: A0 xyHzTz*=yxHzTz*}}{\approx} \bfw_0 \bfw_1 \bfw_2 xx^*x \bfw_3 \stackrel{\eqref{id: A0 xx*xy=xx*x}}{\approx} zz^*z. \]
If $\bfw = \bfw_0 x \bfw_1 x^* \bfw_2 y \bfw_3 y^* \bfw_4$ for some $\bfw_0, \bfw_1, \bfw_2, \bfw_3, \bfw_4 \in \Fmon$, then \[ \bfw \stackrel{\eqref{id: A0 xxHx*=xHx*}}{\approx} \bfw_0 x^2 \bfw_1 x^* \bfw_2 y \bfw_3 y^* \bfw_4 \stackrel{\eqref{id: A0 xyHzTz*=yxHzTz*}}{\approx} \bfw_0 xx^*x \bfw_1 \bfw_2 y \bfw_3 y^* \bfw_4 \stackrel{\eqref{id: A0 xx*xy=xx*x}}{\approx} zz^*z. \qedhere \]
\end{proof}

A word $\bfw \in \Fsem$ is in \textit{$A_0$-standard form} if
\begin{equation}
\bfw=\bfw_1 x \bfw_2 x^*, \label{Dis: A0-standard}
\end{equation}
where $x \in \cXX$, $\bfw_1 = x_1x_2\cdots x_m$, and $\bfw_2 = \bfs_0 \prod_{i=1}^p (\bfc_i\bfs_i)$  for some $m,p \geq 0$ such that the following conditions hold:
\begin{enumerate}[({\Asf}1)]
\item $x_1,x_2,\ldots, x_m \in \cXX$ are such that $\overline{x_1} \order \overline{x_2} \order \cdots \order \overline{x_m} \order \overline{x}$;

\item $\bfs_0, \bfs_1, \ldots, \bfs_p \in \Fmon$ are simple and $\bfc_1, \bfc_2, \ldots, \bfc_p \in \Fsem$ are {\Ablock}s;

\item $x_1, x_2, \ldots, x_m, x, \bfs_0, \bfs_1, \ldots, \bfs_p, \bfc_1, \bfc_2, \ldots, \bfc_p$ are pairwise disjoint;

\item if $\bfw_2 \neq 1$, then
\begin{enumerate}[\ \rm(a)]
\item $p = 0$ with $\bfw_2 = \bfs_0$ and $\bfs_0 \in \cX$; or \label{Asf4 singleton}
\item $p=1$ with $\bfw_2 = \bfc_1$, $\bfs_0 = \bfs_1 = 1$, and $\sfh(\bfc_1) \in \cX$; or \label{Asf4 s0 orderedsquare}
\item $\overline{\sfh(\bfs_0)} \order \overline{\sft(\bfw_2)}$ when $\bfs_0 \neq 1$; or \label{Asf4 s0 nonempty}
\item $\overline{\sft(\bfc_1)} \order \overline{\sft(\bfw_2)}$ when $\bfs_0 = 1$. \label{Asf4 s0 empty}
\end{enumerate}
\end{enumerate}

\begin{remark} \label{Rmk: A0-standard}
The following holds for the word~$\bfw$ in~\eqref{Dis: A0-standard} in $A_0$-standard form:
\begin{enumerate}[\ \rm(i)]
\item If $m=0$, then $\bfw_1 = 1$.
\item If $p = 0$, then $\bfw_2 = \bfs_0$.
\item $\{x, x^*\}$ is the only mixed pair of~$\bfw$ and $x,x^* \notin \sfcon(\bfw_1\bfw_2)$.
\item $\bfw_1$ and $\bfw_2$ are bipartite words such that $\sfcon(\overline{\bfw_1}) \cap \sfcon(\overline{\bfw_2}) = \emptyset$.
\item Each variable in $\cX$ occurs at most twice in $\overline{\bfw}$.
\item If $\bfw_2 \neq 1$ and $|\mathsf{con}(\overline{\bfw_2})| \geq 2$, then $\overline{\sfh(\bfw_2)} \order \overline{\sft(\bfw_2)}$ by condition (\Asf4).
This is because
\begin{itemize}
\item if $\bfs_0 \neq 1$, then $\overline{\sfh(\bfw_2)} = \overline{\sfh(\bfs_0)} \order \overline{\sft(\bfw_2)}$ by condition (\Asf\ref{Asf4 s0 nonempty});
\item if $\bfs_0 = 1$, then since $\overline{\sfh(\bfw_2)} = \overline{\sfh(\bfc_1)} \ordereq \overline{\sft(\bfc_1)}$ due to $\bfc_1$ being an {\Ablock}, we have $\overline{\sfh(\bfw_2)} \ordereq \overline{\sft(\bfc_1)} \order \overline{\sft(\bfw_2)}$ by condition (\Asf\ref{Asf4 s0 empty}).
\end{itemize}
\end{enumerate}
\end{remark}

\begin{lemma} \label{Lem: A0 substitution}
Let $\bfw=\bfw_1 x \bfw_2 x^*$ be the word in~\eqref{Dis: A0-standard} in $A_0$-standard form.
Then there exist substitutions $\alpha_\bfw, \beta_\bfw : \cX \to A_0$ such that
\begin{enumerate}[\ \rm(i)]
\item $\alpha_\bfw(\bfw_1 x \bfw_2) = \rme$ and $\alpha_\bfw(x^*) = \rmf$\up, so that $\alpha_\bfw(\bfw) = \rmef$\up;
\item $\beta_\bfw(\bfw_1 x) = \rme$ and $\beta_\bfw(\bfw_2x^*) = \rmf$\up, so that $\beta_\bfw(\bfw) = \rmef$\up;
\item $\alpha_\bfw(t) = \beta_\bfw(t) = 0$ for all $t \in \cX$ such that $t \notin \sfcon(\overline{\bfw})$.
\end{enumerate}
\end{lemma}

\begin{proof}
It follows from Remark~\ref{Rmk: A0-standard}(iii),(iv) that $\sfcon(\bfw_1) = \mathcal{H}_1 \cup \mathcal{K}_1^*$ and $\sfcon(\bfw_2) = \mathcal{H}_2 \cup \mathcal{K}_2^*$ for some $\mathcal{H}_1,\mathcal{H}_2,\mathcal{K}_1,\mathcal{K}_2 \subseteq \cX$ such that $\mathcal{H}_1, \mathcal{H}_2, \mathcal{K}_1, \mathcal{K}_2, \{ x,x^* \}$ are pairwise disjoint sets.
By symmetry, it suffices to assume that $x \in \cX$, so that $\sfcon(\overline{\bfw}) = \mathcal{H}_1 \cup \mathcal{K}_1 \cup \mathcal{H}_2 \cup \mathcal{K}_2 \cup \{x\}$.
Define
\[
\alpha_\bfw(t) = \begin{cases}
\rme & \text{if $t \in \mathcal{H}_1 \cup \mathcal{H}_2 \cup \{ x \}$}, \\
\rmf & \text{if $t \in \mathcal{K}_1 \cup \mathcal{K}_2$}, \\
0 & \text{otherwise};
\end{cases}
\quad
\beta_\bfw(t) = \begin{cases}
\rme & \text{if $t \in \mathcal{H}_1 \cup \mathcal{K}_2 \cup \{ x \}$}, \\
\rmf & \text{if $t \in \mathcal{K}_1 \cup \mathcal{H}_2$}, \\
0 & \text{otherwise}.
\end{cases}
\]
Then it is routinely checked that the substitutions $\alpha_\bfw$ and $\beta_\bfw$ satisfy (i)--(iii).
\end{proof}

\begin{corollary} \label{C: A0 violate}
For any word~$\bfw$ in $A_0$-standard form and any $z \in \cXX$\up, the identity $\bfw \approx zz^*z$ is not satisfied by $(A_0,\uast)$.
\end{corollary}

\begin{proof}
Under the substitution $\alpha_\bfw : \cX \to A_0$ in Lemma~\ref{Lem: A0 substitution}, we have $\alpha_\bfw(\bfw) = \rmef$ and $\alpha_\bfw(zz^*z) = 0$.
\end{proof}

\begin{lemma} \label{Lem: A0 standard form}
Let $\bfw$ be any mixed word.
Then the identities $\{ \eqref{id: inv}, \eqref{id: A0 basis} \}$ can be used to convert~$\bfw$ into exactly one of the following\up:
\begin{enumerate}[\ \rm(i)]
\item the word $zz^*z$ for any $z \in \cXX$\up;
\item some word in $A_0$-standard form.
\end{enumerate}
\end{lemma}

\begin{proof}
In view of the identities $\{ \eqref{id: inv}, \eqref{id: A0 xxHx*=xHx*}, \eqref{id: A0 xHx*y=y*xHx*}, \eqref{id: A0 xyHzTz*=yxHzTz*} \}$, we may assume that $\bfw = \bfw_1 x \bfw_2 x^*$, where $x \in \cXX$ and $\bfw_1,\bfw_2 \in \Fmon$ are such that $x \notin \sfcon(\bfw_1)$, $x,x^* \notin \sfcon(\bfw_2)$, and $\bfw_2$ is bipartite.
If either $x^* \in \sfcon(\bfw_1)$ or~$\bfw_1$ contains some mixed pair, then by Lemma~\ref{Lem: A0 zz*z}, the identities~\eqref{id: A0 basis} can be used to convert~$\bfw$ into the word $zz^*z$  for any $z \in \cXX$.
Therefore, suppose that $x^* \notin \sfcon(\bfw_1)$ and~$\bfw_1$ is bipartite.
In summary, we may assume that
\begin{enumerate}[\ (a)]
\item $x,x^* \notin \sfcon(\bfw_1\bfw_2)$ and
\item $\bfw_1$ and $\bfw_2$ are bipartite.
\end{enumerate}

Suppose that $\sfcon(\bfw_1) \cap \sfcon(\bfw_2) \neq \emptyset$.
Then the~$x$ in~$\bfw$ is sandwiched between two occurrences of the same variable, one occurring in~$\bfw_1$ and one in~$\bfw_2$.
Therefore, $\bfw_1 = \bfa y \bfb$ and $\bfw_2 = \bfe y \bff$ for some $y \in \cXX$ and $\bfa,\bfb,\bfe,\bff \in \Fmon$ such that $x \notin \sfcon(\bfa\bfb\bfe\bff)$ and $y \notin \sfcon(\bfb\bff)$, whence
\begin{align*}
\bfw & \makebox[0.36in]{$=$} \bfa y \bfb x \bfe y \bff x^* \stackrel{\eqref{id: A0 xxx=xx}}{\approx} \bfa y \bfb x \bfe y^2 \bff x^* \stackrel{\eqref{id: A0 xHyzTx=xHzyTx}}{\approx} \bfa y \bfb \bfe yxy \bff x^* \\
& \makebox[0.36in]{$\stackrel{\eqref{id: A0 xyx=yxy}}{\approx}$} \bfa y \bfb \bfe xyx \bff x^* \stackrel{\eqref{id: A0 xyHzTz*=yxHzTz*}}{\approx} \bfa y \bfb \bfe yx^2 \bff x^* \stackrel{\eqref{id: A0 xxHx*=xHx*}}{\approx} \bfa y \bfb \bfe yx \bff x^*;
\end{align*}
in other words, $x$ can be moved to the right until it is no longer sandwiched by any two occurrences of~$y$.
This process can be repeated until~$x$ is not sandwiched by any two occurrences of the same variable.
Therefore, we may further assume that
\begin{enumerate}[\ (a)]
\item[(c)] $\sfcon(\bfw_1) \cap \sfcon(\bfw_2) = \emptyset$.
\end{enumerate}

Suppose that $\sfcon(\overline{\bfw_1}) \cap \sfcon(\overline{\bfw_2}) \neq \emptyset$.
Then the~$x$ in~$\bfw$ is sandwiched by some mixed pair $\{ y,y^*\}$ with $y \in \sfcon(\bfw_1)$ and $y^* \in \sfcon(\bfw_2)$.
Therefore, $\bfw_1 = \bfa y \bfb$ and $\bfw_2 = \bfe y^* \bff$ for some $\bfa,\bfb,\bfe,\bff \in \Fmon$ such that $x \notin \sfcon(\bfa\bfb\bfe\bff)$ and $y^* \notin \sfcon(\bfe)$.
By~(c), we also have $y \notin \sfcon(\bfe)$.
Then
\begin{align*}
\bfw & \makebox[0.37in]{$=$} \bfa y \bfb x \bfe y^* \bff x^* \stackrel{\eqref{id: A0 xyx*=xy*x*}}{\approx} \bfa y \bfb x (\bfe y^* \bff)^* x^* \stackrel{\eqref{id: inv}}{\approx} \bfa y \bfb x \bff^* y \bfe^* x^* \\
& \makebox[0.37in]{$\stackrel{\eqref{id: A0 xxx=xx}}{\approx}$} \bfa y \bfb x \bff^* y^2 \bfe^* x^* \stackrel{\eqref{id: A0 xHyzTx=xHzyTx}}{\approx} \bfa y \bfb \bff^* yxy \bfe^* x^* \stackrel{\eqref{id: A0 xyx=yxy}}{\approx} \bfa y \bfb \bff^* xyx \bfe^* x^* \\
& \makebox[0.37in]{$\stackrel{\eqref{id: A0 xyHzTz*=yxHzTz*}}{\approx}$} \bfa y \bfb \bff^* yx^2 \bfe^* x^* \stackrel{\eqref{id: A0 xxHx*=xHx*}}{\approx} \bfa y \bfb \bff^* yx \bfe^* x^* \stackrel{\eqref{id: A0 xyx*=xy*x*}}{\approx} \bfa y \bfb \bff^* yx \bfe x^*.
\end{align*}
Since $y,y^* \notin \sfcon(\bfe)$, the variable~$x$ is no longer sandwiched by the mixed pair $\{ y,y^*\}$.
This process can be repeated until~$x$ is not sandwiched by any mixed pair.
Therefore, we may further assume that
\begin{enumerate}[\ (a)]
\item[(d)] $\sfcon(\overline{\bfw_1}) \cap \sfcon(\overline{\bfw_2}) = \emptyset$.
\end{enumerate}

Since the prefix $\bfw_1$ of~$\bfw$ is bipartite by~(b), the identities~\eqref{id: A0 xyHzTz*=yxHzTz*} can be used to put the variables in $\bfw_1$ in order, and the identities~\eqref{id: A0 xxHyTy*=xHyTy*} can be used to reduce the exponent of any non-simple variable to~$1$.
Hence, we may assume that $\bfw_1 = x_1x_2 \cdots x_m$ for some $x_1,x_2,\ldots,x_m \in \cXX$ with $m \geq 0$ such that $\overline{x_1} \order \overline{x_2} \order \cdots \order \overline{x_m}$.
By~(a), we have $\overline{x} \notin \{ \overline{x_1},\overline{x_2},\ldots,\overline{x_m}\}$.
If $\overline{x_m} \norder \overline{x}$, then \[ \bfw = x_1x_2 \cdots x_{m-1} x_m x \bfw_2 x^* \stackrel{\eqref{id: A0 xyHy*=yxHx*}}{\approx} x_1x_2 \cdots x_{n-1} x x_m \bfw_2 x_m^*, \] and the identities~\eqref{id: A0 xyHzTz*=yxHzTz*} can be used to put the variables in the prefix $x_1x_2 \cdots x_{m-1} x$ in an order such that condition (\Asf1) is satisfied.

Since $\bfw_2$ is bipartite by~(b), it can be written in the form $\bfw_2 = \bfs_0 \prod_{i=1}^p (\bfc_i\bfs_i)$, where $\bfs_0, \bfs_1, \ldots, \bfs_p \in \Fmon$ are simple words and $\bfc_1, \bfc_2, \ldots, \bfc_p \in \Fsem$ are connected words with $p \geq 0$ such that
\begin{enumerate}[\ (a)]
\item[(e)] $\bfs_0, \bfs_1, \ldots, \bfs_p, \bfc_1, \bfc_2, \ldots, \bfc_p$ are pairwise disjoint.
\end{enumerate}
Then by Lemma~\ref{Lem: A0 block}, the identities $\{\eqref{id: A0 xxx=xx},\eqref{id: A0 xyx=yxy} \}$ can be used to convert each~$\bfc_i$ into the {\Ablock} word $\widehat{\,\bfc_i\,}$ with $\sfcon(\bfc_i) = \sfcon(\widehat{\,\bfc_i\,})$.
Therefore, condition (\Asf2) is satisfied, and it follows from~(d) and~(e) that condition (\Asf3) is also satisfied.

It remains to address condition (\Asf4).
Suppose that $\bfw_2 = \bfs_0 \prod_{i=1}^p (\bfc_i\bfs_i) \neq 1$.
There are five cases.

\medskip

\noindent{\sc Case~1:}
$p=0$.
Then $\bfw_2 = \bfs_0$.

\smallskip

\noindent{\sc Subcase~1.1:} $|\mathsf{con}(\overline{\bfs_0})|=1$.
Then $\bfw_2 =\bfs_0 \in \{y,y^*\}$ for some $y \in \cX$.
If $\bfw_2=y$, then condition (\Asf\ref{Asf4 singleton}) is satisfied.
If $\bfw_2=y^*$, then \[ \bfw = \bfw_1 x y^* x^* \stackrel{\eqref{id: A0 xyx*=xy*x*}}{\approx} \bfw_1 x (y^*)^* x^* \stackrel{\eqref{id: inv}}{\approx} \bfw_1 xyx^*; \] in other words, the identities $\{\eqref{id: inv},\eqref{id: A0 xyx*=xy*x*}\}$ can be used to convert the $y^*$ in $\bfw$ to $y \in \cX$, so that condition (\Asf\ref{Asf4 singleton}) is satisfied.

\smallskip

\noindent{\sc Subcase~1.2:} $|\mathsf{con}(\overline{\bfs_0})| \geq 2$.
Then $\bfw_2 = \bfs_0 = y_1y_2 \cdots y_r$ for some distinct $y_1,y_2,\ldots,y_r \in \cXX$ with $r \geq 2$.
If $\overline{\sfh(\bfs_0)} \order \overline{\sft(\bfw_2)}$, then condition (\Asf\ref{Asf4 s0 nonempty}) is satisfied.
If $\overline{\sfh(\bfs_0)} \norder \overline{\sft(\bfw_2)}$, so that $\overline{y_1} = \overline{\sfh(\bfs_0)} \rorder \overline{\sft(\bfs_0)} = \overline{y_r}$, then \[ \bfw \stackrel{\eqref{id: A0 xyx*=xy*x*}}{\approx} \bfw_1 x \bfw_2^*x^* = \bfw_1 x \underbrace{y_r^*y_{r-1}^* \cdots y_1^*}_{\bfs_0^* \, = \, \bfw_2^*} x^*, \] where  $\overline{\sfh(\bfs_0^*)} = \overline{y_r} \order \overline{y_1} = \overline{\sft(\bfs_0^*)}  = \overline{\sft(\bfw_2^*)}$, whence condition (\Asf\ref{Asf4 s0 nonempty}) is satisfied.

\medskip

\noindent{\sc Case~2:}
$p \geq 1$ and $\bfs_0 \neq 1 \neq \bfs_p$.
Then $\bfw_2 = \bfs_0 \prod_{i=1}^p (\bfc_i\bfs_i)$, where $\bfs_0 = y_1y_2 \cdots y_r$ and $\bfs_p = z_1z_2\cdots z_s$ for some distinct $y_1,y_2,\ldots,y_r,z_1,z_2,\ldots,z_s \in \cXX$ with $r,s \geq 1$.
If $\overline{\sfh(\bfs_0)} \order \overline{\sft(\bfw_2)}$, then condition (\Asf\ref{Asf4 s0 nonempty}) is satisfied.
Hence suppose that $\overline{\sfh(\bfs_0)} \norder \overline{\sft(\bfw_2)}$, so that $\overline{y_1} = \overline{\sfh(\bfs_0)} \rorder \overline{\sft(\bfs_p)} = \overline{z_s}$.
Then
\[ \bfw \stackrel{\eqref{id: A0 xyx*=xy*x*}}{\approx} \bfw_1 x\bfw_2^*x^* \stackrel{\eqref{id: inv}}{\approx} \bfw_1x \cdot \bfs_p^* \bfc_p^* \bfs_{p-1}^* \bfc_{p-1}^* \cdots \bfs_1^* \bfc_1^* \bfs_0^* \cdot x^*, \]
where $\overline{\sfh(\bfs_p^*)} = \overline{z_s} \order \overline{y_1} = \overline{\sft(\bfs_0^*)}$, and the identities $\{\eqref{id: inv},\eqref{id: A0 xxx=xx},\eqref{id: A0 xyx=yxy} \}$ can be used to convert each $\bfc_i^*$ into an {\Ablock} (see Lemma~\ref{Lem: A0 block}).
Therefore, condition (\Asf\ref{Asf4 s0 nonempty}) is satisfied.

\medskip

\noindent{\sc Case~3:}
$p \geq 1$ and $\bfs_0 \neq 1 = \bfs_p$.
Then $\bfw_2 = \bfs_0\prod_{i=1}^{p-1} (\bfc_i\bfs_i) \bfc_p$, where $\bfs_0 = z_1z_2\cdots z_s$ and $\bfc_p = (y_1y_2 \cdots y_r)^2$ for some distinct $y_1,y_2,\ldots,y_r,z_1,z_2,\ldots,z_s \in \cXX$ with $r,s \geq 1$ such that $\overline{y_1} \order \overline{y_2} \order \cdots \order \overline{y_r}$.
If $\overline{\sfh(\bfs_0)} \order \overline{\sft(\bfw_2)}$, then condition (\Asf\ref{Asf4 s0 nonempty}) is satisfied.
If $\overline{\sfh(\bfs_0)} \norder \overline{\sft(\bfw_2)}$, so that $\overline{z_1} = \overline{\sfh(\bfs_0)} \rorder \overline{\sft(\bfc_p)} = \overline{y_r}$, then
\[ \bfw \stackrel{\eqref{id: A0 xyx*=xy*x*}}{\approx} \bfw_1 x\bfw_2^*x^* \stackrel{\eqref{id: inv}}{\approx} \bfw_1x \cdot \bfc_p^* \bfs_{p-1}^* \bfc_{p-1}^*  \cdots \bfs_1^* \bfc_1^* \bfs_0^* \cdot x^*, \]
and the identities $\{\eqref{id: inv},\eqref{id: A0 xxx=xx},\eqref{id: A0 xyx=yxy} \}$ can be used to convert each $\bfc_i^*$ into an {\Ablock} (see Lemma~\ref{Lem: A0 block}); specifically, $\bfc_p^*$ is converted into $\widehat{\,\bfc_p^*\,} = (y_1^*y_2^* \cdots y_r^*)^2$.
Since $\overline{\sft(\bfc_p^*)} = \overline{y_r} \order \overline{z_1} =\overline{\sft({\bfs}_0^*)}$, condition (\Asf\ref{Asf4 s0 empty}) is satisfied.

\medskip

\noindent{\sc Case~4:}
$p \geq 1$ and $\bfs_0 = 1 \neq \bfs_p$.
Then $\bfw_2 = \prod_{i=1}^p (\bfc_i\bfs_i)$, where $\bfc_1 = (y_1y_2 \cdots y_r)^2$ and $\bfs_p = z_1z_2\cdots z_s$ for some distinct $y_1,y_2,\ldots,y_r,z_1,z_2,\ldots,z_s \in \cXX$ with $r,s \geq 1$ such that $\overline{y_1} \order \overline{y_2} \order \cdots \order \overline{y_r}$.
If $\overline{\sft(\bfc_1)} \order \overline{\sft(\bfw_2)}$, then condition (\Asf\ref{Asf4 s0 empty}) is satisfied.
If $\overline{\sft(\bfc_1)} \norder \overline{\sft(\bfw_2)}$, so that $\overline{y_r} = \overline{\sft(\bfc_1)} \rorder \overline{\sft(\bfs_p)} = \overline{z_s}$, then
\[ \bfw \stackrel{\eqref{id: A0 xyx*=xy*x*}}{\approx} \bfw_1 x\bfw_2^*x^* \stackrel{\eqref{id: inv}}{\approx} \bfw_1x \cdot \bfs_p^* \bfc_p^* \bfs_{p-1}^* \bfc_{p-1}^* \cdots \bfs_1^* \bfc_1^* \cdot x^*, \]
and the identities $\{\eqref{id: inv},\eqref{id: A0 xxx=xx},\eqref{id: A0 xyx=yxy} \}$ can be used to convert each $\bfc_i^*$ into an {\Ablock} (see Lemma~\ref{Lem: A0 block}); specifically, $\bfc_1^*$ is converted into $\widehat{\,\bfc_1^*\,} = (y_1^*y_2^* \cdots y_r^*)^2$.
Since $\overline{\sfh(\bfs_p^*)} = \overline{z_s} \order \overline{y_r} = \overline{\sft(\widehat{\,\bfc_1^*\,})}$, condition (\Asf\ref{Asf4 s0 nonempty}) is satisfied.

\medskip

\noindent{\sc Case~5:} $p \geq 1$ and $\bfs_0 = 1 =\bfs_p$.

\smallskip

\noindent{\sc Subcase~5.1:} $p =1$.
Then $\bfw_2 =\bfc_1=(y_1 y_2 \cdots y_k)^2$ for some $y_1,y_2, \ldots, y_k \in \cXX$ with $k \geq 1$ such that $\overline{y_1} \order \overline{y_2} \order \cdots \order \overline{y_k}$.
If $y_1 \in \cX$, then condition (\Asf\ref{Asf4 s0 orderedsquare}) is satisfied.
Hence suppose that $y_1 \in \cX^*$.
Then \[ \bfw \stackrel{\eqref{id: A0 xyx*=xy*x*}}{\approx} \bfw_1 x\bfw_2^*x^* = \bfw_1x \cdot \bfc_1^* \cdot x^*, \] and the identities $\{\eqref{id: inv},\eqref{id: A0 xxx=xx},\eqref{id: A0 xyx=yxy} \}$ can be used to convert $\bfc_1^*$ into an {\Ablock} $\widehat{\,\bfc_1^*\,} = (y_1^*y_2^* \cdots y_k^*)^2$ (see Lemma~\ref{Lem: A0 block}).
Now since $y_1^* \in \cX$, condition (\Asf\ref{Asf4 s0 orderedsquare}) is satisfied.

\smallskip

\noindent{\sc Subcase~5.2:} $p \geq 2$.
Then $\bfw_2 = \big(\prod_{i=1}^{p-1} (\bfc_i\bfs_i)\big)\bfc_p$, where $\bfc_1 = (y_1y_2 \cdots y_r)^2$ and $\bfc_p = (z_1z_2 \cdots z_s)^2$ for some distinct $y_1,y_2,\ldots,y_r,z_1,z_2,\ldots,z_s \in \cXX$ with $r,s \geq 1$ such that $\overline{y_1} \order \overline{y_2} \order \cdots \order \overline{y_r}$ and $\overline{z_1} \order \overline{z_2} \order \cdots \order \overline{z_s}$.
If $\overline{\sft(\bfc_1)} \order \overline{\sft(\bfw_2)}$, then condition (\Asf\ref{Asf4 s0 empty}) is satisfied.
Hence suppose that $\overline{\sft(\bfc_1)} \norder \overline{\sft(\bfw_2)}$, so that $\overline{z_s} = \overline{\sft(\bfc_p)} \rorder \overline{\sft(\bfc_1)} = \overline{y_r}$.
Then \[ \bfw \stackrel{\eqref{id: A0 xyx*=xy*x*}}{\approx} \bfw_1 x\bfw_2^*x^* \stackrel{\eqref{id: inv}}{\approx} \bfw_1x \cdot \bfc_p^* \bfs_{p-1}^* \bfc_{p-1}^* \cdots \bfs_1^* \bfc_1^* \cdot x^*, \] and the identities $\{\eqref{id: inv},\eqref{id: A0 xxx=xx},\eqref{id: A0 xyx=yxy} \}$ can be used to convert each $\bfc_i^*$ into an {\Ablock} (see Lemma~\ref{Lem: A0 block}); specifically, $\bfc_1^*$ is converted into $\widehat{\,\bfc_1^*\,} = (y_1^*y_2^* \cdots y_r^*)^2$ and~$\bfc_p^*$ is converted into $\widehat{\,\bfc_p^*\,} = (z_1^*z_2^* \cdots z_s^*)^2$.
Since $\overline{\sft(\widehat{\,\bfc_p^*\,})} = \overline{z_s} \order \overline{y_r} = \overline{\sft(\widehat{\,\bfc_1^*\,})}$, condition (\Asf\ref{Asf4 s0 empty}) is satisfied.

\medskip

Consequently, the identities $\{ \eqref{id: inv}, \eqref{id: A0 basis} \}$ can be used to convert~$\bfw$ into either $zz^*z$ or some word~$\widetilde{\bfw}$ in $A_0$-standard form.
But if the identities $\{ \eqref{id: inv}, \eqref{id: A0 basis} \}$ can be used to convert~$\bfw$ into both $zz^*z$ and~$\widetilde{\bfw}$, then that would imply that $(A_0,\uast)$ satisfies the identity $\widetilde{\bfw} \approx zz^*z$, which is impossible by Corollary~\ref{C: A0 violate}.
\end{proof}

\subsection{Proof of Proposition~\ref{Prop: A0 basis}} \label{subsec: A0 proof}

Consider any identity \[ \bfu \approx \bfv \] satisfied by $(A_0, \uast)$.
It suffices to show that $\bfu \approx \bfv$ is deducible from $\{ \eqref{id: inv}, \eqref{id: A0 basis} \}$.
By Lemma~\ref{Lem: A0 bipartite}, this result holds if either~$\bfu$ or~$\bfv$ is bipartite.
Therefore, suppose that~$\bfu$ and~$\bfv$ are both mixed.
By Corollary~\ref{C: A0 violate} and Lemma~\ref{Lem: A0 standard form}, the identities $\{ \eqref{id: inv}, \eqref{id: A0 basis} \}$ can be used to convert~$\bfu$ and~$\bfv$ simultaneously to either $zz^*z$ or words in $A_0$-standard form.
In the former case, $\bfu \approx \bfv$ is deducible from $\{ \eqref{id: inv}, \eqref{id: A0 basis} \}$, whence the proof is complete.
Therefore, it remains to consider the latter case, whence we may assume that~$\bfu$ and~$\bfv$ are in $A_0$-standard form, say \[ \bfu = \bfu_1 x \bfu_2 x^* \quad \text{and} \quad \bfv = \bfv_1 y \bfv_2 y^*, \] where $x,y \in \cXX$, $\bfu_1 = x_1x_2 \cdots x_m$, $\bfu_2 = \bfs_0\prod_{i=1}^p (\bfc_i\bfs_i)$, $\bfv_1 = y_1y_2 \cdots y_n$, and $\bfv_2 = \bft_0\prod_{i=1}^q (\bfd_i\bft_i)$ satisfy conditions (\Asf1)--(\Asf4).

\begin{lemma}\label{Lem: A0 content}
The following holds for the words~$\bfu$ and~$\bfv$\up:
\begin{enumerate}[\ \rm(i)]
\item $\sfcon(\overline{\bfu}) = \sfcon(\overline{\bfv})$\up;
\item $\bfu_1 x =\bfv_1 y$\up;
\item $\sfcon(\overline{\bfu_2}) = \sfcon(\overline{\bfv_2})$.
\end{enumerate}
\end{lemma}

\begin{proof}
(i) Suppose that $\sfcon(\overline{\bfu}) \neq \sfcon(\overline{\bfv})$, say there exists a variable $z \in \sfcon(\overline{\bfu})$ such that $z \notin \sfcon(\overline{\bfv})$.
Then under the substitution $\alpha_\bfv : \cX \to A_0$ in Lemma~\ref{Lem: A0 substitution}, the contradiction $\alpha_\bfv(\bfu) = 0 \neq \rmef = \alpha_\bfv(\bfv)$ is deduced.

\medskip

(ii) Due to condition (\Asf1), the equality $\bfu_1 x =\bfv_1 y$ follows from $\sfcon(\bfu_1 x) = \sfcon(\bfv_1 y)$; to establish the latter, by symmetry, it suffices to verify the inclusion $\sfcon(\bfu_1 x) \subseteq \sfcon(\bfv_1 y)$.
To this end, we need to first show that $y \in \sfcon(\bfu_1x)$.
Since $\overline{y} \in \sfcon(\overline{\bfv}) = \sfcon(\overline{\bfu})$ by part~(i),
\begin{enumerate}[\ (a)]
\item either $y \in \sfcon(\bfu)$ or $y^* \in \sfcon(\bfu)$.
\end{enumerate}
If $y^* \in \sfcon(\bfu_1 x \bfu_2)$, then by Lemma~\ref{Lem: A0 substitution}, we have $\alpha_\bfu(\bfu) = \rmef$ and \[ \alpha_\bfu(\bfv) = \alpha_\bfu(\bfv_1) \cdot \alpha_\bfu(y) \cdot \alpha_\bfu(\bfv_2) \cdot \alpha_\bfu(y^*) = \alpha_\bfu(\bfv_1) \cdot \rmf \cdot \alpha_\bfu(\bfv_2) \cdot \rme = 0, \] which is impossible.
Therefore,
\begin{enumerate}[\ (a)]
\item[(b)] $y^* \notin \sfcon(\bfu_1 x \bfu_2)$, which implies that $y \neq x^*$.
\end{enumerate}
Note that if $y \neq x$, then together with~(b), we have $y^* \notin \sfcon(\bfu_1x\bfu_2x^*) = \sfcon(\bfu)$, so that $y \in \sfcon(\bfu_1x\bfu_2)$ by~(a) and~(b).
On the other hand, if $y = x$, then clearly $y \in \sfcon(\bfu_1x\bfu_2)$.
Therefore, $y \in \sfcon(\bfu_1x\bfu_2)$ either way.
Now if $y \in \sfcon(\bfu_2)$, then by Lemma~\ref{Lem: A0 substitution}, we have $\beta_\bfu(\bfu) = \rmef$ and \[ \beta_\bfu(\bfv) = \beta_\bfu(\bfv_1) \cdot \beta_\bfu(y) \cdot \beta_\bfu(\bfv_2) \cdot \beta_\bfu(y^*) = \beta_\bfu(\bfv_1) \cdot \rmf \cdot \beta_\bfu(\bfv_2) \cdot \rme = 0, \] which is impossible.
Hence $y \notin \sfcon(\bfu_2)$; but since $y \in \sfcon(\bfu_1x\bfu_2)$, we in fact have
\begin{enumerate}[\ (a)]
\item[(c)] $y \in \sfcon(\bfu_1 x)$.
\end{enumerate}

Now we are ready to establish the inclusion $\sfcon(\bfu_1 x) \subseteq \sfcon(\bfv_1 y)$.
Suppose there exists some variable $z \in \sfcon(\bfu_1x)$ such that $z \notin \sfcon(\bfv_1y)$.
Then clearly $z \neq y$.
If $z=y^*$, then $y,y^* \in \sfcon(\bfu_1 x)$ by~(c), so that condition (\Asf1) is contradicted.
Hence
\begin{enumerate}[\ (a)]
\item[(d)] $z \notin \{ y, y^*\}$.
\end{enumerate}
Since $\overline{z} \in \sfcon(\overline{\bfu}) = \sfcon(\overline{\bfv})$ by part~(i), it follows from~(d) that either $z \in \sfcon(\bfv_1\bfv_2)$ or $z^* \in \sfcon(\bfv_1\bfv_2)$.
But since $z \notin \sfcon(\bfv_1)$ by assumption, we have $z \in \sfcon(\bfv_2)$ or $z^* \in \sfcon(\bfv_1)$ or $z^* \in \sfcon(\bfv_2)$.
These three cases are shown in the following to be impossible.
Therefore, the variable~$z$ does not exist, whence the required inclusion $\sfcon(\bfu_1x) \subseteq \sfcon(\bfv_1y)$ is established.

\medskip

\noindent{\sc Case~1:} $z \in \sfcon(\bfv_2)$.
By Lemma~\ref{Lem: A0 substitution}, we have $\beta_\bfv(\bfv_1 y) = \rme$ and $\beta_\bfv(\bfv_2y^*) = \rmf$, so that $\beta_\bfv(\bfv) = \rmef$.
Since $z \in \sfcon(\bfv_2)$, we also have $\beta_\bfv(z) = \rmf$.
Note that \[ \rmef = \beta_\bfv(\bfv) = \beta_\bfv(\bfu) = \beta_\bfv(\bfu_1) \cdot \beta_\bfv(x) \cdot \beta_\bfv(\bfu_2) \cdot \beta_\bfv(x^*), \] so we must have $\beta_\bfv(x^*) = \rmf$ and $\beta_\bfv(x) = \rme$, so that $z \neq x$.
But since $z \in \sfcon(\bfu_1x)$ by assumption, it follows that $\beta_\bfv(\bfu_1x) = \cdots \rmf \cdots \rme = 0$, whence the contradiction $\beta_\bfv(\bfu) = 0$ is deduced.

\medskip

\noindent{\sc Case~2:} $z^* \in \sfcon(\bfv_1)$.
By Lemma~\ref{Lem: A0 substitution}, we have $\alpha_\bfu(\bfu_1 x \bfu_2) = \rme$ and $\alpha_\bfu(x^*) = \rmf$, so that $\alpha_\bfu(\bfu) = \rmef$.
Since $z \in \sfcon(\bfu_1 x)$ by assumption, we also have $\alpha_\bfu(z) = \rme$.
Note that \[ \rmef = \alpha_\bfu(\bfu) = \alpha_\bfu(\bfv) = \alpha_\bfu(\bfv_1) \cdot \alpha_\bfu(y) \cdot \alpha_\bfu(\bfv_2) \cdot \alpha_\bfu(y^*), \] so we must have $\alpha_\bfu(y^*) = \rmf$ and $\alpha_\bfu(y) = \rme$.
But since $z ^*\in \sfcon(\bfv_1)$, it follows that $\alpha_\bfu(\bfv_1y) = \cdots \rmf \cdots \rme = 0$, whence the contradiction $\alpha_\bfu(\bfv) = 0$ is deduced.

\medskip

\noindent{\sc Case~3:} $z^* \in \sfcon(\bfv_2)$.
By Lemma~\ref{Lem: A0 substitution}, we have $\alpha_\bfv(\bfv_1 y \bfv_2) = \rme$ and $\alpha_\bfv(y^*) = \rmf$, so that $\alpha_\bfv(\bfv) = \rmef$.
Since $z^* \in \sfcon(\bfv_2)$, we also have $\alpha_\bfv(z) = \rmf$.
Note that \[ \rmef = \alpha_\bfv(\bfv) = \alpha_\bfv(\bfu) = \alpha_\bfv(\bfu_1) \cdot \alpha_\bfv(x) \cdot \alpha_\bfv(\bfu_2) \cdot \alpha_\bfv(x^*), \] so we must have $\alpha_\bfv(x^*) = \rmf$ and $\alpha_\bfv(x) = \rme$.
But since $z \in \sfcon(\bfu_1x)$ by assumption, it follows that $\alpha_\bfv(\bfu_1x) = \cdots \rmf \cdots \rme = 0$, whence the contradiction $\alpha_\bfv(\bfu) = 0$ is deduced.

\medskip

(iii) This is a consequence of parts~(i) and~(ii).
\end{proof}

Therefore, by Lemma~\ref{Lem: A0 content}, we now have \[ \bfu = { \underbrace{x_1x_2 \cdots x_m}_{\bfu_1} } \cdot x \cdot { \underbrace{\bfs_0\prod_{i=1}^p (\bfc_i\bfs_i)}_{\bfu_2} } \cdot x^* \quad \text{and} \quad \bfv = { \underbrace{x_1x_2 \cdots x_m}_{\bfu_1} } \cdot x \cdot { \underbrace{\bft_0\prod_{i=1}^q (\bfd_i\bft_i)}_{\bfv_2} } \cdot x^*, \] where conditions (\Asf1)--(\Asf4) are satisfied.

\begin{lemma} \label{Lem: A0 con s c}
\begin{enumerate}[\ \rm(i)]
\item $\sfcon(\overline{\bfs_0\bfs_1\cdots \bfs_p})=\sfcon(\overline{\bft_0 \bft_1 \cdots \bft_q})$.
\item $\sfcon(\overline{\bfc_1\bfc_2\cdots \bfc_p})=\sfcon(\overline{\bfd_1 \bfd_2 \cdots \bfd_q})$.
\end{enumerate}
\end{lemma}

\begin{proof}
(i) Let $\bfs=\bfs_0 \bfs_1 \cdots \bfs_p$ and $\bft=\bft_0 \bft_1 \cdots \bft_q.$
Suppose there exists some variable $z \in \sfcon(\overline{\bfs})$ such that $z \notin \sfcon(\overline{\bft})$.
Then \[ \bfu = x_1x_2 \cdots x_m \cdot x \cdot \bfa z^\circledast \bfb \cdot x^*\] for some $\bfa, \bfb \in \Fmon$ and $\circledast \in \{1, *\}$.
By the definition of ${A_0}$-standard form, $\sfocc(z,\overline{\bfu})=1$ and the sets $\sfcon(\overline{x_1x_2 \cdots x_m \cdot x \cdot \bfa}), \{ z \}, \sfcon(\overline{\bfb})$ are pairwise disjoint.
Therefore, if $\varphi: \cXX \to A_0$ is the substitution given by
\[
\varphi(t) = \begin{cases} \rme & \text{if $t \in \sfcon(x_1x_2 \cdots x_m \cdot x \cdot \bfa)$}, \\ \rmef & \text{if $t = z$}, \\ \rmf & \text{otherwise}, \end{cases}
\]
then $\varphi(\bfu) = \rme^m \cdot \rme \cdot \rme^{|\bfa|} (\rmef)^\circledast \rmf^{|\bfb|} \cdot \rme^* = \rmef$.
Now since $z \in \sfcon(\overline{\bfu_2}) = \sfcon(\overline{\bfv_2})$ by Lemma~\ref{Lem: A0 content}(iii) and $z \notin \sfcon(\overline{\bft})$, there exists some $j \in \{ 1,2,\ldots, q\}$ such that $z \in \sfcon(\overline{\bfd_j})$; further, since~$\bfd_j$ is an {\Ablock}, $\sfocc(z, \overline{\bfd_j})=2$.
It follows that $\varphi(\bfv) = \cdots \rmef \cdots \rmef \cdots = 0 \neq \varphi(\bfu)$, which is a contradiction.
Consequently, the variable~$z$ does not exist, so that the inclusion $\sfcon(\overline{\bfs}) \subseteq \sfcon(\overline{\bft})$ holds.
The reverse inclusion $\sfcon(\overline{\bfs}) \supseteq \sfcon(\overline{\bft})$ holds by a symmetrical argument.

(ii) This follows from part~(i) since $\sfcon(\overline{\bfu_2}) = \sfcon(\overline{\bfv_2})$ by Lemma~\ref{Lem: A0 content}(iii).
\end{proof}

\begin{lemma} \label{Lem: A0 yzyz}
\begin{enumerate}[\rm(i)]
\item Suppose that $yzyz \hookrightarrow \bfu_2$ for some $y,z \in \sfcon(\bfc_i)$ with $1 \leq i \leq p$.
Then the longest $\{\overline{y},\overline{z}\}$-{\ssw} of $\bfv_2$ can only be $yzyz$ or $y^*z^*y^*z^*$.
\item Suppose that $yzyz \hookrightarrow \bfv_2$ for some $y,z \in \sfcon(\bfd_i)$ with $1 \leq i \leq q$.
Then the longest $\{\overline{y},\overline{z}\}$-{\ssw} of $\bfu_2$ can only be $yzyz$ or $y^*z^*y^*z^*$.
\end{enumerate}
\end{lemma}

\begin{proof}
(i) By Remark~\ref{Rmk: A0-standard}(v), $yzyz$ is the longest $\{\overline{y},\overline{z}\}$-{\ssw} of~$\bfu$.
Further, $\overline{y} \order \overline{z}$ because $\bfc_i$ is an {\Ablock}.
Hence $\overline{y}, \overline{z} \in \sfcon(\overline{\bfc_i}) \subseteq \sfcon(\overline{\bfd_1 \bfd_2 \cdots \bfd_q})$ by Lemma \ref{Lem: A0 con s c}(ii).
If $\overline{y},\overline{z} \in \sfcon(\overline{\bfd_j})$ for some $j \in \{1,2,\ldots,q\}$, then the longest $\{\overline{y},\overline{z}\}$-{\ssw} of $\bfv_2$ is one of \[ y z y z, \quad y^* z y^* z, \quad y z^* y z^*, \quad y^* z^* y^* z^*; \]
and if $\overline{y} \in \sfcon(\overline{\bfd_j})$ and $\overline{z} \in \sfcon(\overline{\bfd_k})$ for some distinct $j,k \in \{1,2,\ldots,q\}$, then the longest $\{\overline{y},\overline{z}\}$-{\ssw} of $\bfv_2$ is one of
\[ y y z z, \quad y^* y^* z z, \quad y y z^* z^*, \quad y^* y^* z^* z^*, \quad z z y y, \quad z^* z^* y y, \quad z z y^* y^*, \quad z^* z^* y^* y^*. \]
There are four cases to consider.

\medskip

\noindent{\sc Case~1:} $y^*zy^*z$ or $yz^*yz^*$ or $y^* y^* z z$ or $z^*z^*yy$ is a {\ssw} of $\bfv_2$.
Then under the substitution $\alpha_\bfu : \cX \to A_0$ in Lemma~\ref{Lem: A0 substitution}, we have $\alpha_\bfu (\bfu_1 x \bfu_2) = \rme$ and $\alpha_\bfu (x^*) = \rmf$, so that $\alpha_\bfu(\bfu) = \rmef$.
Specifically, $\alpha_\bfu (y) = \alpha_\bfu(z) = \rme$ because $y,z \in \sfcon(\bfu_2)$.
But this implies the contradiction $\alpha_\bfu(\bfv) = 0$.

\medskip

\noindent{\sc Case~2:} either $y y z^* z^* \hookrightarrow \bfv_2$ or $zzy^*y^* \hookrightarrow \bfv_2$.
Then under the substitution $\beta_\bfu : \cX \to A_0$ in Lemma~\ref{Lem: A0 substitution}, we have $\beta_\bfu (\bfu_1 x ) = \rme$ and $\beta_\bfu (\bfu_2x^*) = \rmf$, so that $\beta_\bfu(\bfu) = \rmef$.
Specifically, $\beta_\bfu (y) = \beta_\bfu(z) = \rmf$ because $y,z \in \sfcon(\bfu_2)$.
But this implies the contradiction $\beta_\bfu(\bfv) = 0$.

\medskip

\noindent{\sc Case~3:} either $yyzz \hookrightarrow \bfv_2$ or $y^*y^*z^*z^* \hookrightarrow \bfv_2$.
Then we can write $\bfv_2=\bfa \bfb$ such that $\sfocc(\overline{y}, \overline{\bfa})=2$ and $\sfocc(\overline{z}, \overline{\bfb})=2$.
Let $\varphi: \cXX \rightarrow A_0$ be any substitution that maps each variable in $\sfcon(\bfu_1 x \bfa)$ to~$\rme$ and each variable in $\sfcon(\bfb)$ to $\rmf$.
Then \[\varphi(\bfv) = \varphi(\bfu_1 x \bfa) \cdot \varphi(\bfb) \cdot \varphi(x^*) = \rme \cdot \rme \cdot \rme^* = \rmef. \]
Now depending on whether $yyzz \hookrightarrow \bfv_2$ or $y^*y^*z^*z^* \hookrightarrow \bfv_2$, the pair $(\varphi(y),\varphi(z))$ is either $(\rme,\rmf)$ or $(\rmf,\rme)$; but in either case, \[ \varphi(\bfu) = \cdots \varphi(y) \cdots \varphi(z) \cdots \varphi(y) \cdots \varphi(z) \cdots = 0, \] which is a contradiction.

\medskip

\noindent{\sc Case~4:} either $zzyy \hookrightarrow \bfv_2$ or $z^*z^*y^*y^* \hookrightarrow \bfv_2$.
This is symmetrical to the previous case and so also leads to a contradiction.

\medskip

Since none of the four cases is possible, the longest $\{\overline{y},\overline{z}\}$-{\ssw} of $\bfv_2$ can only be $yzyz$ or $y^*z^*y^*z^*$.

\medskip

(ii) This is symmetrical to part~(i).
\end{proof}

\begin{lemma} \label{Lem: A0 yz}
\begin{enumerate}[\rm(i)]
\item Suppose that $yz \hookrightarrow \bfu_2$ for some $y,z \in \cXX$ that are not in the same {\Ablock}.
Then the $\{\overline{y},\overline{z}\}$-{\ssw}s of $\bfv_2$ of length two can only be $yz$ or $z^*y^*$.
\item Suppose that $yz \hookrightarrow \bfv_2$ for some $y,z \in \cXX$ that are not in the same {\Ablock}.
Then the $\{\overline{y},\overline{z}\}$-{\ssw}s of $\bfu_2$ of length two can only be $yz$ or $z^*y^*$.
\end{enumerate}
\end{lemma}

\begin{proof}
(i) By symmetry, it suffices to assume that within $\bfu_2$, the first~$y$ appears before the first~$z$.
Then by assumption, depending on which of $y$ and $z$ is simple or in an {\Ablock}, the longest $\{\overline{y},\overline{z}\}$-{\ssw} of $\bfu_2$ is $y^rz^s$ for some $r,s \in \{ 1,2\}$.
Since $\bfu$ is in $A_0$-standard form, $\bfu_2 =\bfa y \bfb z \bfe$ for some pairwise disjoint $\bfa, \bfb, \bfe \in \Fmon$ such that $y \notin \sfcon(\bfb\bfe)$ and $z \notin \sfcon(\bfa\bfb)$.
Since $\overline{y},\overline{z} \in \sfcon(\overline{\bfu_2}) = \sfcon(\overline{\bfv_2})$ by Lemma \ref{Lem: A0 con s c}, the $\{\overline{y},\overline{z}\}$-{\ssw}s of $\bfv_2$ of length two can only be \[ yz, \quad y^*z, \quad yz^*, \quad y^*z^*, \quad zy, \quad z^*y, \quad zy^*, \quad z^*y^*. \]
There are three cases to consider.

\medskip

\noindent{\sc Case~1:} either $y^*z \hookrightarrow \bfv_2$ or $z^*y \hookrightarrow \bfv_2$.
Then under the substitution $\alpha_\bfu : \cX \to A_0$ in Lemma~\ref{Lem: A0 substitution}, we have $\alpha_\bfu (\bfu_1 x \bfu_2) = \rme$ and $\alpha_\bfu (x^*) = \rmf$, so that $\alpha_\bfu(\bfu) = \rmef$.
Specifically, $\alpha_\bfu (y) = \alpha_\bfu(z) = \rme$ because $y,z \in \sfcon(\bfu_2)$.
But this implies the contradiction $\alpha_\bfu(\bfv) = 0$.

\medskip

\noindent{\sc Case~2:} either $yz^* \hookrightarrow \bfv_2$ or $zy^* \hookrightarrow \bfv_2$.
Then under the substitution $\beta_\bfu : \cX \to A_0$ in Lemma~\ref{Lem: A0 substitution}, we have $\beta_\bfu (\bfu_1 x ) = \rme$ and $\beta_\bfu (\bfu_2x^*) = \rmf$, so that $\beta_\bfu(\bfu) = \rmef$.
Specifically, $\beta_\bfu (y) = \beta_\bfu(z) = \rmf$ because $y,z \in \sfcon(\bfu_2)$.
But this implies the contradiction $\beta_\bfu(\bfv) = 0$.

\medskip

\noindent{\sc Case~3:} either $y^*z^* \hookrightarrow \bfv_2$ or $zy \hookrightarrow \bfv_2$.
Then under any substitution $\varphi: \cXX \rightarrow A_0$ that maps each variable in $\sfcon(\bfu_1 x \bfa y \bfb)$ to~$\rme$ and each variable in $\sfcon(z \bfe)$ to $\rmf$, we have $\varphi(\bfu) = \varphi (\bfu_1 x \bfa y \bfb) \cdot \varphi (z \bfe) \cdot \varphi (x^*) = \rme \cdot \rmf \cdot \rme^* = \rmef$.
But $\varphi(\bfv) = 0$ is a contradiction.

\medskip

Since none of the three cases is possible, the $\{\overline{y},\overline{z}\}$-{\ssw}s of $\bfv_2$ of length two can only be $yz$ or $z^*y^*$.

\medskip

(ii) This is symmetrical to part~(i).
\end{proof}

\begin{lemma}\label{Lem: A0 head tail}
Suppose that $\bfu_2, \bfv_2 \neq 1$.
Then $\sfh(\bfu_2) =\sfh(\bfv_2)$ and $\sft(\bfu_2) =\sft(\bfv_2)$.
\end{lemma}

\begin{proof}
Recall from Lemma~\ref{Lem: A0 content}(iii) that $\sfcon(\overline{\bfu_2}) = \sfcon(\overline{\bfv_2})$.
First assume that $|\sfcon(\overline{\bfu_2})| = |\sfcon(\overline{\bfv_2})| = 1$, say $\sfcon(\overline{\bfu_2}) = \sfcon(\overline{\bfv_2}) = \{ z \}$ for some $z \in \cX$.
Then each of $\bfu_2$ and $\bfv_2$ can only be simple or an {\Ablock}, so that $\bfu_2, \bfv_2 \in \{ z,z^2\}$ by conditions (\Asf\ref{Asf4 singleton}) and (\Asf\ref{Asf4 s0 orderedsquare}).
Hence $\sfh(\bfu_2)= z =\sfh(\bfv_2)$ and $\sft(\bfu_2)= z =\sft(\bfv_2)$.

Now assume that $|\sfcon(\overline{\bfu_2})| = |\sfcon(\overline{\bfv_2})| \geq 2$.
Let $\headu = \sfh(\bfu_2)$, $\tailu=\sft(\bfu_2)$, $\headv=\sfh(\bfv_2)$, and $\tailv=\sft(\bfv_2)$, so that \begin{equation} \overline{\headu} \order \overline{\tailu} \quad \text{and} \quad \overline{\headv} \order \overline{\tailv} \label{Dis: A0 h<t H<T} \end{equation} by conditions (\Asf\ref{Asf4 s0 nonempty}) and (\Asf\ref{Asf4 s0 empty}).
Recall that $\bfu_2$ and $\bfv_2$ are bipartite words such that $\sfocc(z,\bfu_2), \sfocc(z,\bfv_2) \leq 2$ for all $z \in \cXX$.
There are five cases to consider; in each case, several intermediate results are established to eventually show that $\headu = \headv$ and $\tailu = \tailv$.

\medskip

\noindent{\sc Case~1:} $\sfocc(\headu,\bfu_2) = \sfocc(\tailu,\bfu_2) = 1$.
Then $\headu = \sfh(\bfs_0)$ and $\tailu = \sft(\bfs_p)$, so that $\bfs_0 = \headu \bfa$ and $\bfs_p = \bfb \tailu$ for some $\bfa,\bfb \in \Fmon$:
\begin{equation}
\bfu = x_1 x_2 \cdots x_m \cdot x \cdot \underbrace{\headu \bfa \cdot \bfc_1 \bfs_1 \cdot \bfc_2 \bfs_2 \cdots \bfc_{p-1} \bfs_{p-1} \cdot \bfc_p \bfb \tailu}_{\bfu_2} \cdot \, x^*. \label{Dis: A0 Case1}
\end{equation}

\begin{result} \label{R: Case1: occ}
$\sfocc(\overline{\headu},\overline{\bfv_2}) = \sfocc(\overline{\tailu},\overline{\bfv_2}) = 1$.
\end{result}

\begin{proof}
This holds because $\overline{\headu},\overline{\tailu} \in \sfcon(\overline{\bfs_0\bfs_p}) \subseteq \sfcon(\overline{\bft_0\bft_1 \cdots \bft_q})$ by Lemma~\ref{Lem: A0 con s c}(i).
\end{proof}

\begin{result} \label{R: Case1: ssw}
The longest $\{\overline{\headu},\overline{\tailu}\}$-{\ssw} of $\bfv_2$ is $\headu\tailu$.
\end{result}

\begin{proof}
Since $\headu\tailu \hookrightarrow \bfu_2$, it follows from Lemma~\ref{Lem: A0 yz}(i) and Result~\ref{R: Case1: occ} that the longest $\{\overline{\headu},\overline{\tailu}\}$-{\ssw} of $\bfv_2$ is either $\headu\tailu$ or $\tailu^*\headu^*$.
Seeking a contradiction, suppose that $\tailu^*\headu^*$ is the longest $\{\overline{\headu},\overline{\tailu}\}$-{\ssw} of $\bfv_2$.
Note that $\tailu^*$ does not occur in any {\Ablock} in $\bfv_2$ due to $\sfocc(\overline{\tailu},\overline{\bfv_2}) = 1$ by Result~\ref{R: Case1: occ}.
It follows that if $\headv \neq \tailu^*$, so that $\headv\tailu^* \hookrightarrow \bfv_2$, then by Lemma~\ref{Lem: A0 yz}(ii), either $\headv\tailu^* \hookrightarrow \bfu_2$ or $\tailu\headv^* \hookrightarrow \bfu_2$; but neither {\ssw} is possible in view of \eqref{Dis: A0 Case1}.
By a symmetrical argument, it is impossible for $\tailv \neq \headu^*$.
Therefore, $\headv = \tailu^*$ and $\tailv = \headu^*$.
It follows that $\sfocc(\overline{\headv},\overline{\bfv_2}) = \sfocc(\overline{\tailu},\overline{\bfv_2}) = 1$ and $\headv = \sfh(\bfv_2) = \sfh(\bft_0)$, so that $\bft_0 \neq 1$.
Given that $\bfv$ is in $A_0$-standard form, we have $\overline{\headv} \order \overline{\tailv}$ by condition (\Asf\ref{Asf4 s0 nonempty}); but this implies that $\overline{\tailu} = \overline{\headv} \order \overline{\tailv} = \overline{\headu}$, which contradicts~\eqref{Dis: A0 h<t H<T}.
\end{proof}

\begin{result} \label{R: Case1: h=H t=T}
$\headu = \headv$ and $\tailu = \tailv$.
\end{result}

\begin{proof}
Suppose that $\headu \neq \headv$.
Then it follows from Result~\ref{R: Case1: ssw} that $\headv\headu \hookrightarrow \bfv_2$, and $\headu$ is not in any {\Ablock} in $\bfv_2$ due to Result~\ref{R: Case1: occ}.
Therefore, by Lemma~\ref{Lem: A0 yz}(ii), either $\headv\headu \hookrightarrow \bfu_2$ or $\headu^*\headv^* \hookrightarrow \bfu_2$; but neither {\ssw} is possible in view of \eqref{Dis: A0 Case1}.
By a symmetrical argument, it is impossible for $\tailu \neq \tailv$.
\end{proof}

\noindent{\sc Case~2:} $\sfocc(\headu,\bfu_2) = 2$ and $\sfocc(\tailu,\bfu_2) = 1$.
Then $\bfs_0 = 1$, $\headu = \sfh(\bfc_1)$, and $\tailu = \sft(\bfs_p)$, so that $\bfc_1 = \headu\bfa\headu \bfa$ and $\bfs_p = \bfb \tailu$ for some $\bfa,\bfb \in \Fmon$:
\begin{equation}
\bfu = x_1 x_2 \cdots x_m \cdot x \cdot \underbrace{\headu \bfa\headu \bfa \, \bfs_1 \cdot \bfc_2 \bfs_2 \cdots \bfc_{p-1} \bfs_{p-1} \cdot \bfc_p \bfb \tailu}_{\bfu_2} \cdot \, x^*. \label{Dis: A0 Case2}
\end{equation}

\begin{result} \label{R: Case2: occ}
$\sfocc(\overline{\headu},\overline{\bfv_2}) = 2$ and $\sfocc(\overline{\tailu},\overline{\bfv_2}) = 1$.
\end{result}

\begin{proof}
This holds because $\overline{\headu} \in \sfcon(\overline{\bfc_1}) \subseteq \sfcon(\overline{\bfd_1\bfd_2 \cdots \bfd_q})$ and $\overline{\tailu} \in \sfcon(\overline{\bfs_p}) \subseteq \sfcon(\overline{\bft_0\bft_1 \cdots \bft_q})$ by Lemma \ref{Lem: A0 con s c}.
\end{proof}

\begin{result} \label{R: Case2: ssw}
The longest $\{\overline{\headu},\overline{\tailu}\}$-{\ssw} of $\bfv_2$ is $\headu^2\tailu$.
\end{result}

\begin{proof}
Since $\headu^2\tailu \hookrightarrow \bfu_2$, it follows from Lemma~\ref{Lem: A0 yz}(i) and Result~\ref{R: Case2: occ} that the longest $\{\overline{\headu},\overline{\tailu}\}$-{\ssw} of $\bfv_2$ is either $\headu^2\tailu$ or $\tailu^*(\headu^*)^2$.
Seeking a contradiction, suppose that $\tailu^*(\headu^*)^2$ is the longest $\{\overline{\headu},\overline{\tailu}\}$-{\ssw} of $\bfv_2$.
Note that $\tailu^*$ does not occur in any {\Ablock} in $\bfv_2$ due to $\sfocc(\overline{\tailu},\overline{\bfv_2}) = 1$ by Result~\ref{R: Case2: occ}.
It follows that if $\headv \neq \tailu^*$, so that $\headv\tailu^* \hookrightarrow \bfv_2$, then by Lemma~\ref{Lem: A0 yz}(ii), either $\headv\tailu^* \hookrightarrow \bfu_2$ or $\tailu\headv^* \hookrightarrow \bfu_2$; but neither {\ssw} is possible in view of \eqref{Dis: A0 Case2}.
Hence $\headv = \tailu^*$.

Now suppose that $\tailv \neq \headu^*$, so that $\headu^*\tailv \hookrightarrow \bfv_2$.
Then since $\bfv_2$ is bipartite, we have $\headu,\tailv^* \notin \sfcon(\bfv_2)$.
There are two cases.
\begin{enumerate}[(a)]
\item $\headu^*$ and $\tailv$ are not in the same {\Ablock} in $\bfv_2$.
Then by Lemma~\ref{Lem: A0 yz}(ii), either $\headu^*\tailv \hookrightarrow \bfu_2$ or $\tailv^*\headu \hookrightarrow \bfu_2$.
But it is clear from \eqref{Dis: A0 Case2} that $\headu^*\tailv \not\hookrightarrow \bfu_2$, so only $\tailv^*\headu \hookrightarrow \bfu_2$ holds.
Specifically, $\tailv^*\headu \hookrightarrow \bfc_1 = \headu\bfa\headu\bfa$, so that $\headu\tailv^*\headu\tailv^* \hookrightarrow \bfu_2$.
Therefore, by Lemma~\ref{Lem: A0 yzyz}(i), either $\headu\tailv^*\headu\tailv^* \hookrightarrow \bfv_2$ or $\headu^*\tailv\headu^*\tailv \hookrightarrow \bfv_2$; but the former contradicts $\headu,\tailv^* \notin \sfcon(\bfv_2)$, while the latter contradicts the assumption of the present case.

\item $\headu^*$ and $\tailv$ are in the same {\Ablock} in $\bfv_2$.
Specifically, since $\tailv = \sft(\bfv_2)$, we have $\headu^*,\tailv \in \sfcon(\bfd_q)$, so that $\headu^*\tailv\headu^*\tailv \hookrightarrow \bfv_2$ and $\overline{\headu} \order \overline{\tailv}$.
By Lemma~\ref{Lem: A0 yzyz}(ii), either $\headu^*\tailv\headu^*\tailv \hookrightarrow \bfu_2$ or $\headu\tailv^*\headu\tailv^* \hookrightarrow \bfu_2$.
It follows from \eqref{Dis: A0 Case2} that $\headu^*\tailv\headu^*\tailv \not\hookrightarrow \bfu_2$, so only $\headu\tailv^*\headu\tailv^* \hookrightarrow \bfu_2$ holds, whence $\tailv^* \in \sfcon(\bfc_1)$.
Since $\bfs_0 = 1$ and $\bfu$ is in $A_0$-standard form, by condition (\Asf\ref{Asf4 s0 empty}), we have $\overline{\tailv} \ordereq \overline{\sft(\bfc_1)} \order \overline{\sft(\bfu_2)} = \overline{\tailu}$.
Now the first variable in~$\bfv_2$ is simple because $\sfh(\bfv_2) = \headv = \tailu^*$ and $\sfocc(\overline{\tailu},\overline{\bfv_2}) = 1$; hence $\bft_0 \neq 1$.
Since $\bfv$ is in $A_0$-standard form, by condition (\Asf\ref{Asf4 s0 nonempty}), we have $\overline{\tailu} = \overline{\headv} = \overline{\sfh(\bft_0)} \order \overline{\sft(\bfv_2)} = \overline{\tailv}$, which is a contradiction.
\end{enumerate}
Since neither~(a) nor~(b) is possible, we have $\tailv = \headu^*$.
As observed in~(b), the first variable in~$\bfv_2$ is simple because $\sfh(\bfv_2) = \headv = \tailu^*$ and $\sfocc(\overline{\tailu},\overline{\bfv_2}) = 1$, thus $\bft_0 \neq 1$.
Given that $\bfv$ is in $A_0$-standard form, by condition (\Asf\ref{Asf4 s0 nonempty}), we have $\overline{\tailu} = \overline{\headv} = \overline{\sfh(\bft_0)} \order \overline{\sft(\bfv_2)} = \overline{\tailv} = \overline{\headu}$, which contradicts \eqref{Dis: A0 h<t H<T}.
\end{proof}

\begin{result} \label{R: Case2: h=H t=T}
$\headu = \headv$ and $\tailu = \tailv$.
\end{result}

\begin{proof}
First, suppose that $\tailu \neq \tailv$.
Then $\tailu\tailv \hookrightarrow \bfv_2$ by Result~\ref{R: Case2: ssw}, and $\tailu$ is not in any {\Ablock} in $\bfv_2$ due to Result~\ref{R: Case2: occ}.
Therefore, by Lemma~\ref{Lem: A0 yz}(ii), either $\tailu\tailv \hookrightarrow \bfu_2$ or $\tailv^*\tailu^* \hookrightarrow \bfu_2$; but neither {\ssw} is possible in view of \eqref{Dis: A0 Case2}.
Hence $\tailu = \tailv$.

Now suppose that $\headu \neq \headv$.
Then $\headv\headu \hookrightarrow \bfv_2$ by Result~\ref{R: Case2: ssw}.
Since $\bfv_2$ is bipartite, we have $\headu^*,\headv^* \notin \sfcon(\bfv_2)$.
There are two cases.
\begin{enumerate}[(a)]
\item $\headv$ and $\headu$ are not in the same {\Ablock} in $\bfv_2$.
Then by Lemma~\ref{Lem: A0 yz}(ii), either $\headv\headu \hookrightarrow \bfu_2$ or $\headu^*\headv^* \hookrightarrow \bfu_2$.
But it is clear from \eqref{Dis: A0 Case2} that $\headu^*\headv^* \not\hookrightarrow \bfu_2$, so only $\headv\headu \hookrightarrow \bfu_2$ holds.
Specifically, $\headv\headu \hookrightarrow \bfc_1 = \headu\bfa\headu\bfa$, so that $\headu\headv\headu\headv \hookrightarrow \bfu_2$.
Therefore, by Lemma~\ref{Lem: A0 yzyz}(i), either $\headu\headv\headu\headv \hookrightarrow \bfv_2$ or $\headu^*\headv^*\headu^*\headv^* \hookrightarrow \bfv_2$; but the former contradicts the assumption of the present case, while the latter contradicts $\headu^*,\headv^* \notin \sfcon(\bfv_2)$.

\item $\headv$ and $\headu$ are in the same {\Ablock} in $\bfv_2$.
Then $\headv\headu\headv\headu \hookrightarrow \bfv_2$.
Hence by Lemma~\ref{Lem: A0 yzyz}(ii), either $\headv\headu\headv\headu \hookrightarrow \bfu_2$ or $\headv^*\headu^*\headv^*\headu^* \hookrightarrow \bfu_2$; but neither {\ssw} is possible in view of \eqref{Dis: A0 Case2}.
\end{enumerate}
Since neither~(a) nor~(b) is possible, we have $\headu = \headv$.
\end{proof}

\medskip

\noindent{\sc Case~3:} $\sfocc(\headu,\bfu_2) = 1$ and $\sfocc(\tailu,\bfu_2) = 2$.
Then $\headu = \sfh(\bfs_0)$, $\tailu = \sft(\bfc_p)$, and $\bfs_p = 1$, so that $\bfs_0 = \headu\bfa$ and $\bfc_p = \bfb\tailu\bfb\tailu$ for some $\bfa,\bfb \in \Fmon$:
\begin{equation}
\bfu = x_1 x_2 \cdots x_m \cdot x \cdot \underbrace{\headu\bfa \cdot \bfc_1 \bfs_1 \cdot \bfc_2 \bfs_2 \cdots \bfc_{p-1} \bfs_{p-1} \cdot \bfb\tailu\bfb\tailu}_{\bfu_2} \cdot \, x^*. \label{Dis: A0 Case3}
\end{equation}

\begin{result} \label{R: Case3: occ}
$\sfocc(\overline{\headu},\overline{\bfv_2}) = 1$ and $\sfocc(\overline{\tailu},\overline{\bfv_2}) = 2$.
\end{result}

\begin{proof}
This holds because $\overline{\headu} \in \sfcon(\overline{\bfs_0}) \subseteq \sfcon(\overline{\bft_0\bft_1 \cdots \bft_q})$ and $\overline{\tailu} \in \sfcon(\overline{\bfc_p}) \subseteq \sfcon(\overline{\bfd_1\bfd_2 \cdots \bfd_q})$ by Lemma~\ref{Lem: A0 con s c}.
\end{proof}

\begin{result} \label{R: Case3: ssw}
The longest $\{\overline{\headu},\overline{\tailu}\}$-{\ssw} of $\bfv_2$ is $\headu\tailu^2$.
\end{result}

\begin{proof}
Since $\headu\tailu^2 \hookrightarrow \bfu_2$, it follows from Lemma~\ref{Lem: A0 yz}(i) and Result~\ref{R: Case3: occ} that the longest $\{\overline{\headu},\overline{\tailu}\}$-{\ssw} of $\bfv_2$ is either $\headu\tailu^2$ or $(\tailu^*)^2\headu^*$.
Seeking a contradiction, suppose that $(\tailu^*)^2\headu^*$ is the longest $\{\overline{\headu},\overline{\tailu}\}$-{\ssw} of $\bfv_2$.
Note that $\headu^*$ does not occur in any {\Ablock} in $\bfv_2$ due to $\sfocc(\overline{\headu},\overline{\bfv_2}) = 1$ by Result~\ref{R: Case3: occ}.
It follows that if $\tailv \neq \headu^*$, so that $\headu^*\tailv \hookrightarrow \bfv_2$, then by Lemma~\ref{Lem: A0 yz}(ii), either $\headu^*\tailv \hookrightarrow \bfu_2$ or $\tailv^*\headu \hookrightarrow \bfu_2$; but neither {\ssw} is possible in view of \eqref{Dis: A0 Case3}.
Hence $\tailv = \headu^*$.

Now suppose that $\headv \neq \tailu^*$, so that $\headv\tailu^* \hookrightarrow \bfv_2$.
Then since $\bfv_2$ is bipartite, $\headv^*,\tailu \notin \sfcon(\bfv_2)$.
There are two cases.
\begin{enumerate}[(a)]
\item $\headv$ and $\tailu^*$ are not in the same {\Ablock} in $\bfv_2$.
Then by Lemma~\ref{Lem: A0 yz}(ii), either $\headv\tailu^* \hookrightarrow \bfu_2$ or $\tailu\headv^* \hookrightarrow \bfu_2$.
But it is clear from \eqref{Dis: A0 Case3} that $\headv\tailu^* \not\hookrightarrow \bfu_2$, so only $\tailu\headv^* \hookrightarrow \bfu_2$ holds.
Specifically, $\tailu\headv^* \hookrightarrow \bfc_p = \bfb\tailu\bfb\tailu$, so that $\headv^*\tailu\headv^*\tailu \hookrightarrow \bfu_2$.
Therefore, by Lemma~\ref{Lem: A0 yzyz}(i), either $\headv^*\tailu\headv^*\tailu \hookrightarrow \bfv_2$ or $\headv\tailu^*\headv\tailu^* \hookrightarrow \bfv_2$; but the former contradicts $\headv^*,\tailu \notin \sfcon(\bfv_2)$, while the latter contradicts the assumption of the present case.

\item $\headv$ and $\tailu^*$ are in the same {\Ablock} in $\bfv_2$.
Specifically, since $\headv = \sfh(\bfv_2)$, we have $\bft_0=1$ and $\headv,\tailu^* \in \sfcon(\bfd_1)$.
Given that $\bfv$ is in $A_0$-standard form, by condition (\Asf\ref{Asf4 s0 empty}), we have $\overline{t} \ordereq \overline{\sft(\bfd_1)} \order \overline{\sft(\bfv_2)} = \overline{\tailv} = \overline{\headu}$, but this contradicts \eqref{Dis: A0 h<t H<T}.
\end{enumerate}
Since neither~(a) nor~(b) is possible, we must have $\headv = \tailu^*$.
Now since $\overline{\headv} \order \overline{\tailv}$ by \eqref{Dis: A0 h<t H<T}, we have $\overline{\tailu} = \overline{\headv} \order \overline{\tailv} = \overline{\headu}$; but this contradicts $\overline{\headu} \order \overline{\tailu}$ in \eqref{Dis: A0 h<t H<T}.
\end{proof}

\begin{result} \label{R: Case3: h=H t=T}
$\headu = \headv$ and $\tailu = \tailv$.
\end{result}

\begin{proof}
First, suppose that $\headu \neq \headv$.
Then $\headv\headu \hookrightarrow \bfv_2$ by Result~\ref{R: Case3: ssw}, and $\headu$ is not in any {\Ablock} in $\bfv_2$ due to Result~\ref{R: Case3: occ}.
Therefore, by Lemma~\ref{Lem: A0 yz}(ii), either $\headv\headu \hookrightarrow \bfu_2$ or $\headu^*\headv^* \hookrightarrow \bfu_2$; but this is impossible in view of \eqref{Dis: A0 Case3}.
Hence $\headu = \headv$.

Now suppose that $\tailu \neq \tailv$.
Then $\tailu\tailv \hookrightarrow \bfv_2$ by Result~\ref{R: Case3: ssw}, and $\tailu^*,\tailv^* \notin \sfcon(\bfv_2)$ due to $\bfv_2$ being bipartite.
If $\tailu$ and $\tailv$ are in the same {\Ablock} in $\bfv_2$, so that $\tailu\tailv\tailu\tailv \hookrightarrow \bfv_2$, then by Lemma~\ref{Lem: A0 yzyz}(ii), either $\tailu\tailv\tailu\tailv \hookrightarrow \bfu_2$ or $\tailu^*\tailv^*\tailu^*\tailv^* \hookrightarrow \bfu_2$; but neither {\ssw} is possible in view of \eqref{Dis: A0 Case3}.
Therefore, $\tailu$ and $\tailv$ are not in the same {\Ablock} in $\bfv_2$.
Then by Lemma~\ref{Lem: A0 yz}(ii), either $\tailu\tailv \hookrightarrow \bfu_2$ or $\tailv^*\tailu^* \hookrightarrow \bfu_2$.
But it is clear from \eqref{Dis: A0 Case3} that $\tailv^*\tailu^* \not\hookrightarrow \bfu_2$, so only $\tailu\tailv \hookrightarrow \bfu_2$ holds.
Specifically, $\tailu\tailv \hookrightarrow \bfc_p = \bfb\tailu\bfb\tailu$, so that $\tailv\tailu\tailv\tailu \hookrightarrow \bfu_2$.
Hence by Lemma~\ref{Lem: A0 yzyz}(i), either $\tailv\tailu\tailv\tailu \hookrightarrow \bfv_2$ or $\tailv^*\tailu^*\tailv^*\tailu^* \hookrightarrow \bfv_2$; but the former contradicts $\tailu$ and $\tailv$ not being in the same {\Ablock} in $\bfv_2$, while the latter contradicts $\tailu^*,\tailv^* \notin \sfcon(\bfv_2)$.
\end{proof}

\medskip

\noindent{\sc Case~4:} $\sfocc(\headu,\bfu_2) = \sfocc(\tailu,\bfu_2) = 2$ with $p \geq 2$.
Then $\bfs_0 = 1$, $\headu = \sfh(\bfc_1)$, $\tailu = \sft(\bfc_p)$, and $\bfs_p = 1$, so that $\bfc_1 = \headu\bfa\headu\bfa$ and $\bfc_p = \bfb\tailu\bfb\tailu$ for some $\bfa,\bfb \in \Fmon$:
\begin{equation}
\bfu = x_1 x_2 \cdots x_m \cdot x \cdot \underbrace{\headu\bfa\headu\bfa \, \bfs_1 \cdot \bfc_2 \bfs_2 \cdots \bfc_{p-1} \bfs_{p-1} \cdot \bfb\tailu\bfb\tailu}_{\bfu_2} \cdot \, x^*. \label{Dis: A0 Case4}
\end{equation}

\begin{result} \label{R: Case4: occ}
$\sfocc(\overline{\headu},\overline{\bfv_2}) = \sfocc(\overline{\tailu},\overline{\bfv_2}) = 2$.
\end{result}

\begin{proof}
This holds because $\overline{\headu},\overline{\tailu} \in \sfcon(\overline{\bfc_1\bfc_p}) \subseteq \sfcon(\overline{\bfd_1\bfd_2 \cdots \bfd_q})$ by Lemma~\ref{Lem: A0 con s c}(ii).
\end{proof}

\begin{result} \label{R: Case4: headu2tailu2 or tailu*2headu*2}
The longest $\{\overline{\headu},\overline{\tailu}\}$-{\ssw} of $\bfv_2$ is either $\headu^2\tailu^2$ or $(\tailu^*)^2(\headu^*)^2$.
\end{result}

\begin{proof}
Since $\headu^2\tailu^2 \hookrightarrow \bfu_2$, it follows from Lemma~\ref{Lem: A0 yz}(i) and Result~\ref{R: Case4: occ} that the longest $\{\overline{\headu},\overline{\tailu}\}$-{\ssw} of $\bfv_2$ is one of the following six words: $\headu^2\tailu^2$,\, $\headu\tailu\headu\tailu$,\, $\tailu\headu\tailu\headu$,\, $(\tailu^*)^2(\headu^*)^2$,\, $\tailu^*\headu^*\tailu^*\headu^*$, and $\headu^*\tailu^*\headu^*\tailu^*$.
If either $\headu\tailu\headu\tailu \hookrightarrow \bfv_2$ or $\headu^*\tailu^*\headu^*\tailu^* \hookrightarrow \bfv_2$, then it follows from Lemma~\ref{Lem: A0 yzyz}(ii) that either $\headu\tailu\headu\tailu \hookrightarrow \bfu_2$ or $\headu^*\tailu^*\headu^*\tailu^* \hookrightarrow \bfu_2$; but neither {\ssw} is possible in view of \eqref{Dis: A0 Case4}.
If $\tailu\headu\tailu\headu \hookrightarrow \bfv_2$, then $\tailu$ and $\headu$ are in the same {\Ablock} in $\bfv_2$, whence $\overline{\tailu} \order \overline{\headu}$; but this contradicts \eqref{Dis: A0 h<t H<T}.
A similar contradiction is obtained if $\tailu^*\headu^*\tailu^*\headu^* \hookrightarrow \bfv_2$.
\end{proof}

\begin{result} \label{R: Case4: tailu<tailv}
Suppose that $(\tailu^*)^2(\headu^*)^2 \hookrightarrow \bfv_2$ and $\headv \neq \tailu^*$.
Then $\overline{\tailu} \order \overline{\tailv}$.
\end{result}

\begin{proof}
By assumption, $\headv\tailu^* \hookrightarrow \bfv_2$.
Since $\bfv_2$ is bipartite, we have $\headv^*,\tailu \notin \sfcon(\bfv_2)$.
Suppose that $\headv$ and $\tailu^*$ are not in the same {\Ablock} in $\bfv_2$.
Then by Lemma~\ref{Lem: A0 yz}(ii), either $\headv\tailu^* \hookrightarrow \bfu_2$ or $\tailu\headv^* \hookrightarrow \bfu_2$.
But it is clear from \eqref{Dis: A0 Case4} that $\headv\tailu^* \not\hookrightarrow \bfu_2$, so only $\tailu\headv^* \hookrightarrow \bfu_2$ holds.
Specifically, $\tailu\headv^*  \hookrightarrow\bfc_p = \bfb\tailu\bfb\tailu$, so that $\headv^*\tailu\headv^*\tailu \hookrightarrow \bfu_2$.
Therefore, by Lemma~\ref{Lem: A0 yzyz}(i), either $\headv^*\tailu\headv^*\tailu \hookrightarrow \bfv_2$ or $\headv\tailu^*\headv\tailu^* \hookrightarrow \bfv_2$; but the former contradicts $\headv^*,\tailu \notin \sfcon(\bfv_2)$, while the latter contradicts $\headv$ and $\tailu^*$ not being in the same {\Ablock} in $\bfv_2$.

Therefore, $\headv$ and $\tailu^*$ are in the same {\Ablock} in $\bfv_2$.
Specifically, since $\headv = \sfh(\bfv_2)$, we have $\headv,\tailu^* \in \sfcon(\bfd_1)$ and $\bft_0 = 1$.
Given that $\bfv$ is in $A_0$-standard form, it follows from condition (\Asf\ref{Asf4 s0 empty}) that $\overline{\tailu} \ordereq \overline{\sft(\bfd_1)} \order \overline{\sft(\bfv_2)} = \overline{\tailv}$.
\end{proof}

\begin{result} \label{R: Case4: tailv<tailu}
Suppose that $(\tailu^*)^2(\headu^*)^2 \hookrightarrow \bfv_2$ and $\tailv \neq \headu^*$.
Then $\overline{\tailv} \order \overline{\tailu}$.
\end{result}

\begin{proof}
By assumption, $\headu^*\tailv \hookrightarrow \bfv_2$.
There are two cases.
\begin{enumerate}[(a)]
\item $\tailv$ and $\headu^*$ are not in the same {\Ablock} in $\bfv_2$.
Then by Lemma~\ref{Lem: A0 yz}(ii), either $\headu^*\tailv \hookrightarrow \bfu_2$ or $\tailv^*\headu \hookrightarrow \bfu_2$.
It is clear from \eqref{Dis: A0 Case4} that $\headu^*\tailv \not\hookrightarrow \bfu_2$, so only $\tailv^*\headu \hookrightarrow \bfu_2$ holds.
Specifically, $\tailv^*\headu \hookrightarrow \bfc_1 = \headu\bfa\headu\bfa$.

\item $\tailv$ and $\headu^*$ are in the same {\Ablock} in $\bfv_2$.
Then $\headu^*\tailv\headu^*\tailv \hookrightarrow \bfv_2$.
By Lemma~\ref{Lem: A0 yzyz}(ii), either $\headu^*\tailv\headu^*\tailv \hookrightarrow \bfu_2$ or $\headu\tailv^*\headu\tailv^* \hookrightarrow \bfu_2$.
It is clear from \eqref{Dis: A0 Case4} that $\headu^*\tailv\headu^*\tailv \not\hookrightarrow \bfu_2$, so only $\headu\tailv^*\headu\tailv^* \hookrightarrow \bfu_2$ holds, whence $\headu\tailv^*\headu\tailv^* \hookrightarrow \bfc_1 = \headu\bfa\headu\bfa$.
\end{enumerate}
Therefore, in any case, we have $\tailv^* \in \sfcon(\bfc_1)$.
Since $\bfs_0 = 1$ and $\bfu$ is in $A_0$-standard form, by condition (\Asf\ref{Asf4 s0 empty}), we have $\overline{\tailv} \ordereq \overline{\sft(\bfc_1)} \order \overline{\sft(\bfu_2)} = \overline{\tailu}$.
\end{proof}

\begin{result} \label{R: Case4: ssw}
The longest $\{\overline{\headu},\overline{\tailu}\}$-{\ssw} of $\bfv_2$ is $\headu^2\tailu^2$.
\end{result}

\begin{proof}
By Result~\ref{R: Case4: headu2tailu2 or tailu*2headu*2}, the longest $\{\overline{\headu},\overline{\tailu}\}$-{\ssw} of $\bfv_2$ is either $\headu^2\tailu^2$ or $(\tailu^*)^2(\headu^*)^2$.
Seeking a contradiction, suppose that $(\tailu^*)^2(\headu^*)^2 \hookrightarrow \bfv_2$.
If $\headv \neq \tailu^*$, then $\overline{\tailu} \order \overline{\tailv}$ by Result~\ref{R: Case4: tailu<tailv}, whence $\tailv = \headu^*$ by Result~\ref{R: Case4: tailv<tailu}; but this implies that $\overline{\tailu} \order \overline{\tailv} = \overline{\headu}$, which contradicts \eqref{Dis: A0 h<t H<T}.
If $\tailv \neq \headu^*$, then $\overline{\tailv} \order \overline{\tailu}$ by Result~\ref{R: Case4: tailv<tailu}, whence $\headv = \tailu^*$ by Result~\ref{R: Case4: tailu<tailv}; but this implies that $\overline{\tailv} \order \overline{\tailu} = \overline{\headv}$, which contradicts \eqref{Dis: A0 h<t H<T} again.
Therefore, we must have $\headv = \tailu^*$ and $\tailv = \headu^*$.
Now since $\overline{\headv} \order \overline{\tailv}$ by \eqref{Dis: A0 h<t H<T}, we have $\overline{\tailu} = \overline{\headv} \order \overline{\tailv} = \overline{\headu}$; but this contradicts $\overline{\headu} \order \overline{\tailu}$ in \eqref{Dis: A0 h<t H<T}.
\end{proof}

\begin{result} \label{R: Case4: h=H t=T}
$\headu = \headv$ and $\tailu = \tailv$.
\end{result}

\begin{proof}
Seeking a contradiction, suppose that $\headu \neq \headv$.
Then $\headv\headu \hookrightarrow \bfv_2$ by Result~\ref{R: Case4: ssw}, and $\headv^*,\headu^* \notin \sfcon(\bfv_2)$ due to~$\bfv_2$ being bipartite.
If $\headv$ and $\headu$ are in the same {\Ablock} in $\bfv_2$, so that $\headv\headu\headv\headu \hookrightarrow \bfv_2$, then by Lemma~\ref{Lem: A0 yzyz}(ii), either $\headv\headu\headv\headu \hookrightarrow \bfu_2$ or $\headv^*\headu^*\headv^*\headu^* \hookrightarrow \bfu_2$; but neither {\ssw} is possible in view of \eqref{Dis: A0 Case4}.
Therefore, $\headv$ and $\headu$ are not in the same {\Ablock} in $\bfv_2$.
Then by Lemma~\ref{Lem: A0 yz}(ii), either $\headv\headu \hookrightarrow \bfu_2$ or $\headu^*\headv^* \hookrightarrow \bfu_2$.
It is clear from \eqref{Dis: A0 Case4} that $\headu^*\headv^* \not\hookrightarrow \bfu_2$, so only $\headv\headu \hookrightarrow \bfu_2$ holds.
Specifically, $\headv\headu \hookrightarrow \bfc_1 = \headu\bfa\headu\bfa$, so that $\headu\headv\headu\headv \hookrightarrow \bfu_2$.
Therefore, by Lemma~\ref{Lem: A0 yzyz}(i), either $\headu\headv\headu\headv \hookrightarrow \bfv_2$ or $\headu^*\headv^*\headu^*\headv^* \hookrightarrow \bfv_2$; but the former contradicts $\headv$ and $\headu$ not being in the same {\Ablock} in $\bfv_2$, while the latter contradicts $\headv^*,\headu^* \notin \sfcon(\bfv_2)$.

A symmetrical argument shows that the assumption $\tailu \neq \tailv$ also leads to a contradiction.
\end{proof}

\medskip

\noindent{\sc Case~5:} $\sfocc(\headu,\bfu_2) = \sfocc(\tailu,\bfu_2) = 2$ with $p = 1$.
Then $\bfs_0 = 1$, $\headu = \sfh(\bfc_1)$, $\tailu = \sft(\bfc_1)$, and $\bfs_1 = 1$, so that $\bfc_1 = \headu\bfa\tailu\headu\bfa\tailu$ for some $\bfa \in \Fmon$:
\begin{equation}
\bfu = x_1 x_2 \cdots x_m \cdot x \cdot \underbrace{\,\headu\bfa\tailu\headu\bfa\tailu\,}_{\bfu_2} \cdot \, x^*. \label{Dis: A0 Case5}
\end{equation}
Note that due to condition (\Asf\ref{Asf4 s0 orderedsquare}), we have $h \in \cX$.

\begin{result} \label{R: Case5: occ}
$\sfocc(\overline{y},\overline{\bfv_2}) = 2$ for all $y \in \sfcon(\bfv_2)$.
\end{result}

\begin{proof}
Since $\bfs_0 = \bfs_1 = 1$ and $\sfcon(\overline{\bft_0\bft_1\cdots\bft_q}) = \sfcon(\overline{\bfs_0\bfs_1})$ by Lemma~\ref{Lem: A0 con s c}(i), we have $\bft_0 = \bft_1 = \cdots = \bft_q = 1$.
Therefore, $\bfv_2 = \bfd_1 \bfd_2 \cdots \bfd_q$ and the result follows.
\end{proof}

\begin{result} \label{R: Case5: ssw}
The longest $\{\overline{\headu},\overline{\tailu}\}$-{\ssw} of $\bfv_2$ is $\headu\tailu\headu\tailu$.
\end{result}

\begin{proof}
Since $\headu\tailu\headu\tailu \hookrightarrow \bfu_2$, it follows from Lemma~\ref{Lem: A0 yzyz}(i) that the longest $\{\overline{\headu},\overline{\tailu}\}$-{\ssw} of $\bfv_2$ is either $\headu\tailu\headu\tailu$ or $\headu^*\tailu^*\headu^*\tailu^*$.
Seeking a contradiction, suppose that $\headu^*\tailu^*\headu^*\tailu^* \hookrightarrow \bfv_2$.

First consider the case when $\headv \neq \headu^*$, so that $\headv\headu^* \hookrightarrow \bfv_2$.
Then $\headv^*,\headu \notin \sfcon(\bfv_2)$ due to $\bfv_2$ being bipartite.
If $\headv$ and $\headu^*$ are in the same {\Ablock} in $\bfv_2$, so that $\headv\headu^*\headv\headu^* \hookrightarrow \bfv_2$, then it follows from Lemma~\ref{Lem: A0 yzyz}(ii) that either $\headv\headu^*\headv\headu^* \hookrightarrow \bfu_2$ or $\headv^*\headu\headv^*\headu \hookrightarrow \bfu_2$; but neither {\ssw} is possible in view of \eqref{Dis: A0 Case5}.
Therefore, $\headv$ and $\headu^*$ are not in the same {\Ablock} in $\bfv_2$.
Then by Lemma~\ref{Lem: A0 yz}(ii), either $\headv\headu^* \hookrightarrow \bfu_2$ or $\headu\headv^* \hookrightarrow \bfu_2$.
It is clear from \eqref{Dis: A0 Case5} that $\headv\headu^* \not\hookrightarrow \bfu_2$, so only $\headu\headv^* \hookrightarrow \bfu_2$ holds.
It follows that $\headu\headv^*\headu\headv^* \hookrightarrow \bfu_2$, so that by Lemma~\ref{Lem: A0 yzyz}(i), either $\headu\headv^*\headu\headv^* \hookrightarrow \bfv_2$ or $\headu^*\headv\headu^*\headv \hookrightarrow \bfv_2$; but the former contradicts $\headv^*,\headu \notin \sfcon(\bfv_2)$, while the latter contradicts $\headv$ and $\headu^*$ not being in the same {\Ablock} in~$\bfv_2$.

Therefore, $\headv = \headu^*$.
By a symmetrical argument, we have $\tailv = \tailu^*$.
It follows that $\headv\tailv\headv\tailv = \headu^*\tailu^*\headu^*\tailu^* \hookrightarrow \bfv_2$, so that $\bfv_2 = \bfd_1 = \headu^* \bfb \tailu^*\headu^* \bfb \tailu^*$ for some $\bfb \in \Fmon$ (with $q=1$).
As deduced in the proof of Result~\ref{R: Case5: occ}, we have $\bft_0 = \bft_1 = 1$.
Since~$\bfv$ is in $A_0$-standard form, by condition (\Asf\ref{Asf4 s0 orderedsquare}), we have $\headu^* = \sfh(\bfd_1) \in \cX$, which contradicts the observation $\headu \in \cX$ made after \eqref{Dis: A0 Case5}.
\end{proof}

\begin{result} \label{R: Case5: h=H t=T}
$\headu = \headv$ and $\tailu = \tailv$.
\end{result}

\begin{proof}
Seeking a contradiction, suppose that $\headu \neq \headv$.
Then $\headv\headu \hookrightarrow \bfv_2$ by Result~\ref{R: Case5: ssw}, and $\headv^*,\headu^* \notin \sfcon(\bfv_2)$ due to $\bfv_2$ being bipartite.
If $\headv$ and $\headu$ are in the same {\Ablock} in $\bfv_2$, so that $\headv\headu\headv\headu \hookrightarrow \bfv_2$, then it follows from Lemma~\ref{Lem: A0 yzyz}(ii) that either $\headv\headu\headv\headu \hookrightarrow \bfu_2$ or $\headv^*\headu^*\headv^*\headu^* \hookrightarrow \bfu_2$; but neither {\ssw} is possible in view of \eqref{Dis: A0 Case5}.
Therefore, $\headv$ and $\headu$ are not in the same {\Ablock} in $\bfv_2$.
Then by Lemma~\ref{Lem: A0 yz}(ii), either $\headv\headu \hookrightarrow \bfu_2$ or $\headu^*\headv^* \hookrightarrow \bfu_2$.
It is clear from \eqref{Dis: A0 Case5} that $\headu^*\headv^* \not\hookrightarrow \bfu_2$, so only $\headv\headu \hookrightarrow \bfu_2$ holds.
It follows that $\headu\headv\headu\headv \hookrightarrow \bfu_2$, so that by Lemma~\ref{Lem: A0 yzyz}(i), either $\headu\headv\headu\headv \hookrightarrow \bfv_2$ or $\headu^*\headv^*\headu^*\headv^* \hookrightarrow \bfv_2$; but the former contradicts $\headv$ and $\headu$ not being in the same {\Ablock} in $\bfv_2$, while the latter contradicts $\headv^*,\headu^* \notin \sfcon(\bfv_2)$.

A contradiction can be similarly deduced if $\tailu \neq \tailv$.
\end{proof}
In conclusion, we have $\headu = \headv$ and $\tailu = \tailv$ in all five cases (Results~\ref{R: Case1: h=H t=T}, \ref{R: Case2: h=H t=T}, \ref{R: Case3: h=H t=T}, \ref{R: Case4: h=H t=T}, and~\ref{R: Case5: h=H t=T}).
The proof of Lemma~\ref{Lem: A0 head tail} is thus complete.
\end{proof}

\begin{lemma} \label{Lem: A0 satisfies u2=v2}
The identity $\bfu_2 \approx \bfv_2$ is satisfied by $(A_0,\uast)$.
\end{lemma}

\begin{proof}
Recall from Lemma~\ref{Lem: A0 content}(iii) that $\sfcon(\overline{\bfu_2}) = \sfcon(\overline{\bfv_2})$.
As shown in the beginning of the proof of Lemma~\ref{Lem: A0 head tail}, if $|\sfcon(\overline{\bfu_2})| = |\sfcon(\overline{\bfv_2})| = 1$, then $\bfu_2, \bfv_2 \in \{ z,z^2\}$ for some $z \in \cX$, so that $\bfu_2 = \bfv_2$ by Lemma~\ref{Lem: A0 con s c} and conditions (\Asf\ref{Asf4 singleton}) and (\Asf\ref{Asf4 s0 orderedsquare}).
Therefore, it suffices to assume that $|\sfcon(\overline{\bfu_2})| = |\sfcon(\overline{\bfv_2})| \geq 2$, so that by Lemma~\ref{Lem: A0 head tail}, we have $\headu = \sfh(\bfu_2) = \sfh(\bfv_2)$ and $\tailu = \sft(\bfu_2) = \sft(\bfv_2)$ with $\overline{\headu} \order \overline{\tailu}$.
Specifically, $\bfu_2 = \headu \bfa \tailu$ and $\bfv_2 = \headu \bfb \tailu$ for some $\bfa,\bfb \in \Fmon$.

Seeking a contradiction, suppose that $\bfu_2 \approx \bfv_2$ is not satisfied by $(A_0,\uast)$.
Then there exists a substitution $\psi: \sfcon(\overline{\bfu_2}) \to A_0$ such that $\psi(\bfu_2) \neq \psi(\bfv_2)$, so that \begin{equation} \psi(\headu) \cdot \psi(\bfa) \cdot \psi(\tailu) \neq \psi(\headu) \cdot \psi(\bfb) \cdot \psi(\tailu). \label{Dis: A0 satisfies u2=v2} \end{equation}
Clearly, $\psi(\headu) \neq 0 \neq \psi(\tailu)$.
Further, it is routinely checked that
\begin{align*}
\rme \cdot A_0^1 \cdot \rme & = \{ 0,\rme\}, & \rme \cdot A_0^1 \cdot \rmf & = \{ 0, \rmef\}, & \rme \cdot A_0^1 \cdot \rmef & = \{ 0, \rmef \}, \\
\rmf \cdot A_0^1 \cdot \rme & = \{ 0\}, & \rmf \cdot A_0^1 \cdot \rmf & = \{ 0, \rmf\}, & \rmf \cdot A_0^1 \cdot \rmef & = \{ 0\}, \\
\rmef \cdot A_0^1 \cdot \rme & = \{ 0\}, & \rmef \cdot A_0^1 \cdot \rmf & = \{ 0, \rmef\}, & \rmef \cdot A_0^1 \cdot \rmef & = \{ 0\}.
\end{align*}
Thus for \eqref{Dis: A0 satisfies u2=v2} to hold, we need $(\psi(\headu),\psi(\tailu)) \in \{ (\rme,\rme), (\rme,\rmf), (\rme,\rmef), (\rmf,\rmf), (\rmef,\rmf) \}$, whence $\{ \psi(\bfu_2), \psi(\bfv_2) \}$ can be $\{ 0, \rme \}$, $\{ 0, \rmf \}$, or $\{ 0,\rmef\}$.
Generality is not lost by assuming that $\psi(\bfu_2) = 0$ and $\psi(\bfv_2) \in \{ \rme,\rmf,\rmef\}$.
Now extend $\psi$ to the substitution $\Psi$ that maps every $z \in \sfcon(\overline{\bfu_2}) = \sfcon(\overline{\bfv_2})$ to $\psi(z)$ and every $z \in \{ x_1,x_2,\ldots,x_m,x\}$ to~$\rme$.
Then $\Psi(\bfu) \neq \Psi(\bfv)$ because
\begin{align*}
\Psi(\bfu) & = \Psi(x_1) \cdot \Psi(x_2) \cdots \Psi(x_m) \cdot \Psi(x) \cdot \psi(\bfu_2) \cdot \Psi(x)^* = \rme \cdot 0 \cdot \rmf = 0  \\
\text{and } \ \Psi(\bfv) & = \Psi(x_1) \cdot \Psi(x_2) \cdots \Psi(x_m) \cdot \Psi(x) \cdot \psi(\bfv_2) \cdot \Psi(x)^* = \rme \cdot \psi(\bfv_2) \cdot \rmf = \rmef;
\end{align*}
but this is impossible given that $\bfu \approx \bfv$ is satisfied by $(A_0,\uast)$.
\end{proof}

Since the words~$\bfu_2$ and~$\bfv_2$ are bipartite, it follows from Lemmas~\ref{Lem: A0 bipartite} and~\ref{Lem: A0 satisfies u2=v2} that $\bfu_2 \approx \bfv_2$ is deducible from $\{ \eqref{id: inv}, \eqref{id: A0 basis} \}$.
Since \[ \bfu = x_1x_2 \cdots x_m \cdot x \cdot \bfu_2 \cdot x^* \quad \text{and} \quad \bfv = x_1x_2 \cdots x_m \cdot x \cdot \bfv_2 \cdot x^*, \] the identity $\bfu \approx \bfv$ is also deducible from $\{ \eqref{id: inv}, \eqref{id: A0 basis} \}$.
The proof of Proposition~\ref{Prop: A0 basis} is thus complete.

\section{The {\insem} $(S_4,\uast)$} \label{sec: B0}

The {\insem} $(S_4,\uast)$ is isomorphic to the semigroup \[ B_0 = \langle \rma,\rme,\rmf \,|\, \rma\rmf = \rme\rma = \rma, \,  \rme^2=\rme, \, \rmf^2 = \rmf, \, \rme\rmf = \rmf\rme = 0 \rangle = \{ 0, \rma,\rme,\rmf \} \] with the operation~$\uast$ that interchanges $\rme$ and $\rmf$ and fixes every other element.
\[
\begin{tabular}{c | c c c c}
$B_0$  & 0 & $\rma$ & $\rme$ & $\rmf$ \\ \hline
0      & 0 & 0      & 0      & 0 \\
$\rma$ & 0 & 0      & 0      & $\rma$ \\
$\rme$ & 0 & $\rma$ & $\rme$ & 0 \\
$\rmf$ & 0 & 0      & 0      & $\rmf$ \\ \hline\hline
$x$    & 0 & $\rma$ & $\rme$ & $\rmf$ \\ \hline
$x^*$  & 0 & $\rma$ & $\rmf$ & $\rme$
\end{tabular}
\]
The {\insem} $(B_0,\uast)$ is isomorphic to the {\inv} subsemigroup of $(B_2^1,{}^\sk)$ that consists of the elements \[ \begin{bmatrix}\,0\,\,\,0\,\\\,0\,\,\,0\,\end{bmatrix}, \begin{bmatrix}\,1\,\,\,0\,\\\,0\,\,\,0\,\end{bmatrix}, \begin{bmatrix}\,0\,\,\,1\,\\\,0\,\,\,0\,\end{bmatrix}, \begin{bmatrix}\,0\,\,\,0\,\\\,0\,\,\,1\,\end{bmatrix}. \]

The {\insem} $(B_0,\uast)$ belongs to the variety $\Var(A_0,\uast)$ generated by $(A_0,\uast)$ \cite[Proposition~3.1]{Lee22} and so satisfies the identities \eqref{id: A0 basis} of $(A_0,\uast)$.
In this section, it is shown that the identities of $(B_0,\uast)$ are axiomatized by \eqref{id: A0 basis} and one additional identity.

\begin{proposition} \label{Prop: B0 basis}
The identities \eqref{id: A0 basis} and \begin{equation} x^2 y^2 \approx y^2 x^2 \label{id: B0 xxyy=yyxx} \end{equation} constitute an identity basis for $(B_0,\uast)$.
\end{proposition}

It is easily checked that $(B_0,\uast)$ satisfies the identities $\{ \eqref{id: A0 basis}, \eqref{id: B0 xxyy=yyxx} \}$.
In Section~\ref{subsec: B0 identities}, some information on identities of $(B_0,\uast)$ are given.
In Section~\ref{subsec: B0 forms}, it is shown that the identities of $(B_0,\uast)$ can be used to convert every mixed word into one of two specific forms.
Based on these results, it is shown in Section~\ref{subsec: B0 proof} that every identity of $(B_0,\uast)$ is deducible from $\{ \eqref{id: inv}, \eqref{id: A0 basis},\eqref{id: B0 xxyy=yyxx} \}$.
This completes the proof of Proposition~\ref{Prop: B0 basis}.

\begin{corollary} \label{Cor: B0 simplified}
The identities
\begin{equation} \label{id: B0 simplified}
\begin{split}
x^3 & \approx x^2, \quad xyx \approx x^2y^2, \quad x^2x^* \approx xx^*, \quad x^2yx^* \approx xyx^*, \\
& xy^*x^* \approx xyx^*, \quad xx^* \approx yy^*, \quad xyx^*z \approx z^*xyx^*
\end{split}
\end{equation}
constitute an identity basis for $(B_0,\uast)$.
\end{corollary}

\begin{proof}
It is routine to check, say with Prover9~\cite{McC05}, that the identities $\{ \eqref{id: inv}, \eqref{id: A0 basis}, \eqref{id: B0 xxyy=yyxx} \}$ and $\{ \eqref{id: inv}, \eqref{id: B0 simplified} \}$ are deducible from one another.
\end{proof}

\begin{remark}
\begin{enumerate}[\ \rm(i)]
\item The variety $\Var\,B_0$ is defined within $\Var\,A_0$ by the identity~\eqref{id: B0 xxyy=yyxx}, and $\Var\,B_0$ is the unique maximal subvariety of $\Var\,A_0$ \cite[Lemma~4.2]{Lee04}; in other words, the interval $[\Var\,B_0,\, \Var\,A_0]$ is a chain of length two.
\item In contrast, although the variety $\Var(B_0,\uast)$ is also defined within $\Var(A_0,\uast)$ by the identity~\eqref{id: B0 xxyy=yyxx} (Proposition~\ref{Prop: B0 basis}), the interval $[\Var(B_0,\uast),\Var(A_0,\uast)]$ contains an infinite descending chain \cite[Theorem~1.5]{Lee22}.
\end{enumerate}
\end{remark}

\subsection{Some identities of $(B_0,\uast)$} \label{subsec: B0 identities}

\begin{lemma} \label{Lem: id: B0 plain}
The identities $\{ \eqref{id: A0 xxx=xx}, \eqref{id: A0 xyx=yxy}, \eqref{id: B0 xxyy=yyxx} \}$ constitute an identity basis for the semigroup $B_0$.
\end{lemma}

\begin{proof}
The identities of~$A_0$, together with~\eqref{id: B0 xxyy=yyxx}, form an identity basis for~$B_0$ \cite[Section~4]{Lee04}.
The present lemma then follows from Lemma~\ref{Lem: id: A0 plain}.
\end{proof}

\begin{lemma} \label{Lem: B0 bipartite}
Let $\bfu \approx \bfv$ be any identity of $(B_0,\uast)$ such that either~$\bfu$ or~$\bfv$ is bipartite.
Then $\bfu \approx \bfv$ is deducible from $\{ \eqref{id: inv}, \eqref{id: A0 basis}, \eqref{id: B0 xxyy=yyxx} \}$.
\end{lemma}

\begin{proof}
Since $(S\ell_3,{}^\sk)$ is isomorphic to the {\inv} subsemigroup $(\{0,\rme,\rmf\},\uast)$ of $(B_0,\uast)$, the identity $\bfu \approx \bfv$ is satisfied by $(S\ell_3,{}^\sk)$.
Since either~$\bfu$ or~$\bfv$ is bipartite, by Lemma~\ref{Lem: A0 twisted}, both $\bfu$ and $\bfv$ are bipartite with $\sfcon(\bfu) = \sfcon(\bfv)$.
It follows from Lemma~\ref{Lem: A0 mixed bipartite} that $(B_0,\uast)$ satisfies the plain identity $\overline{\bfu} \approx \overline{\bfv}$.
By Lemma~\ref{Lem: id: B0 plain}, the identities $\{ \eqref{id: A0 xxx=xx}, \eqref{id: A0 xyx=yxy}, \eqref{id: B0 xxyy=yyxx} \}$ constitute an identity basis for~$B_0$, so that $\overline{\bfu} \approx \overline{\bfv}$ is deducible from $\{ \eqref{id: A0 xxx=xx}, \eqref{id: A0 xyx=yxy}, \eqref{id: B0 xxyy=yyxx} \}$.
It then follows from Lemma~\ref{Lem: A0 mixed bipartite} that $\bfu \approx \bfv$ is deducible from $\{ \eqref{id: inv}, \eqref{id: A0 basis}, \eqref{id: B0 xxyy=yyxx} \}$.
\end{proof}

An \textit{\Bblock} is a word of the form \[ \bfc = y_1^2 y_2^2 \cdots y_k^2, \] where $y_1,y_2,\ldots,y_k \in \cXX$ are such that $\overline{y_1} \order \overline{y_2} \order \cdots \order \overline{y_k}$ in~$\cX$ and $k \geq 1$.
Note that every {\Bblock} is bipartite.

\begin{lemma}
\label{Lem: B0 block}
Let $\bfw_1,\bfw_2,\ldots,\bfw_m \in \Fsem$ be any pairwise disjoint bipartite connected words such that $\sfcon(\bfw_1\bfw_2\cdots\bfw_m) = \{ y_1,y_2,\ldots,y_k\}$ and $\overline{y_1} \order \overline{y_2} \order \cdots \order \overline{y_k}$ in~$\cX$.
Then the identities $\{\eqref{id: A0 xxx=xx},\eqref{id: A0 xyx=yxy},\eqref{id: B0 xxyy=yyxx} \}$ can be used to convert the product $\bfw_1\bfw_2\cdots\bfw_m$ into the {\Bblock} $\bfc = y_1^2 y_2^2 \cdots y_k^2$.
\end{lemma}

\begin{proof}
By Lemma~\ref{Lem: A0 block}, the identities $\{\eqref{id: A0 xxx=xx},\eqref{id: A0 xyx=yxy} \}$ can be used to convert each~$\bfw_i$ into some {\Ablock} $\bfc_i$ with $\sfcon(\bfw_i) = \sfcon(\bfc_i)$.
Since $\bfc_1,\bfc_2,\ldots,\bfc_m$ are {\Ablock}s, we have $\bfc_i \stackrel{\eqref{id: A0 xxx=xx}}{\approx} \bfc_i^2$ and $\bfc_i\bfc_j \stackrel{\eqref{id: B0 xxyy=yyxx}}{\approx} \bfc_j\bfc_i$.
Hence \[ \bfw_1\bfw_2\cdots\bfw_m \stackrel{\eqref{id: A0 xxx=xx},\eqref{id: A0 xyx=yxy}}{\approx} \bfc_1\bfc_2\cdots\bfc_m \stackrel{\eqref{id: A0 xxx=xx}}{\approx} \bfc_1^2\bfc_2^2\cdots\bfc_m^2 \stackrel{\eqref{id: B0 xxyy=yyxx}}{\approx} (\bfc_1\bfc_2\cdots\bfc_m)^2. \]
Since $(\bfc_1\bfc_2\cdots\bfc_m)^2$ is a bipartite connected word with content $\{ y_1,y_2,\ldots,y_k\}$, by Lemma~\ref{Lem: A0 block}, the identities $\{\eqref{id: A0 xxx=xx},\eqref{id: A0 xyx=yxy} \}$ can be used to convert it into the {\Ablock} $(y_1 y_2 \cdots y_k)^2$.
Hence \[ (y_1 y_2 \cdots y_k)^2 \stackrel{\eqref{id: A0 xxx=xx}}{\approx} (y_1^2 y_2^2 \cdots y_k^2)^2 \stackrel{\eqref{id: B0 xxyy=yyxx}}{\approx} y_1^4 y_2^4 \cdots y_k^4 \stackrel{\eqref{id: A0 xxx=xx}}{\approx} y_1^2 y_2^2 \cdots y_k^2. \qedhere \]
\end{proof}

\subsection{Some special forms of words} \label{subsec: B0 forms}

It is easily checked that for any substitution $\varphi: \cX \to B_0$ and any variable $z \in \cX$, we have $\varphi(zz^*) = 0$ in $B_0$.
Therefore, in the $\Var(B_0,\uast)$-free algebra over~$\cX$, the class $[zz^*]$ containing $zz^*$ is its zero element.
This phenomenon is equivalent to the following result, whose justification is routine.

\begin{lemma} \label{Lem: B0 xx*y=xx*}
The identities \begin{equation} xx^*y \approx xx^*, \quad yxx^* \approx xx^*, \quad xx^* \approx yy^* \label{id: B0 xx*y=xx*} \end{equation} are deducible from $\{ \eqref{id: inv}, \eqref{id: A0 basis},\eqref{id: B0 xxyy=yyxx} \}$.
\end{lemma}

Words of other possible forms in the class $[zz^*]$ are listed in the following result.

\begin{lemma} \label{Lem: B0 zz*}
Let $\bfw \in \Fsem$.
Suppose that one of the following conditions holds\up:
\begin{enumerate}[\ \rm(a)]
\item $xx^*x \hookrightarrow \bfw$ for some $x \in \cXX$\up;
\item $xx^*yy^* \hookrightarrow \bfw$ for some $x,y \in \cXX$\up;
\item $\bfw = \bfa x \bfb x^* \bfe$ for some $x \in \cXX$ and $\bfa,\bfb,\bfe \in \Fmon$ such that for each $y \in \sfcon(\bfb)$\up, we have $\sfocc(y,\bfw) \geq 2$.
\end{enumerate}
Then the identities $\{ \eqref{id: A0 basis},\eqref{id: B0 xxyy=yyxx} \}$ can be used to convert~$\bfw$ into the word $zz^*$ for any $z \in \cXX$.
\end{lemma}

\begin{proof}
By Lemma~\ref{Lem: B0 xx*y=xx*}, it suffices to convert~$\bfw$ into the word~$zz^*$, using the identities $\{ \eqref{id: A0 basis},\eqref{id: B0 xxyy=yyxx}, \eqref{id: B0 xx*y=xx*} \}$.
If either~(a) or~(b) holds, then by Lemma~\ref{Lem: A0 zz*z}, \[ \bfw \stackrel{\eqref{id: A0 basis}}{\approx} zz^*z \stackrel{\eqref{id: B0 xx*y=xx*}}{\approx} zz^*. \]
Thus suppose~(c) holds.
By assumption, $\bfb = y_1y_2 \cdots y_m$ for some $y_1,y_2,\ldots,y_m \in \cXX$ with $m \geq 0$ such that $\sfocc(y_i,\bfw) \geq 2$ for all~$i$.
Then by Lemma~\ref{Lem: B0 xx*y=xx*}, \[ \bfw \stackrel{\eqref{id: A0 xxHx*=xHx*}}{\approx} \bfa x^2 \bfb x^* \bfe \stackrel{\eqref{id: A0 xxx=xx}}{\approx} \bfa x^2 y_1^2y_2^2 \cdots y_m^2 x^* \bfe \stackrel{\eqref{id: B0 xxyy=yyxx}}{\approx} \bfa y_1^2y_2^2 \cdots y_m^2 x^2 x^* \bfe \stackrel{\eqref{id: B0 xx*y=xx*}}{\approx} zz^*. \qedhere \]
\end{proof}

A word $\bfw \in \Fsem$ is in \textit{$B_0$-standard form} if
\begin{equation}
\bfw=\bfw_1 x \bfw_2 x^*, \label{Dis: B0-standard}
\end{equation}
where $x \in \cXX$, $\bfw_1 = x_1x_2\cdots x_m$, and $\bfw_2 = \bfs_0 \prod_{i=1}^p (\bfc_i\bfs_i)$  for some $m,p \geq 0$ such that the following conditions hold:
\begin{enumerate}[({\Bsf}1)]
\item $x_1,x_2,\ldots, x_m \in \cXX$ are such that $\overline{x_1} \order \overline{x_2} \order \cdots \order \overline{x_m} \order \overline{x}$;
\item $\bfs_0, \bfs_1, \ldots, \bfs_p \in \Fsem$ are simple and $\bfc_1, \bfc_2, \ldots, \bfc_p \in \Fsem$ are {\Bblock}s;
\item $x_1, x_2, \ldots, x_m, x, \bfs_0, \bfs_1, \ldots, \bfs_p, \bfc_1, \bfc_2, \ldots, \bfc_p$ are pairwise disjoint;
\item either
\begin{enumerate}[\ \rm(a)]
\item $p = 0$ with $\bfw_2 = \bfs_0$ and $\bfs_0 \in \cX$; or \label{Bsf4 singleton}
\item $\overline{\sfh(\bfw_2)} \order \overline{\sft(\bfw_2)}$. \label{Bsf4 s0 orderedsquare}
\end{enumerate}
\end{enumerate}

\begin{remark} \label{Rmk: B0-standard}
The following holds for the word~$\bfw$ in~\eqref{Dis: B0-standard} in $B_0$-standard form:
\begin{enumerate}[\ \rm(i)]
\item If $m=0$, then $\bfw_1 = 1$.
\item If $p = 0$, then $\bfw_2 = \bfs_0 \in \Fsem$; in particular, $\bfw_2$ always contains some simple variable and so is nonempty.
\item $\{x, x^*\}$ is the only mixed pair of~$\bfw$ and $x,x^* \notin \sfcon(\bfw_1\bfw_2)$.
\item $\bfw_1$ and $\bfw_2$ are bipartite words such that $\sfcon(\overline{\bfw_1}) \cap \sfcon(\overline{\bfw_2}) = \emptyset$.
\item Each variable in $\cX$ occurs at most twice in $\overline{\bfw}$.
\end{enumerate}
\end{remark}

\begin{lemma} \label{Lem: B0 substitution}
Let $\bfw=\bfw_1 x \bfw_2 x^*$ be the word in~\eqref{Dis: B0-standard} in $B_0$-standard form and $z \in \cXX$ be any simple variable in~$\bfw$\up, so that $z \in \sfcon(\bfs_0 \bfs_1 \cdots \bfs_p)$ and $\bfw_2=\bfa z \bfb$ for some $\bfa, \bfb \in \Fmon$.
Then there exists a substitution $\gamma_\bfw^z: \cX \to B_0$ such that
\begin{enumerate}[\ \rm(i)]
\item $\gamma_\bfw^z(\bfw_1 x \bfa) = \rme$\up, $\gamma_\bfw^z(z) = \rma$\up,  and $\gamma_\bfw^z(\bfb x^*) = \rmf$\up, so that $\gamma_\bfw^z(\bfw) = \rma$\up;
\item $\gamma_\bfw^z(s) = 0$ for all $s \in \cX$ such that $s \notin \sfcon(\overline{\bfw})$.
\end{enumerate}
\end{lemma}

\begin{proof}
It follows from Remark~\ref{Rmk: B0-standard}(iii),(iv) that $\sfcon(\bfw_1) = \mathcal{H}_1 \cup \mathcal{K}_1^*$, $\sfcon(\bfa) = \mathcal{H}_2 \cup \mathcal{K}_2^*$, and $\sfcon(\bfb) = \mathcal{H}_3 \cup \mathcal{K}_3^*$ for some $\mathcal{H}_1,\mathcal{H}_2,\mathcal{H}_3, \mathcal{K}_1,\mathcal{K}_2, \mathcal{K}_3 \subseteq \cX$ such that $\mathcal{H}_1, \mathcal{H}_2, \mathcal{H}_3, \mathcal{K}_1, \mathcal{K}_2, \mathcal{K}_3, \{ x,x^* \}, \{ z \}$ are pairwise disjoint sets.
By symmetry, it suffices to assume that $x \in \cX$, so that $\sfcon(\overline{\bfw}) = \mathcal{H}_1 \cup \mathcal{K}_1 \cup \mathcal{H}_2 \cup \mathcal{K}_2 \cup \mathcal{H}_3 \cup \mathcal{K}_3 \cup \{x, \overline{z} \}$.
Define
\[
\gamma_\bfw^z(s) = \begin{cases}
\rme & \text{if $s \in \mathcal{H}_1 \cup \mathcal{H}_2 \cup \mathcal{K}_3 \cup \{ x \}$}, \\
\rma & \text{if $s = \overline{z}$}, \\
\rmf & \text{if $s \in \mathcal{K}_1 \cup \mathcal{K}_2 \cup \mathcal{H}_3$}, \\
0 & \text{otherwise}.
\end{cases}
\]
Then it is routinely checked that the substitution $\gamma_\bfw^z$ satisfies~(i) and~(ii).
\end{proof}

\begin{corollary} \label{C: B0 violate}
For any word~$\bfw$ in $B_0$-standard form and any $z \in \cXX$\up, the identity $\bfw \approx zz^*$ is not satisfied by $(B_0,\uast)$.
\end{corollary}

\begin{proof}
Let $\bfw = \bfw_1 x \bfw_2 x^*$ be the word in \eqref{Dis: B0-standard} in $B_0$-standard form.
Then by Remark~\ref{Rmk: B0-standard}(ii), the word~$\bfw_2$ contains some simple variable $s \in \cXX$, so that $\bfw_2 = \bfa s \bfb$ for some $\bfa,\bfb \in \Fmon$.
Under the substitution $\gamma_\bfw^s : \cX \to B_0$ in Lemma~\ref{Lem: B0 substitution}, we have $\gamma_\bfw^s(\bfw) = \rma$ and $\gamma_\bfw^s(zz^*) = 0$.
\end{proof}

\begin{lemma} \label{Lem: B0 standard form}
Let $\bfw$ be any mixed word.
Then the identities $\{ \eqref{id: inv}, \eqref{id: A0 basis}, \eqref{id: B0 xxyy=yyxx} \}$ can be used to convert~$\bfw$ into exactly one of the following\up:
\begin{enumerate}[\ \rm(i)]
\item the word $zz^*$ for any $z \in \cXX$\up;
\item some word in $B_0$-standard form.
\end{enumerate}
\end{lemma}

\begin{proof}
By Lemma~\ref{Lem: B0 xx*y=xx*}, it suffices to convert~$\bfw$ into~(i) or~(ii), using the identities $\{ \eqref{id: inv}, \eqref{id: A0 basis},\eqref{id: B0 xxyy=yyxx}, \eqref{id: B0 xx*y=xx*} \}$.
By Lemma~\ref{Lem: A0 standard form}, the identities $\{ \eqref{id: inv}, \eqref{id: A0 basis} \}$ can first be used to convert~$\bfw$ into either $zz^*z$ or some word in $A_0$-standard form.
In the former case, the first identity in \eqref{id: B0 xx*y=xx*} can be used to convert $zz^*z$ into $zz^*$.
Therefore, it remains to assume that $\bfw = \bfw_1 x \bfw_2 x^*$, where $\bfw_1 = x_1x_2 \cdots x_m$ and $\bfw_2 = \bfs_0 \prod_{i=1}^p (\bfc_i\bfs_i)$, satisfies conditions (\Asf1)--(\Asf4).
Then conditions (\Bsf1) and (\Bsf3) hold because they coincide with conditions (\Asf1) and (\Asf3).

By condition (\Asf2), $\bfs_0, \bfs_1, \ldots, \bfs_p \in \Fmon$ are simple and $\bfc_1, \bfc_2, \ldots, \bfc_p \in \Fsem$ are {\Ablock}s.
By Lemma~\ref{Lem: B0 block}, each~$\bfc_i$ can be converted by $\{ \eqref{id: A0 basis},\eqref{id: B0 xxyy=yyxx} \}$ into some {\Bblock} $y_{i,1}^2 y_{i,2}^2 \cdots y_{i,h_i}^2$.
If $\bfs_0 = \bfs_1 = \cdots = \bfs_p = 1$, then by Lemma~\ref{Lem: B0 zz*}, the identities $\{ \eqref{id: inv}, \eqref{id: A0 basis} \}$ can be used to convert~$\bfw$ into $zz^*$.
Therefore, assume that $\bfs_0, \bfs_1, \ldots, \bfs_p$ are not all empty.
If $\bfs_i = 1$ for some $i \in \{ 1,2,\ldots,p-1\}$, so that the {\Bblock}s $\bfc_i$ and $\bfc_{i+1}$ are adjacent, then the identity~\eqref{id: B0 xxyy=yyxx} can be used to arrange the squares $y_{i,1}^2,y_{i,2}^2,\ldots,y_{i,h_i}^2,y_{i+1,1}^2,y_{i+1,2}^2,\ldots,y_{i+1,h_{i+1}}^2$ in the product $\bfc_i\bfc_{i+1}$ in order, resulting in a single {\Bblock}.
Hence we may assume that for each $i \in \{ 1,2,\ldots,p-1\}$, the words $\bfc_i$ and $\bfc_{i+1}$ are separated due to $\bfs_i \neq 1$.
If $\bfs_0 = 1$, so that $x$ is adjacent to the {\Bblock} $\bfc_1$, then the identities $\{ \eqref{id: A0 basis},\eqref{id: B0 xxyy=yyxx} \}$ can be used to move~$\bfc_1$ to the left of~$x$ and turn it into a simple word:
\begin{align*}
\bfw & \makebox[0.37in]{$\stackrel{\eqref{id: A0 xxHx*=xHx*}}{\approx}$} x_1x_2 \cdots x_m \cdot x^2 \cdot \overbrace{y_{1,1}^2 y_{1,2}^2 \cdots y_{1,h_1}^2}^{\bfc_1} \cdot\, \bfs_1 \bigg(\prod_{i=2}^p (\bfc_i\bfs_i) \bigg) x^* \\
& \makebox[0.37in]{$\stackrel{\eqref{id: B0 xxyy=yyxx}}{\approx}$} x_1x_2 \cdots x_m \cdot y_{1,1}^2 y_{1,2}^2 \cdots y_{1,h_1}^2 \cdot x^2 \cdot \bfs_1 \bigg(\prod_{i=2}^p (\bfc_i\bfs_i) \bigg) x^* \\
& \makebox[0.37in]{$\stackrel{\eqref{id: A0 xxHyTy*=xHyTy*}}{\approx}$} x_1x_2 \cdots x_m \cdot y_{1,1} y_{1,2} \cdots y_{1,h_1} \cdot x^2 \cdot \bfs_1 \bigg(\prod_{i=2}^p (\bfc_i\bfs_i) \bigg) x^* \\
& \makebox[0.37in]{$\stackrel{\eqref{id: A0 xxHx*=xHx*}}{\approx}$} x_1x_2 \cdots x_m \cdot y_{1,1} y_{1,2} \cdots y_{1,h_1} \cdot x \cdot \bfs_1 \bigg(\prod_{i=2}^p (\bfc_i\bfs_i) \bigg) x^*;
\end{align*}
by the arguments in the proof of Lemma~\ref{Lem: A0 standard form}, we may assume that \[ \overline{x_1} \order \overline{x_2} \order \cdots \order \overline{x_m} \order \overline{y_{1,1}} \order \overline{y_{1,2}} \order \cdots \order \overline{y_{1,h_1}} \order \overline{x}. \]
If $\bfs_p = 1$, so that $x^*$ is adjacent to the {\Bblock} $\bfc_p$, then the identities $\{ \eqref{id: A0 basis},\eqref{id: B0 xxyy=yyxx} \}$ can be used to move~$\bfc_p$ to the left of~$x$ and turn it into a simple word:
\begin{align*}
\bfw & \makebox[0.37in]{$\stackrel{\eqref{id: A0 xxHx*=xHx*}}{\approx}$} x_1x_2 \cdots x_m \cdot x \cdot \bfs_0 \bigg(\prod_{i=1}^{p-1} (\bfc_i\bfs_i) \bigg) \overbrace{y_{p,1}^2 y_{p,2}^2 \cdots y_{p,h_p}^2}^{\bfc_p} \cdot\, (x^*)^2 \\
& \makebox[0.37in]{$\stackrel{\eqref{id: B0 xxyy=yyxx}}{\approx}$} x_1x_2 \cdots x_m \cdot x \cdot \bfs_0 \bigg(\prod_{i=1}^{p-1} (\bfc_i\bfs_i) \bigg) (x^*)^2 y_{p,1}^2 y_{p,2}^2 \cdots y_{p,h_p}^2 \\
& \makebox[0.37in]{$\stackrel{\eqref{id: A0 xHx*y=y*xHx*}}{\approx}$} x_1x_2 \cdots x_m \cdot (y_{p,1}^*)^2 (y_{p,2}^*)^2 \cdots (y_{p,h_p}^*)^2 \cdot x \cdot \bfs_0 \bigg(\prod_{i=1}^{p-1} (\bfc_i\bfs_i) \bigg) (x^*)^2  \\
& \makebox[0.37in]{$\stackrel{\eqref{id: A0 xxHyTy*=xHyTy*}}{\approx}$} x_1x_2 \cdots x_m \cdot y_{p,1}^* y_{p,2}^* \cdots y_{p,h_p}^* \cdot x \cdot \bfs_0 \bigg(\prod_{i=1}^{p-1} (\bfc_i\bfs_i) \bigg) (x^*)^2 \\
& \makebox[0.37in]{$\stackrel{\eqref{id: A0 xxHx*=xHx*}}{\approx}$} x_1x_2 \cdots x_m \cdot y_{p,1}^* y_{p,2}^* \cdots y_{p,h_p}^* \cdot x \cdot \bfs_0 \bigg(\prod_{i=1}^{p-1} (\bfc_i\bfs_i) \bigg) x^*;
\end{align*}
by repeating the arguments in the proof of Lemma~\ref{Lem: A0 standard form}, we may assume that \[ \overline{x_1} \order \overline{x_2} \order \cdots \order \overline{x_m} \order \overline{y_{p,1}} \order \overline{y_{p,2}} \order \cdots \order \overline{y_{p,h_p}} \order \overline{x}. \]
Therefore, we may assume that $\bfs_0,\bfs_p \neq 1$.
It follows that $\bfs_0, \bfs_1, \ldots, \bfs_p \in \Fsem$, so that condition (\Bsf2) is satisfied.

It remains to address condition (\Bsf4).
If~$\bfw_2$ is a single variable, so that $p=0$ with $\bfw_2 = \bfs_0 \in \cXX$, then the identity~\eqref{id: A0 xyx*=xy*x*} can be used to convert~$\bfs_0$ into a variable in~$\cX$, whence condition (\Bsf\ref{Bsf4 singleton}) is satisfied.
Hence, assume that~$\bfw_2$ is not a single variable, so that $\overline{\sfh(\bfw_2)} \neq \overline{\sft(\bfw_2)}$ by conditions (\Bsf2) and (\Bsf3).
In this case, since each~$\bfs_i$ is a nonempty simple word, we have $\bfs_i = s_{i,1} s_{i,2} \cdots s_{i,k_i}$ for some $s_{i,1}, s_{i,2}, \ldots, s_{i,k_i} \in \cXX$ such that $\overline{s_{i,1}}, \overline{s_{i,2}}, \ldots, \overline{s_{i,k_i}}$ are distinct.
If $\overline{\sfh(\bfw_2)} \order \overline{\sft(\bfw_2)}$, then condition (\Bsf\ref{Bsf4 s0 orderedsquare}) is satisfied.
If $\overline{\sfh(\bfw_2)} \norder \overline{\sft(\bfw_2)}$, so that $\overline{s_{p,k_p}} \order \overline{s_{0,1}}$, then
\[ \bfw \stackrel{\eqref{id: A0 xyx*=xy*x*}}{\approx} \bfw_1 \cdot x \cdot \bigg(\bfs_0\prod_{i=1}^p (\bfc_i\bfs_i) \bigg)^* x^* \stackrel{\eqref{id: inv}}{\approx} \bfw_1 \cdot x \cdot \bigg(\prod_{i=p}^1 (\bfs_i^*\bfc_i^*) \bigg) \bfs_0^* \cdot x^*, \]
where the identities $\{ \eqref{id: inv}, \eqref{id: B0 xxyy=yyxx} \}$ can be used to convert $\bfs_i^*$ and $\bfc_i^*$ into the simple word $s_{i,k_i}^* s_{i,k_i-1}^* \cdots s_{i,1}^*$ and the {\Bblock} $(y_{i,1}^*)^2 (y_{i,2}^*)^2 \cdots (y_{i,h_i}^*)^2$, respectively; thus, condition (\Bsf\ref{Bsf4 s0 orderedsquare}) is satisfied.

Consequently, the identities $\{ \eqref{id: inv}, \eqref{id: A0 basis}, \eqref{id: B0 xxyy=yyxx} \}$ can be used to convert~$\bfw$ into either $zz^*$ or some word~$\widetilde{\bfw}$ in $B_0$-standard form.
But if the identities $\{ \eqref{id: inv}, \eqref{id: A0 basis}, \eqref{id: B0 xxyy=yyxx} \}$ can be used to convert~$\bfw$ into both $zz^*$ and~$\widetilde{\bfw}$, then that would imply that $(B_0,\uast)$ satisfies the identity $\widetilde{\bfw} \approx zz^*$, which is impossible by Corollary~\ref{C: B0 violate}.
\end{proof}

\subsection{Proof of Proposition~\ref{Prop: B0 basis}} \label{subsec: B0 proof}

Consider any identity \[ \bfu \approx \bfv \] satisfied by $(B_0, \uast)$.
It suffices to show that $\bfu \approx \bfv$ is deducible from $\{ \eqref{id: inv}, \eqref{id: A0 basis}, \eqref{id: B0 xxyy=yyxx} \}$.
By Lemma~\ref{Lem: B0 bipartite}, this result holds if either~$\bfu$ or~$\bfv$ is bipartite.
Therefore, suppose that~$\bfu$ and~$\bfv$ are both mixed.
By Corollary~\ref{C: B0 violate} and Lemma~\ref{Lem: B0 standard form}, the identities $\{ \eqref{id: inv}, \eqref{id: A0 basis}, \eqref{id: B0 xxyy=yyxx} \}$ can be used to convert~$\bfu$ and~$\bfv$ simultaneously to either $zz^*$ or words in $B_0$-standard form.
In the former case, $\bfu \approx \bfv$ is deducible from $\{ \eqref{id: inv}, \eqref{id: A0 basis}, \eqref{id: B0 xxyy=yyxx} \}$, whence the proof is complete.
Therefore, it remains to consider the latter case, whence we may assume that~$\bfu$ and~$\bfv$ are in $B_0$-standard form, say \[ \bfu = \bfu_1 x \bfu_2 x^* \quad \text{and} \quad \bfv = \bfv_1 y \bfv_2 y^*, \] where $x,y \in \cXX$, $\bfu_1 = x_1x_2 \cdots x_m$, $\bfu_2 = \bfs_0\prod_{i=1}^p (\bfc_i\bfs_i)$, $\bfv_1 = y_1y_2 \cdots y_n$, and $\bfv_2 = \bft_0\prod_{i=1}^q (\bfd_i\bft_i)$ satisfy conditions (\Bsf1)--(\Bsf4).

\begin{lemma}\label{Lem: B0 content}
The following holds for the words~$\bfu$ and~$\bfv$\up:
\begin{enumerate}[\ \rm(i)]
\item $\sfcon(\overline{\bfu}) = \sfcon(\overline{\bfv})$\up;
\item $\bfu_1 x =\bfv_1 y$\up;
\item $\sfcon(\overline{\bfu_2}) = \sfcon(\overline{\bfv_2})$.
\end{enumerate}
\end{lemma}

\begin{proof}
(i) Suppose that $\sfcon(\overline{\bfu}) \neq \sfcon(\overline{\bfv})$, say there exists a variable $t \in \sfcon(\overline{\bfv})$ such that $t \notin \sfcon(\overline{\bfu})$.
Then by Remark~\ref{Rmk: B0-standard}(ii), the word $\bfu_2$ contains some simple variable $z \in \cXX$, so that $\bfu_2 = \bfa z \bfb$ for some $\bfa,\bfb \in \Fmon$.
Under the substitution $\gamma_\bfu^z : \cX \to B_0$ in Lemma~\ref{Lem: B0 substitution}, the contradiction $\gamma_\bfu^z(\bfu) = \rma \neq 0 = \gamma_\bfu^z(\bfv)$ is deduced.

\medskip

(ii) Due to condition (\Bsf1), the equality $\bfu_1 x =\bfv_1 y$ follows from $\sfcon(\bfu_1 x) = \sfcon(\bfv_1 y)$; to establish the latter, by symmetry, it suffices to verify the inclusion $\sfcon(\bfu_1 x) \subseteq \sfcon(\bfv_1 y)$.
To this end, we need to first show that $y \in \sfcon(\bfu_1x)$.
Since $\overline{y} \in \sfcon(\overline{\bfv}) = \sfcon(\overline{\bfu})$ by part~(i),
\begin{enumerate}[\ (a)]
\item[(a)] either $y \in \sfcon(\bfu)$ or $y^* \in \sfcon(\bfu)$.
\end{enumerate}
Now since $\bfs_p \neq 1$, the variable $z = \sft(\bfs_p) = \sft(\bfu_2)$ is simple in~$\bfu$, so that $\bfu_2 = \bfa z$ for some $\bfa \in \Fmon$ such that $\overline{z} \notin \sfcon(\overline{\bfa}) $.
Then under the substitution $\gamma_\bfu^z : \cX \to B_0$ in Lemma~\ref{Lem: B0 substitution}, we have \[ \gamma_\bfu^z(\bfu) = \gamma_\bfu^z(\bfu_1 x \bfa) \cdot \gamma_\bfu^z(z) \cdot \gamma_\bfu^z(x^*) = \rme \cdot \rma \cdot \rme^* = \rma. \]
If $y^* \in \sfcon(\bfu_1 x \bfu_2)$, then
\begin{align*}
\gamma_\bfu^z(\bfv) & = \gamma_\bfu^z(\bfv_1) \cdot \gamma_\bfu^z(y) \cdot \gamma_\bfu^z(\bfv_2) \cdot \gamma_\bfu^z(y^*) \\
& = \begin{cases} \gamma_\bfu^z(\bfv_1) \cdot \rma \cdot \gamma_\bfu^z(\bfv_2) \cdot \rma & \text{if $y^* = z$} \\ \gamma_\bfu^z(\bfv_1) \cdot \rmf \cdot \gamma_\bfu^z(\bfv_2) \cdot \rme & \text{if $y^* \neq z$} \end{cases} \\
& = 0,
\end{align*}
which is impossible.
Therefore,
\begin{enumerate}[\ (a)]
\item[(b)] $y^* \notin \sfcon(\bfu_1 x \bfu_2)$, which implies that $y \neq x^*$.
\end{enumerate}
Note that if $y \neq x$, then together with~(b), we have $y^* \notin \sfcon(\bfu_1x\bfu_2x^*) = \sfcon(\bfu)$, so that $y \in \sfcon(\bfu_1x\bfu_2)$ by~(a) and~(b).
On the other hand, if $y = x$, then clearly $y \in \sfcon(\bfu_1x\bfu_2)$.
Therefore, $y \in \sfcon(\bfu_1x\bfu_2)$ either way.

Seeking a contradiction, suppose that $y \in \sfcon(\bfu_2)$.
Then $\sfocc(y,\bfu_2) \in \{1,2\}$ by condition (\Bsf2).
If $\sfocc(y,\bfu_2)=1$, then under the substitution $\gamma_\bfu^y : \cX \to B_0$ in Lemma~\ref{Lem: B0 substitution}, we have $\gamma_\bfu^y (\bfu) = \rma$ and \[ \gamma_\bfu^y(\bfv) = \gamma_\bfu^y(\bfv_1) \cdot \gamma_\bfu^y(y) \cdot \gamma_\bfu^y(\bfv_2) \cdot \gamma_\bfu^y(y^*)= \gamma_\bfu^y(\bfv_1) \cdot \rma \cdot \gamma_\bfu^y(\bfv_2) \cdot \rma = 0, \] which is impossible.
If $\sfocc(y,\bfu_2)=2$, so that $\bfu_2 = h \bfb$ for some $\bfb \in \Fsem$ with $y \in \sfcon(\bfb)$ and $h = \sfh(\bfs_0)$ being simple in~$\bfu_2$, then under the substitution $\gamma_\bfu^h : \cX \to B_0$ in Lemma~\ref{Lem: B0 substitution}, we have $\gamma_\bfu^h (\bfu) = \rma$ and \[ \gamma_\bfu^h(\bfv) = \gamma_\bfu^h(\bfv_1) \cdot \gamma_\bfu^h(y) \cdot \gamma_\bfu^h(\bfv_2) \cdot \gamma_\bfu^h(y^*) = \gamma_\bfu^h(\bfv_1) \cdot \rmf \cdot \gamma_\bfu^h(\bfv_2) \cdot \rme = 0, \] which again is impossible.
Thus, $y \notin \sfcon(\bfu_2)$; but since $y \in \sfcon(\bfu_1x\bfu_2)$, together with~(b), we have
\begin{enumerate}[\ (a)]
\item[(c)] $y \in \sfcon(\bfu_1 x)$ and $y,y^* \notin \sfcon(\bfu_2)$.
\end{enumerate}
By a symmetrical argument,
\begin{enumerate}[\ (a)]
\item[(d)] $x \in \sfcon(\bfv_1 y)$ and $x,x^* \notin \sfcon(\bfv_2)$.
\end{enumerate}

Now we are ready to establish the inclusion $\sfcon(\bfu_1 x) \subseteq \sfcon(\bfv_1 y)$.
Suppose there exists some variable $z \in \sfcon(\bfu_1x)$ such that $z \notin \sfcon(\bfv_1y)$.
Then clearly $z \neq y$.
But if $z=y^*$, then it follows from~(c) that $y,y^* \in \sfcon(\bfu_1 x)$, whence condition (\Bsf1) is contradicted.
Hence
\begin{enumerate}[\ (a)]
\item[(e)] $z \notin \{ y, y^*\}$.
\end{enumerate}
Since $\overline{z} \in \sfcon(\overline{\bfu}) = \sfcon(\overline{\bfv})$ by part~(i), it follows from~(e) that either $z \in \sfcon(\bfv_1\bfv_2)$ or $z^* \in \sfcon(\bfv_1\bfv_2)$.
But since $z \notin \sfcon(\bfv_1)$ by assumption, we have $z \in \sfcon(\bfv_2)$ or $z^* \in \sfcon(\bfv_1)$ or $z^* \in \sfcon(\bfv_2)$.
These three cases are shown in the following to be impossible.
Therefore, the variable~$z$ does not exist, whence the required inclusion $\sfcon(\bfu_1x) \subseteq \sfcon(\bfv_1y)$ is established.

\medskip

\noindent{\sc Case~1:} $z \in \sfcon(\bfv_2)$.
Then $z \notin \{ x,x^*\}$ by~(d).
But since $z \in \sfcon(\bfu_1x)$ by assumption, we have
\begin{enumerate}[\ (a)]
\item[(f)] $z \in \sfcon(\bfu_1)$.
\end{enumerate}
By condition (\Bsf2), we have $\sfocc(z,\bfv_2) \in \{1,2\}$, so there are two subcases.

\smallskip

\noindent{\sc Subcase~1.1:} $\sfocc(z, \bfv_2) = 1$.
Then $\bfv_2=\bfa z \bfb$ for some $\bfa, \bfb \in \Fmon$ such that $\overline{z} \notin \sfcon(\overline{\bfa\bfb})$.
Hence under the substitution $\gamma_\bfv^z : \cX \to B_0$ in Lemma~\ref{Lem: B0 substitution}, we have \[ \gamma_\bfv^z(\bfv) =\gamma_\bfv^z(\bfv_1 y \bfa) \cdot \gamma_\bfv^z(z) \cdot \gamma_\bfv^z(\bfb y^*)=\rme \cdot \rma \cdot \rmf=\rma. \]
It follows that $\gamma_\bfv^z(x) = \rme$ because $x \in \sfcon(\bfv_1 y)$ by~(d), and $\gamma_\bfv^z(\bfu_1) = \cdots \rma \cdots$ because $z \in \sfcon(\bfu_1)$ by~(f).
Therefore, \[ \gamma_\bfv^z(\bfu) = \gamma_\bfv^z (\bfu_1) \cdot \gamma_\bfv^z(x) \cdot \gamma_\bfv^z(\bfu_2x^*) = \cdots \rma \cdots \rme  \cdot \gamma_\bfv^z(\bfu_2x^*) = 0, \] which implies a contradiction.

\smallskip

\noindent{\sc Subcase~1.2:} $\sfocc(z, \bfv_2) = 2$.
Since $\bft_0 \neq 1$, the variable $h = \sfh(\bft_0) = \sfh(\bfv_2)$ is simple in~$\bfv$, so that $\bfv_2 = h\bfa z^2\bfb$ for some $\bfa,\bfb \in \Fmon$ such that $\bfa,\bfb, h, z$ are pairwise disjoint.
Then under the substitution $\gamma_\bfv^h : \cX \to B_0$ in Lemma~\ref{Lem: B0 substitution}, we have \[ \gamma_\bfv^h(\bfv) = \gamma_\bfv^h(\bfv_1 y) \cdot \gamma_\bfv^h(h) \cdot \gamma_\bfv^h(\bfa z^2 \bfb y^*) = \rme \cdot \rma \cdot \rmf = \rma. \]
It follows that $\gamma_\bfv^h(x) = \rme$ because $x \in \sfcon(\bfv_1 y)$ by~(d), and $\gamma_\bfv^h(z) = \rmf$.
Further, since $z \in \sfcon(\bfu_1)$ by~(f), we have $\gamma_\bfv^h(\bfu_1) = \cdots \rmf \cdots$.
Therefore, \[ \gamma_\bfv^h(\bfu) = \gamma_\bfv^h(\bfu_1) \cdot \gamma_\bfv^h(x) \cdot \gamma_\bfv^h(\bfu_2x^*) = \cdots \rmf \cdots \rme \cdot \gamma_\bfv^h(\bfu_2x^*) = 0, \] which implies a contradiction.

\medskip

\noindent{\sc Case~2:} $z^* \in \sfcon(\bfv_1)$.
Since $\bfs_0 \neq 1$, the variable $h = \sfh(\bfs_0) = \sfh(\bfu_2)$ is simple in~$\bfu$, so that $\bfu_2 = h\bfb$ for some $\bfb \in \Fmon$ such that $\overline{h} \notin \sfcon(\overline{\bfb})$.
Then under the substitution $\gamma_\bfu^h : \cX \to B_0$ in Lemma~\ref{Lem: B0 substitution}, we have \[ \gamma_\bfu^h(\bfu) = \gamma_\bfu^h(\bfu_1 x) \cdot \gamma_\bfu^h(h) \cdot \gamma_\bfu^h(\bfb x^*) = \rme \cdot \rma \cdot \rmf = \rma. \]
It follows that $\gamma_\bfu^h(y) = \rme$ because $y \in \sfcon(\bfu_1x)$ by~(c), and that $\gamma_\bfu^h(z) = \rme$ because $z \in \sfcon(\bfu_1x)$ by assumption.
Further, since $z^* \in \sfcon(\bfv_1)$, we have $\gamma_\bfu^h(\bfv_1) = \cdots \rmf \cdots$.
Therefore, \[ \gamma_\bfu^h(\bfv) = \gamma_\bfu^h(\bfv_1) \cdot \gamma_\bfu^h(y) \cdot \gamma_\bfu^h(\bfv_2y^*) = \cdots \rmf \cdots \rme \cdot \gamma_\bfu^h(\bfv_2y^*) = 0, \] which implies a contradiction.

\medskip

\noindent{\sc Case~3:} $z^* \in \sfcon(\bfv_2)$.
Then $z \notin \{ x,x^*\}$ by~(d).
But since $z \in \sfcon(\bfu_1x)$ by assumption, we have
\begin{enumerate}[\ (a)]
\item[(g)] $z \in \sfcon(\bfu_1)$.
\end{enumerate}
By condition (\Bsf2), we have $\sfocc(z^*,\bfv_2) \in \{1,2\}$, so there are two subcases.

\smallskip

\noindent{\sc Subcase~3.1:} $\sfocc(z^*,\bfv_2) = 1$.
Then $\bfv_2=\bfa z^* \bfb$ for some $\bfa, \bfb \in \Fmon$ such that $\overline{z} \notin \sfcon(\overline{\bfa\bfb})$.
Hence under the substitution $\gamma_\bfv^{z^*} : \cX \to B_0$ in Lemma~\ref{Lem: B0 substitution}, we have \[ \gamma_\bfv^{z^*}\!(\bfv) =\gamma_\bfv^{z^*}\!(\bfv_1 y \bfa) \cdot \gamma_\bfv^{z^*}\!(z^*) \cdot \gamma_\bfv^{z^*}\!(\bfb y^*)=\rme \cdot \rma \cdot \rmf=\rma. \]
It follows that $\gamma_\bfv^{z^*}\!(x) = \rme$ because $x \in \sfcon(\bfv_1 y)$ by~(d), and $\gamma_\bfv^{z^*}\!(\bfu_1) = \cdots \rma \cdots$ because $z \in \sfcon(\bfu_1)$ by~(g).
Therefore, \[ \gamma_\bfv^{z^*}\!(\bfu) = \gamma_\bfv^{z^*}\!(\bfu_1) \cdot \gamma_\bfv^{z^*}\!(x) \cdot \gamma_\bfv^{z^*}\!(\bfu_2x^*) = \cdots \rma \cdots \rme  \cdot \gamma_\bfv^{z^*}\!(\bfu_2x^*) = 0, \] which implies a contradiction.

\smallskip

\noindent{\sc Subcase~3.2:} $\sfocc(z^*,\bfv_2) = 2$.
Since $\bft_q \neq 1$, the variable $t = \sft(\bft_q) = \sft(\bfv_2)$ is simple in~$\bfv$, so that $\bfv_2 = \bfa (z^*)^2\bfb t$ for some $\bfa,\bfb \in \Fmon$ such that $\bfa,\bfb, z, t$ are pairwise disjoint.
Then under the substitution $\gamma_\bfv^t : \cX \to B_0$ in Lemma~\ref{Lem: B0 substitution}, we have \[ \gamma_\bfv^t(\bfv) = \gamma_\bfv^t\big(\bfv_1 y \bfa (z^*)^2 \bfb\big) \cdot \gamma_\bfv^t(t) \cdot \gamma_\bfv^t(y^*) = \rme \cdot \rma \cdot \rmf = \rma. \]
It follows that $\gamma_\bfv^t(x) = \rme$ because $x \in \sfcon(\bfv_1 y)$ by~(d), and $\gamma_\bfv^t(z) = \rmf$.
Further, since $z \in \sfcon(\bfu_1)$ by~(g), we have $\gamma_\bfv^t(\bfu_1) = \cdots \rmf \cdots$.
Therefore, \[ \gamma_\bfv^t(\bfu) = \gamma_\bfv^t(\bfu_1) \cdot \gamma_\bfv^t(x) \cdot \gamma_\bfv^t(\bfu_2x^*) = \cdots \rmf \cdots \rme \cdot \gamma_\bfv^h(\bfu_2x^*) = 0, \] which implies a contradiction.

\medskip

(iii) This is a consequence of parts~(i) and~(ii).
\end{proof}

Therefore, by Lemma~\ref{Lem: B0 content}, we now have \[ \bfu = { \underbrace{x_1x_2 \cdots x_m}_{\bfu_1} } \cdot x \cdot { \underbrace{\bfs_0\prod_{i=1}^p (\bfc_i\bfs_i)}_{\bfu_2} } \cdot x^* \quad \text{and} \quad \bfv = { \underbrace{x_1x_2 \cdots x_m}_{\bfu_1} } \cdot x \cdot { \underbrace{\bft_0\prod_{i=1}^q (\bfd_i\bft_i)}_{\bfv_2} } \cdot x^*, \] where conditions (\Bsf1)--(\Bsf4) are satisfied.
In the remainder of this section, it is shown that $\bfu_2 = \bfv_2$ (Lemma~\ref{Lem: B0 u2=v2}), so that $\bfu = \bfv$.
The identity $\bfu \approx \bfv$ is thus vacuously deducible from $\{ \eqref{id: inv}, \eqref{id: A0 basis}, \eqref{id: B0 xxyy=yyxx} \}$, whence the proof of Proposition~\ref{Prop: B0 basis} is complete.

\begin{lemma} \label{Lem: B0 con s c}
\begin{enumerate}[\ \rm(i)]
\item $\sfcon(\overline{\bfs_0\bfs_1\cdots \bfs_p})=\sfcon(\overline{\bft_0 \bft_1 \cdots \bft_q})$.
\item $\sfcon(\overline{\bfc_1\bfc_2\cdots \bfc_p})=\sfcon(\overline{\bfd_1 \bfd_2 \cdots \bfd_q})$.
\end{enumerate}
\end{lemma}

\begin{proof}
(i) Let $\bfs=\bfs_0 \bfs_1 \cdots \bfs_p$ and $\bft=\bft_0 \bft_1 \cdots \bft_q$.
Suppose that $\sfcon(\overline{\bfs}) \nsubseteq \sfcon(\overline{\bft})$.
Then there exists some $z \in \sfcon(\bfs)$ such that $\overline{z} \in \sfcon(\overline{\bfs})$ and $\overline{z} \notin \sfcon(\overline{\bft})$.
Therefore, under the substitution $\gamma_\bfu^z: \cX \to B_0$ in Lemma~\ref{Lem: B0 substitution}, we have $\gamma_\bfu^z(\bfu) = \rma$.
On the other hand, since $\overline{z} \in \sfcon(\overline{\bfu_2}) = \sfcon(\overline{\bfv_2})$ by Lemma~\ref{Lem: B0 content}(iii) but $\overline{z} \notin \sfcon(\overline{\bft})$ by assumption, we have $\overline{z} \in \sfcon(\overline{\bfd_i})$ for some $i \in \{ 1,2,\ldots,q\}$.
Since $\bfd_i$ is an {\Bblock}, we have $\gamma_\bfu^z(\bfd_i) = \cdots \rma^2 \cdots = 0$, whence the contradiction $\gamma_\bfu^z(\bfv) = 0$ is deduced.
Hence the variable~$z$ does not exist, so that the inclusion $\sfcon(\overline{\bfs}) \subseteq \sfcon(\overline{\bft})$ holds.
The reverse inclusion $\sfcon(\overline{\bfs}) \supseteq \sfcon(\overline{\bft})$ holds by a symmetrical argument.

(ii) This follows from part~(i) since $\sfcon(\overline{\bfu_2}) = \sfcon(\overline{\bfv_2})$ by Lemma~\ref{Lem: B0 content}(iii).
\end{proof}

\begin{lemma} \label{Lem: B0 yz}
\begin{enumerate}[\ \rm(i)]
\item Suppose that $y z \hookrightarrow \bfu_2$ for some $y,z \in \sfcon(\bfs_0 \bfs_1 \cdots \bfs_p)$.
Then the longest $\{\overline{y}, \overline{z}\}$-{\ssw} of $\bfv_2$ can only be $y z$ or $z^* y^*$.
\item Suppose that $y z \hookrightarrow \bfv_2$ for some $y,z \in \sfcon(\bft_0 \bft_1 \cdots \bft_q)$.
Then the longest $\{\overline{y}, \overline{z}\}$-{\ssw} of $\bfu_2$ can only be $y z$ or $z^* y^*$.
\end{enumerate}
\end{lemma}

\begin{proof}
(i) By assumption, \[ \bfu = \bfu_1 \cdot x \cdot \underbrace{\bfa y \bfb z \bfe}_{\bfu_2} \cdot x^* \] for some $\bfa, \bfb, \bfe \in \Fmon$ such that $\bfu_1, \bfa,\bfb,\bfe, x,y,z$ are pairwise disjoint.
Since $\overline{y}, \overline{z} \in \sfcon(\overline{\bfs_0 \bfs_1 \cdots \bfs_p}) = \sfcon(\overline{\bft_0 \bft_1 \cdots \bft_q})$ by Lemma~\ref{Lem: B0 con s c}(i), we have $\sfocc(\overline{y}, \overline{\bfv})=\sfocc(\overline{z}, \overline{\bfv})=1$ by conditions (\Bsf2) and (\Bsf3).
Therefore, the longest $\{\overline{y}, \overline{z}\}$-{\ssw} of $\bfv_2$ can only be one of \[y z, \quad y z^*, \quad y^* z,  \quad y^* z^*, \quad z y, \quad zy^*, \quad z^* y,  \quad z^* y^*.\]
There are four cases to consider.

\medskip

\noindent{\sc Case~1:} $y z^* \hookrightarrow \bfv_2$.
Then \[ \bfv = \bfu_1 \cdot x \cdot \underbrace{\bff y \bfg z^* \bfh}_{\bfv_2} \cdot\, x^* \] for some $\bff, \bfg, \bfh \in \Fmon$ such that $\bfu_1, \bff,\bfg,\bfh, x,y,z$ are pairwise disjoint.
Under the substitution $\gamma_\bfv^y : \cX \to B_0$ in Lemma~\ref{Lem: B0 substitution}, we have \[ \gamma_\bfv^y(\bfv) = \gamma_\bfv^y(\bfu_1x\bff) \cdot \gamma_\bfv^y(y) \cdot \gamma_\bfv^y(\bfg z^* \bfh x^*) = \rme \cdot \rma \cdot \rmf = \rma. \]
But since $\gamma_\bfv^y(y) =\rma$ and $\gamma_\bfv^y(z) =\rme$, we deduce the contradiction \[ \gamma_\bfv^y(\bfu) = \gamma_\bfv^y(\bfu_1 x \bfa) \cdot \rma \cdot \gamma_\bfv^y(\bfb) \cdot \rme \cdot \gamma_\bfv^y(\bfe x^*) = 0. \]

\medskip

\noindent{\sc Case~2:} $z y^* \hookrightarrow \bfv_2$.
Then \[ \bfv = \bfu_1 \cdot x \cdot \underbrace{\bff z \bfg y^* \bfh}_{\bfv_2} \cdot\, x^* \] for some $\bff, \bfg, \bfh \in \Fmon$ such that $\bfu_1, \bff,\bfg,\bfh, x,y,z$ are pairwise disjoint.
Under the substitution $\gamma_\bfv^{y^*} \! : \cX \to B_0$ in Lemma~\ref{Lem: B0 substitution}, we have \[ \gamma_\bfv^{y^*}\!(\bfv) = \gamma_\bfv^{y^*}\!(\bfu_1x\bff z \bfg) \cdot \gamma_\bfv^{y^*}\!(y^*) \cdot \gamma_\bfv^{y^*}\!(\bfh x^*) = \rme \cdot \rma \cdot \rmf = \rma. \]
But since $\gamma_\bfv^{y^*}\!(y) =\rma$ and $\gamma_\bfv^{y^*}\!(z) =\rme$, we deduce the contradiction \[ \gamma_\bfv^{y^*}\!(\bfu) = \gamma_\bfv^{y^*}(\bfu_1 x \bfa) \cdot \rma \cdot \gamma_\bfv^{y^*}\!(\bfb) \cdot \rme \cdot \gamma_\bfv^{y^*}\!(\bfe x^*) = 0. \]

\medskip

\noindent{\sc Case~3:} $y^* z \hookrightarrow \bfv_2$ or $z y \hookrightarrow \bfv_2$.
Under the substitution $\gamma_\bfu^z : \cX \to B_0$ in Lemma~\ref{Lem: B0 substitution}, we have \[ \gamma_\bfu^z(\bfu) = \gamma_\bfu^z(\bfu_1 x \bfa y \bfb) \cdot \gamma_\bfu^z(z) \cdot \gamma_\bfu^z(\bfe x^*) = \rme \cdot \rma \cdot \rmf = \rma. \]
But since $\gamma_\bfu^z(y) =\rme$ and $\gamma_\bfu^z(z) =\rma$, we deduce the contradiction
\begin{align*}
\gamma_\bfu^z(\bfv) & = \begin{cases} \cdots \rme^* \cdots \rma \cdots & \text{if $y^* z \hookrightarrow \bfv_2$} \\ \cdots \rma \cdots \rme \cdots & \text{if $zy \hookrightarrow \bfv_2$} \end{cases} \\ & = 0.
\end{align*}

\medskip

\noindent{\sc Case~4:} $z^* y \hookrightarrow \bfv_2$ or $y^* z^* \hookrightarrow \bfv_2$.
Then \[ \bfv = \bfu_1 \cdot x \cdot \underbrace{\bff z^* \bfg y \bfh}_{\bfv_2} \cdot\, x^* \quad \text{or} \quad \bfv = \bfu_1 \cdot x \cdot \underbrace{\bff y^* \bfg z^* \bfh}_{\bfv_2} \cdot\, x^* \] for some $\bff,\bfg,\bfh \in \Fmon$ such that $\bfu_1, \bff,\bfg,\bfh, x,y,z$ are pairwise disjoint.
Under the substitution $\gamma_\bfv^{z^*} \! : \cX \to B_0$ in Lemma~\ref{Lem: B0 substitution}, we have
\begin{align*}
\gamma_\bfv^{z^*}\!(\bfv) & =
\begin{cases} \gamma_\bfv^{z^*}\!(\bfu_1 x \bff) \cdot \gamma_\bfv^{z^*}\!(z^*) \cdot \gamma_\bfv^{z^*}\!(\bfg y \bfh x^*) & \text{if $z^* y \hookrightarrow \bfv_2$} \\
\gamma_\bfv^{z^*}\!(\bfu_1 x \bff y^* \bfg) \cdot \gamma_\bfv^{z^*}\!(z^*) \cdot \gamma_\bfv^{z^*}\!(\bfh x^*) & \text{if $y^* z^* \hookrightarrow \bfv_2$}
\end{cases} \\ & = \rme \cdot \rma \cdot \rmf = \rma.
\end{align*}
But since $\gamma_\bfv^{z^*}\!(y) =\rmf$ and $\gamma_\bfv^{z^*}\!(z) =\rma$, we deduce the contradiction \[ \gamma_\bfv^{z^*}\!(\bfu) = \gamma_\bfv^{z^*}\!(\bfu_1 x \bfa) \cdot \rmf \cdot \gamma_\bfv^{z^*}\!(\bfb) \cdot \rma \cdot \gamma_\bfv^{z^*}\!(\bfe x^*) = 0. \]

\medskip

Since none of the four cases is possible, the longest $\{\overline{y},\overline{z}\}$-{\ssw}s of $\bfv_2$ can only be $yz$ or $z^*y^*$.

(ii) This is symmetrical to part~(i).
\end{proof}

\begin{lemma} \label{Lem: B0 yyz}
\begin{enumerate}[\ \rm(i)]
\item Suppose that $y^2z \hookrightarrow \bfu_2$ for some $y \in \sfcon(\bfc_1 \bfc_2 \cdots \bfc_p)$ and $z \in \sfcon(\bfs_0 \bfs_1 \cdots \bfs_p)$.
Then the longest $\{\overline{y},\overline{z}\}$-{\ssw} of $\bfv_2$ can only be $y^2z$ or $z^* (y^*)^2$.
\item Suppose that $y^2z \hookrightarrow \bfv_2$ for some $y \in \sfcon(\bfd_1 \bfd_2 \cdots \bfd_q)$ and $z \in \sfcon(\bft_0 \bft_1 \cdots \bft_q)$.
Then the longest $\{\overline{y},\overline{z}\}$-{\ssw} of $\bfu_2$ can only be $y^2z$ or $z^* (y^*)^2$.
\end{enumerate}
\end{lemma}

\begin{proof}
(i) Since the prefix $\bfs_0$ of $\bfu_2$ consists of simple variables of~$\bfu$, it follows from the assumption that \[ \bfu = \bfu_1 \cdot x \cdot \underbrace{\bfa y^2 \bfb z \bfe}_{\bfu_2} \cdot\, x^*\] for some $\bfa \in \Fsem$ and $\bfb, \bfe \in \Fmon$ such that $\bfu_1, \bfa,\bfb,\bfe,x,y,z$ are pairwise disjoint.
Since $\overline{z} \in \sfcon(\overline{\bfs_0 \bfs_1 \cdots \bfs_p}) = \sfcon(\overline{\bft_0 \bft_1 \cdots \bft_q})$ and $\overline{y} \in \sfcon(\overline{\bfc_1 \bfc_2 \cdots \bfc_p}) = \sfcon(\overline{\bfd_1 \bfd_2 \cdots \bfd_q})$ by Lemma~\ref{Lem: B0 con s c}, we have $\sfocc(\overline{y}, \overline{\bfv_2})=2$ and $\sfocc(\overline{z}, \overline{\bfv_2})=1$.
Hence the longest $\{\overline{y}, \overline{z}\}$-{\ssw} of $\bfv_2$ can only be one of \[ y^2z, \quad y^2z^*, \quad (y^*)^2z, \quad (y^*)^2z^*, \quad z y^2, \quad z (y^*)^2, \quad  z^* y^2,  \quad z^* (y^*)^2. \]
There are two cases to consider.

\medskip

\noindent{\sc Case~1:} $y^2z^* \hookrightarrow \bfv_2$ or $z (y^*)^2 \hookrightarrow \bfv_2$.
The variable $h = \sfh(\bfa)$ is simple in~$\bfu$, so that $\bfa = h\bff$ for some $\bff \in \Fmon$.
Then under the substitution $\gamma_\bfu^h : \cX \to B_0$ in Lemma~\ref{Lem: B0 substitution}, we have \[ \gamma_\bfu^h(\bfu) = \gamma_\bfu^h(\bfu_1 x) \cdot \gamma_\bfu^h(h) \cdot \gamma_\bfu^h(\bff y^2 \bfb z \bfe x^*) = \rme \cdot \rma \cdot \rmf = \rma. \]
But since $\gamma_\bfu^h(y) = \gamma_\bfu^h(z) =\rmf$, we deduce the contradiction
\begin{align*}
\gamma_\bfu^h(\bfv) & = \begin{cases} \cdots \rmf^2 \cdots \rmf^* \cdots & \text{if $y^2z^* \hookrightarrow \bfv_2$} \\ \cdots \rmf \cdots (\rmf^*)^2 \cdots & \text{if $z (y^*)^2 \hookrightarrow \bfv_2$} \end{cases} \\ & = 0.
\end{align*}

\medskip

\noindent{\sc Case~2:} $(y^*)^2 z \hookrightarrow \bfv_2$ or $(y^*)^2 z^* \hookrightarrow \bfv_2$ or $z y^2 \hookrightarrow \bfv_2$ or $z^* y^2 \hookrightarrow \bfv_2$.
Then under the substitution $\gamma_\bfu^z : \cX \to B_0$ in Lemma~\ref{Lem: B0 substitution}, we have \[ \gamma_\bfu^z(\bfu) = \gamma_\bfu^z(\bfu_1 x \bfa y^2 \bfb) \cdot \gamma_\bfu^z(z) \cdot \gamma_\bfu^z(\bfe x^*) = \rme \cdot \rma \cdot \rmf = \rma. \]
But since $\gamma_\bfu^z(y) =\rme$ and $\gamma_\bfu^z(z) =\rma$, we deduce the contradiction
\begin{align*}
\gamma_\bfu^z(\bfv) & =
\begin{cases}
\cdots (\rme^*)^2 \cdots \rma \cdots & \text{if $(y^*)^2 z \hookrightarrow \bfv_2$ or $(y^*)^2 z^* \hookrightarrow \bfv_2$} \\
\cdots \rma \cdots \rme^2 \cdots & \text{if $z y^2 \hookrightarrow \bfv_2$ or $z^* y^2 \hookrightarrow \bfv_2$}
\end{cases} \\ & = 0.
\end{align*}

\medskip

Since none of the two cases is possible, the longest $\{\overline{y},\overline{z}\}$-{\ssw}s of $\bfv_2$ can only be $y^2z$ or $z^* (y^*)^2$.

(ii) This is symmetrical to part~(i).
\end{proof}

\begin{lemma} \label{Lem: B0 yzz}
\begin{enumerate}[\ \rm(i)]
\item Suppose that $yz^2 \hookrightarrow \bfu_2$ for some $y \in \sfcon(\bfs_0 \bfs_1 \cdots \bfs_p)$ and $z \in \sfcon(\bfc_1 \bfc_2 \cdots \bfc_p)$.
Then the longest $\{\overline{y},\overline{z}\}$-{\ssw} of $\bfv_2$ can only be $yz^2$ or $(z^*)^2y^*$.
\item Suppose that $yz^2 \hookrightarrow \bfv_2$ for some $y \in \sfcon(\bft_0 \bft_1 \cdots \bft_q)$ and $z \in \sfcon(\bfd_1 \bfd_2 \cdots \bfd_q)$.
Then the longest $\{\overline{y},\overline{z}\}$-{\ssw} of $\bfu_2$ can only be $yz^2$ or $(z^*)^2y^*$.
\end{enumerate}
\end{lemma}

\begin{proof}
This is very similar to the proof of Lemma~\ref{Lem: B0 yyz}, but details are given for the sake of completeness.

(i) Since the suffix $\bfs_p$ of $\bfu_2$ consists of simple variables of~$\bfu$, it follows from the assumption that \[ \bfu = \bfu_1 \cdot x \cdot \underbrace{\bfa y \bfb z^2 \bfe}_{\bfu_2} \cdot\, x^*\] for some $\bfa,\bfb \in \Fmon$ and $\bfe \in \Fsem$ such that $\bfu_1, \bfa,\bfb,\bfe,x,y,z$ are pairwise disjoint.
Since $\overline{y} \in \sfcon(\overline{\bfs_0 \bfs_1 \cdots \bfs_p}) = \sfcon(\overline{\bft_0 \bft_1 \cdots \bft_q})$ and $\overline{z} \in \sfcon(\overline{\bfc_1 \bfc_2 \cdots \bfc_p}) = \sfcon(\overline{\bfd_1 \bfd_2 \cdots \bfd_q})$ by Lemma~\ref{Lem: B0 con s c}, we have $\sfocc(\overline{y}, \overline{\bfv_2})=1$ and $\sfocc(\overline{z}, \overline{\bfv_2})=2$.
Hence the longest $\{\overline{y}, \overline{z}\}$-{\ssw} of $\bfv_2$ can only be one of \[ yz^2, \quad y(z^*)^2, \quad y^*z^2, \quad y^*(z^*)^2, \quad z^2y, \quad z^2y^*, \quad (z^*)^2y,  \quad (z^*)^2y^*. \]
There are two cases to consider.

\medskip

\noindent{\sc Case~1:} $y^*z^2 \hookrightarrow \bfv_2$ or $(z^*)^2y \hookrightarrow \bfv_2$.
The variable $t = \sft(\bfe)$ is simple in~$\bfu$, so that $\bfe = \bff t$ for some $\bff \in \Fmon$.
Then under the substitution $\gamma_\bfu^t : \cX \to B_0$ in Lemma~\ref{Lem: B0 substitution}, we have \[ \gamma_\bfu^t(\bfu) = \gamma_\bfu^t(\bfu_1 x \bfa y \bfb z^2 \bff) \cdot \gamma_\bfu^t(t) \cdot \gamma_\bfu^t(x^*) = \rme \cdot \rma \cdot \rmf = \rma. \]
But since $\gamma_\bfu^t(y) = \gamma_\bfu^t(z) =\rme$, we deduce the contradiction
\begin{align*}
\gamma_\bfu^t(\bfv) & = \begin{cases} \cdots \rme^* \cdots \rme^2 \cdots & \text{if $y^*z^2 \hookrightarrow \bfv_2$} \\ \cdots (\rme^*)^2 \cdots \rme \cdots & \text{if $(z^*)^2y \hookrightarrow \bfv_2$} \end{cases} \\ & = 0.
\end{align*}

\medskip

\noindent{\sc Case~2:} $y(z^*)^2 \hookrightarrow \bfv_2$ or $y^*(z^*)^2 \hookrightarrow \bfv_2$ or $z^2y \hookrightarrow \bfv_2$ or $z^2y^* \hookrightarrow \bfv_2$.
Then under the substitution $\gamma_\bfu^y : \cX \to B_0$ in Lemma~\ref{Lem: B0 substitution}, we have \[ \gamma_\bfu^y(\bfu) = \gamma_\bfu^y(\bfu_1 x \bfa) \cdot \gamma_\bfu^y(y) \cdot \gamma_\bfu^y(\bfb z^2 \bfe x^*) = \rme \cdot \rma \cdot \rmf = \rma. \]
But since $\gamma_\bfu^y(y) =\rma$ and $\gamma_\bfu^y(z) =\rmf$, we deduce the contradiction
\begin{align*}
\gamma_\bfu^y(\bfv) & =
\begin{cases}
\cdots \rma \cdots (\rmf^*)^2 \cdots & \text{if $y(z^*)^2 \hookrightarrow \bfv_2$ or $y^*(z^*)^2 \hookrightarrow \bfv_2$} \\
\cdots \rmf^2 \cdots \rma \cdots & \text{if $z^2y \hookrightarrow \bfv_2$ or $z^2y^* \hookrightarrow \bfv_2$}
\end{cases} \\ & = 0.
\end{align*}

\medskip

Since none of the two cases is possible, the longest $\{\overline{y},\overline{z}\}$-{\ssw}s of $\bfv_2$ can only be $yz^2$ or $(z^*)^2y^*$.

(ii) This is symmetrical to part~(i).
\end{proof}

\begin{lemma} \label{Lem: B0 head tail}
$\sfh(\bfu_2) =\sfh(\bfv_2)$ and $\sft(\bfu_2) =\sft(\bfv_2)$.
\end{lemma}

\begin{proof}
Recall that $\sfcon(\overline{\bfu_2}) = \sfcon(\overline{\bfv_2})$ by Lemma~\ref{Lem: B0 content}(iii).
First, suppose that $|\sfcon(\overline{\bfu_2})| = |\sfcon(\overline{\bfv_2})| = 1$, say $\sfcon(\overline{\bfu_2}) = \sfcon(\overline{\bfv_2}) = \{ z \}$ for some $z \in \cX$.
Then it follows from condition (\Bsf\ref{Bsf4 singleton}) that $\bfu_2=\bfv_2 =z$, whence $\sfh(\bfu_2)= z =\sfh(\bfv_2)$ and $\sft(\bfu_2)= z =\sft(\bfv_2)$.

Hence, it remains to assume $|\sfcon(\overline{\bfu_2})| = |\sfcon(\overline{\bfv_2})| \geq 2$.
Let $\headu = \sfh(\bfu_2) = \sfh(\bfs_0)$, $\tailu=\sft(\bfu_2) = \sfh(\bfs_p)$, $\headv=\sfh(\bfv_2) = \sfh(\bft_0)$, and $\tailv=\sft(\bfv_2) = \sfh(\bft_q)$, so that
\begin{equation}
\overline{\headu} \order \overline{\tailu} \quad \text{and} \quad \overline{\headv} \order \overline{\tailv} \label{Dis: B0 h<t H<T}
\end{equation}
by condition (\Bsf\ref{Bsf4 s0 orderedsquare}).
Then
\begin{equation}
\bfu = \bfu_1 \cdot x \cdot \underbrace{\headu \bfa \tailu}_{\bfu_2} \cdot\, x^* \quad \text{and} \quad \bfv = \bfu_1 \cdot x \cdot \underbrace{\headv \bfb \tailv}_{\bfv_2} \cdot\, x^* \label{Dis: B0 u2 v2}
\end{equation}
for some $\bfa,\bfb \in \Fmon$ such that $\overline{\headu\tailu} \notin \sfcon(\overline{\bfu_1x\bfa})$ and $\overline{\headv\tailv} \notin \sfcon(\overline{\bfu_1x\bfb})$.
Since $\headu\tailu \hookrightarrow \bfu_2$ with $\headu,\tailu \in \sfcon(\bfs_0\bfs_p)$, it follows from Lemma~\ref{Lem: B0 yz}(i) that the longest $\{\overline{\headu},\overline{\tailu}\}$-{\ssw} of $\bfv_2$ is $\headu\tailu$ or $\tailu^*\headu^*$.

Seeking a contradiction, suppose that $\tailu^*\headu^*$ is the longest $\{\overline{\headu},\overline{\tailu}\}$-{\ssw} of~$\bfv_2$.
If $\headv \neq \tailu^*$, so that $\headv \tailu^* \hookrightarrow \bfv_2$, then by Lemma~\ref{Lem: B0 yz}(ii), either $\headv \tailu^* \hookrightarrow \bfu_2$ or $\tailu \headv^* \hookrightarrow \bfu_2$; but neither {\ssw} is possible in view of \eqref{Dis: B0 u2 v2}.
Hence $\headv = \tailu^*$.
By a symmetrical argument, we deduce $\tailv = \headu^*$.
Since $\overline{\headu} \order \overline{\tailu}$ by \eqref{Dis: B0 h<t H<T}, we have $\overline{\tailv} = \overline{\headu} \order \overline{\tailu} = \overline{\headv}$; but this contradicts $\overline{\headv} \order \overline{\tailv}$ from \eqref{Dis: B0 h<t H<T}.

Therefore, $\headu\tailu$ is the longest $\{\overline{\headu},\overline{\tailu}\}$-{\ssw} of~$\bfv_2$.
If $\headv \neq \headu$, so that $\headv\headu \hookrightarrow \bfv_2$, then by Lemma~\ref{Lem: B0 yz}(ii), either $\headv \headu \hookrightarrow \bfu_2$ or $\headu^* \headv^* \hookrightarrow \bfu_2$; but neither {\ssw} is possible in view of \eqref{Dis: B0 u2 v2}.
Hence $\headv = \headu$.
By a symmetrical argument, we deduce $\tailv = \tailu$.
\end{proof}

\begin{lemma}\label{Lem: B0 u2=v2}
$\bfu_2=\bfv_2$.
\end{lemma}

\begin{proof}
Recall that \[ \bfu_2 = \bfs_0 \prod_{i=1}^p (\bfc_i\bfs_i) \quad \text{and} \quad \bfv_2 = \bft_0 \prod_{i=1}^q (\bfd_i\bft_i), \] and $\sfcon(\overline{\bfu_2}) = \sfcon(\overline{\bfv_2})$ by Lemma~\ref{Lem: B0 content}(iii).
If $|\sfcon(\overline{\bfu_2})| = |\sfcon(\overline{\bfv_2})| = 1$, then as shown in the proof of Lemma~\ref{Lem: B0 head tail}, we have $\bfu_2= \bfv_2$.
Therefore, it suffices to assume that $|\sfcon(\overline{\bfu_2})| = |\sfcon(\overline{\bfv_2})| \geq 2$.
By Lemma~\ref{Lem: B0 head tail}, we have
\begin{enumerate}[\ (a)]
\item[(a)] $\sfh(\bfs_0) = \sfh(\bfu_2) = \sfh(\bfv_2) = \sfh(\bft_0)$ and $\sft(\bfs_p) = \sft(\bfu_2) = \sft(\bfv_2) = \sft(\bft_q)$.
\end{enumerate}
Suppose that $z \in \sfcon(\bfs_0\bfs_1 \cdots \bfs_p)$ with $z \neq \sfh(\bfs_0), \sft(\bfs_p)$, so that $\sfh(\bfs_0) z \hookrightarrow \bfu_2$.
Then by Lemma~\ref{Lem: B0 yz}(i), the longest $\{\overline{\sfh(\bfs_0)}, \overline{z}\}$-{\ssw} of $\bfv_2$ can only be $\sfh(\bfs_0) z$ or $z^*(\sfh(\bfs_0))^*$.
But since~$\bfv_2$ is a bipartite word (see Remark~\ref{Rmk: B0-standard}(iv)) that contains the variable $\sfh(\bft_0) = \sfh(\bfs_0)$, it cannot contain the variable $(\sfh(\bfs_0))^*$.
Hence the longest $\{\overline{\sfh(\bfs_0)}, \overline{z}\}$-{\ssw} of $\bfv_2$ must be $\sfh(\bfs_0) z$, so that $z \in \sfcon(\bft_0\bft_1 \cdots \bft_q)$.
Therefore, the inclusion $\sfcon(\bfs_0\bfs_1 \cdots \bfs_p) \subseteq \sfcon(\bft_0\bft_1 \cdots \bft_q)$ holds.
The reverse inclusion $\sfcon(\bfs_0\bfs_1 \cdots \bfs_p) \supseteq \sfcon(\bft_0\bft_1 \cdots \bft_q)$ is established by a symmetrical argument, so that $\sfcon(\bfs_0\bfs_1 \cdots \bfs_p) = \sfcon(\bft_0\bft_1 \cdots \bft_q)$.
Further, since $\sfh(\bfs_0) = \sfh(\bft_0)$ and $\sft(\bfs_p) = \sft(\bft_q)$ by~(a), it is easy to show by Lemma~\ref{Lem: B0 yz} that $\bfs_0\bfs_1 \cdots \bfs_p = \bft_0\bft_1 \cdots \bft_q$.
Hence
\begin{enumerate}[\ (a)]
\item[(b)] $\bfs_0\bfs_1 \cdots \bfs_p = \bft_0\bft_1 \cdots \bft_q = z_1z_2 \cdots z_r$ for some distinct $z_1,z_2,\ldots,z_r \in \cXX$,
\end{enumerate}
where $z_1 = \sfh(\bfs_0) = \sfh(\bft_0)$ and $z_r = \sft(\bfs_p) = \sft(\bft_q)$.

Now it follows from Lemma~\ref{Lem: B0 con s c} that $p=0$ if and only if $q=0$, so there are two cases: $p=q=0$ and $p,q \geq 1$.
If $p=q=0$, then $\bfu_2 = \bfs_0 = \bft_0 = \bfv_2$, so the proof is complete.
Hence it remains to assume $p,q \geq 1$.

Seeking a contradiction, suppose that $\bfs_0 \neq \bft_0$.
Then by~(b), either $\bfs_0$ is a proper prefix of $\bft_0$ or $\bft_0$ is a proper prefix of $\bfs_0$.
By symmetry, suppose that $\bfs_0$ is a proper prefix of $\bft_0$, so that $\bfs_0 = z_1z_2 \cdots z_k$ and $\bft_0 = z_1z_2 \cdots z_\ell$ with $k < \ell \leq r$.
Then
\begin{align*}
\bfu_2 & = z_1z_2 \cdots z_k \cdot \bfc_1 \cdot z_{k+1} z_{k+2} \cdots \\ \text{and} \quad \bfv_2 & = z_1z_2 \cdots z_k \cdot z_{k+1} z_{k+2} \cdots z_\ell \cdot \bfd_1 \cdot z_{\ell+1} z_{\ell+2} \cdots.
\end{align*}
Since~$\bfc_1$ is an {\Bblock}, it begins with $y^2$ for some $y \in \cXX$.
Then $y^2z_{k+1} \hookrightarrow \bfu_2$ but $y^2z_{k+1} \not\hookrightarrow \bfv_2$ and $z_{k+1}^*(y^*)^2 \not\hookrightarrow \bfv_2$, which is impossible in view of Lemma~\ref{Lem: B0 yyz}(i).
Therefore, $\bfs_0 = z_1z_2 \cdots z_k = \bft_0$, so that
\begin{align*}
\bfu_2 & = z_1z_2 \cdots z_k \cdot \bfc_1 \cdot z_{k+1} z_{k+2} \cdots \\ \text{and} \quad \bfv_2 & = z_1z_2 \cdots z_k \cdot \bfd_1 \cdot z_{k+1} z_{k+2} \cdots.
\end{align*}

Seeking a contradiction, suppose that $\sfcon(\bfc_1) \neq \sfcon(\bfd_1)$, say $y \in \sfcon(\bfc_1) \backslash \sfcon(\bfd_1)$.
Then since~$\bfc_1$ is an {\Bblock}, $\bfc_1 = \bfa y^2 \bfb$ for some $\bfa,\bfb \in \Fmon$.
Hence $y^2z_{k+1} \hookrightarrow \bfu_2$ but $y^2z_{k+1} \not\hookrightarrow \bfv_2$ and $z_{k+1}^*(y^*)^2 \not\hookrightarrow \bfv_2$, which is impossible in view of Lemma~\ref{Lem: B0 yyz}(i).
Therefore, $\sfcon(\bfc_1) = \sfcon(\bfd_1)$.
Since $\bfc_1$ and $\bfd_1$ are {\Bblock}s, we have $\bfc_1 = \bfd_1$.

Without loss of generality, assume that $p \leq q$.
The arguments in the previous two paragraphs can be repeated to show that  $\bfc_i=\bfd_i$ and $\bfs_i = \bft_i$ and for all $i = 1,2,\ldots,p-1$.
Hence \[ \bfu_2 = \bfs_0 \bigg(\prod_{i=1}^{p-1} (\bfc_i\bfs_i) \bigg) \bfc_p \bfs_p \quad \text{and} \quad \bfv_2 = \bfs_0 \bigg(\prod_{i=1}^{p-1} (\bfc_i\bfs_i) \bigg) \bfd_p \bft_p \bfd_{p+1} \bft_{p+1} \cdots \bfd_q \bft_q. \]
Arguments that are dual to those from the previous two paragraphs (with the use of Lemma~\ref{Lem: B0 yzz} instead of Lemma~\ref{Lem: B0 yyz}) can be repeated to show that $\bfs_p = \bft_q$ and then $\bfc_p = \bfd_q$.
It follows from~(b) that $p=q$ and $\bfu_2 = \bfv_2$.
\end{proof}

\section{{\Insem}s of order up to four} \label{sec: tables}

\newcommand{\horbar}{\\\underline{\hspace{0.35in}}\\}
\newcommand{\horbarii}{\\\underline{\hspace{0.2in}}\\}
\newcommand{\horbariii}{\\\underline{\hspace{0.3in}}\\}

Multiplication tables of {\insem}s of order up to four are given in this section.
For a more compact presentation, the column/row headers are omitted from each multiplication table.
For instance, the {\insem} $S=\{ 1,2,3\}$ given by the multiplication table
\[
\begin{tabular}[b]{c | c c c}
$S$   & 1 & 2 & 3 \\ \hline
1     & 1 & 2 & 3 \\
2     & 2 & 3 & 1 \\
3     & 3 & 1 & 2 \\ \hline\hline
$x$   & 1 & 2 & 3 \\ \hline
$x^*$ & 1 & 3 & 2
\end{tabular}
\]
is abbreviated to
\[
\begin{tabular}[b]{c c c}
1 & 2 & 3 \\
2 & 3 & 1 \\
3 & 1 & 2 \\ \hline
1 & 3 & 2
\end{tabular}
\]

\subsection{{\Insem}s of order up to three}

Up to isomorphism, there are three {\insem}s of order two and 15 {\insem}s of order three, all of which are commutative \cite[Section~4]{Lee22}.
\begin{align*}
&
\begin{array}[c]{llllllllll}
\shortstack{ 1 1 \\ 1 1 \horbarii 1 2} &
\shortstack{ 1 1 \\ 1 2 \horbarii 1 2} &
\shortstack{ 1 2 \\ 2 1 \horbarii 1 2} &
& & & & & & \\[0.1in]
\shortstack{ 1 1 1 \\ 1 1 1 \\ 1 1 1 \horbariii 1 2 3} &
\shortstack{ 1 1 1 \\ 1 1 1 \\ 1 1 1 \horbariii 1 3 2} &
\shortstack{ 1 1 1 \\ 1 1 1 \\ 1 1 2 \horbariii 1 2 3} &
\shortstack{ 1 1 1 \\ 1 1 1 \\ 1 1 3 \horbariii 1 2 3} &
\shortstack{ 1 1 1 \\ 1 1 2 \\ 1 2 3 \horbariii 1 2 3} &
\shortstack{ 1 1 1 \\ 1 2 1 \\ 1 1 3 \horbariii 1 2 3} &
\shortstack{ 1 1 1 \\ 1 2 1 \\ 1 1 3 \horbariii 1 3 2} &
\shortstack{ 1 1 1 \\ 1 2 2 \\ 1 2 2 \horbariii 1 2 3} &
\shortstack{ 1 1 1 \\ 1 2 2 \\ 1 2 3 \horbariii 1 2 3} &
\shortstack{ 1 1 1 \\ 1 2 3 \\ 1 3 2 \horbariii 1 2 3} \\[0.1in]
\shortstack{ 1 1 3 \\ 1 1 3 \\ 3 3 1 \horbariii 1 2 3} &
\shortstack{ 1 1 3 \\ 1 2 3 \\ 3 3 1 \horbariii 1 2 3} &
\shortstack{ 1 2 2 \\ 2 1 1 \\ 2 1 1 \horbariii 1 2 3} &
\shortstack{ 1 2 3 \\ 2 3 1 \\ 3 1 2 \horbariii 1 2 3} &
\shortstack{ 1 2 3 \\ 2 3 1 \\ 3 1 2 \horbariii 1 3 2} &
& & & &
\end{array}
\end{align*}

\subsection{{\Insem}s of order four}

Up to isomorphism, there are 83 {\insem}s of order four.
{\scriptsize
\begin{align*}
&
\begin{array}[c]{cccccccccc}
\shortstack{ 1 1 1 1 \\ 1 1 1 1 \\ 1 1 1 1 \\ 1 1 1 1 \horbar 1 2 3 4} &
\shortstack{ 1 1 1 1 \\ 1 1 1 1 \\ 1 1 1 1 \\ 1 1 1 1 \horbar 1 2 4 3} &
\shortstack{ 1 1 1 1 \\ 1 1 1 1 \\ 1 1 1 1 \\ 1 1 1 2 \horbar 1 2 3 4} &
\shortstack{ 1 1 1 1 \\ 1 1 1 1 \\ 1 1 1 1 \\ 1 1 1 4 \horbar 1 2 3 4} &
\shortstack{ 1 1 1 1 \\ 1 1 1 1 \\ 1 1 1 1 \\ 1 1 1 4 \horbar 1 3 2 4} &
\shortstack{ 1 1 1 1 \\ 1 1 1 1 \\ 1 1 1 1 \\ 1 1 2 1 \horbar 1 2 4 3} &
\shortstack{ 1 1 1 1 \\ 1 1 1 1 \\ 1 1 1 2 \\ 1 1 2 1 \horbar 1 2 3 4} &
\shortstack{ 1 1 1 1 \\ 1 1 1 1 \\ 1 1 1 2 \\ 1 1 2 1 \horbar 1 2 4 3} &
\shortstack{ 1 1 1 1 \\ 1 1 1 1 \\ 1 1 1 2 \\ 1 1 2 2 \horbar 1 2 3 4} &
\shortstack{ 1 1 1 1 \\ 1 1 1 1 \\ 1 1 1 2 \\ 1 1 2 3 \horbar 1 2 3 4}
\end{array} \\[0.03in]
&
\begin{array}[c]{cccccccccc}
\shortstack{ 1 1 1 1 \\ 1 1 1 1 \\ 1 1 1 3 \\ 1 1 3 4 \horbar 1 2 3 4} &
\shortstack{ 1 1 1 1 \\ 1 1 1 1 \\ 1 1 1 3 \\ 1 2 1 4 \horbar 1 3 2 4} &
\shortstack{ 1 1 1 1 \\ 1 1 1 1 \\ 1 1 2 1 \\ 1 1 1 2 \horbar 1 2 3 4} &
\shortstack{ 1 1 1 1 \\ 1 1 1 1 \\ 1 1 2 1 \\ 1 1 1 2 \horbar 1 2 4 3} &
\shortstack{ 1 1 1 1 \\ 1 1 1 1 \\ 1 1 2 1 \\ 1 1 1 4 \horbar 1 2 3 4} &
\shortstack{ 1 1 1 1 \\ 1 1 1 1 \\ 1 1 2 1 \\ 1 1 2 2 \horbar 1 2 4 3} &
\shortstack{ 1 1 1 1 \\ 1 1 1 1 \\ 1 1 2 2 \\ 1 1 2 2 \horbar 1 2 3 4} &
\shortstack{ 1 1 1 1 \\ 1 1 1 1 \\ 1 1 2 2 \\ 1 1 2 2 \horbar 1 2 4 3} &
\shortstack{ 1 1 1 1 \\ 1 1 1 1 \\ 1 1 3 1 \\ 1 1 1 4 \horbar 1 2 3 4} &
\shortstack{ 1 1 1 1 \\ 1 1 1 1 \\ 1 1 3 1 \\ 1 1 1 4 \horbar 1 2 4 3}
\end{array} \\[0.03in]
&
\begin{array}[c]{cccccccccc}
\shortstack{ 1 1 1 1 \\ 1 1 1 1 \\ 1 1 3 3 \\ 1 1 3 3 \horbar 1 2 3 4 } &
\shortstack{ 1 1 1 1 \\ 1 1 1 1 \\ 1 1 3 3 \\ 1 1 3 4 \horbar 1 2 3 4 } &
\shortstack{ 1 1 1 1 \\ 1 1 1 1 \\ 1 1 3 4 \\ 1 1 4 3 \horbar 1 2 3 4 } &
\shortstack{ 1 1 1 1 \\ 1 1 1 2 \\ 1 1 1 2 \\ 1 2 2 4 \horbar 1 2 3 4 } &
\shortstack{ 1 1 1 1 \\ 1 1 1 2 \\ 1 1 1 3 \\ 1 2 3 4 \horbar 1 2 3 4 } &
\shortstack{ 1 1 1 1 \\ 1 1 1 2 \\ 1 1 1 3 \\ 1 2 3 4 \horbar 1 3 2 4 } &
\shortstack{ 1 1 1 1 \\ 1 1 1 2 \\ 1 1 2 3 \\ 1 2 3 4 \horbar 1 2 3 4 } &
\shortstack{ 1 1 1 1 \\ 1 1 1 2 \\ 1 1 3 1 \\ 1 2 1 4 \horbar 1 2 3 4 } &
\shortstack{ 1 1 1 1 \\ 1 1 1 2 \\ 1 1 3 3 \\ 1 2 3 4 \horbar 1 2 3 4 } &
\shortstack{ 1 1 1 1 \\ 1 1 1 2 \\ 1 2 3 1 \\ 1 1 1 4 \horbar 1 2 4 3 }
\end{array} \\[0.03in]
&
\begin{array}[c]{cccccccccc}
\shortstack{ 1 1 1 1 \\ 1 1 1 2 \\ 1 2 3 2 \\ 1 1 1 4 \horbar 1 2 4 3 } &
\shortstack{ 1 1 1 1 \\ 1 1 2 2 \\ 1 2 3 3 \\ 1 2 3 3 \horbar 1 2 3 4 } &
\shortstack{ 1 1 1 1 \\ 1 1 2 2 \\ 1 2 3 3 \\ 1 2 3 4 \horbar 1 2 3 4 } &
\shortstack{ 1 1 1 1 \\ 1 1 2 2 \\ 1 2 3 4 \\ 1 2 4 3 \horbar 1 2 3 4 } &
\shortstack{ 1 1 1 1 \\ 1 2 1 1 \\ 1 1 3 1 \\ 1 1 1 4 \horbar 1 2 3 4 } &
\shortstack{ 1 1 1 1 \\ 1 2 1 1 \\ 1 1 3 1 \\ 1 1 1 4 \horbar 1 2 4 3 } &
\shortstack{ 1 1 1 1 \\ 1 2 1 1 \\ 1 1 3 3 \\ 1 1 3 3 \horbar 1 2 3 4 } &
\shortstack{ 1 1 1 1 \\ 1 2 1 1 \\ 1 1 3 3 \\ 1 1 3 4 \horbar 1 2 3 4 } &
\shortstack{ 1 1 1 1 \\ 1 2 1 1 \\ 1 1 3 4 \\ 1 1 4 3 \horbar 1 2 3 4 } &
\shortstack{ 1 1 1 1 \\ 1 2 1 2 \\ 1 1 3 3 \\ 1 2 3 4 \horbar 1 2 3 4 }
\end{array} \\[0.03in]
&
\begin{array}[c]{cccccccccc}
\shortstack{ 1 1 1 1 \\ 1 2 1 2 \\ 1 1 3 3 \\ 1 2 3 4 \horbar 1 3 2 4 } &
\shortstack{ 1 1 1 1 \\ 1 2 2 2 \\ 1 2 2 2 \\ 1 2 2 2 \horbar 1 2 3 4 } &
\shortstack{ 1 1 1 1 \\ 1 2 2 2 \\ 1 2 2 2 \\ 1 2 2 2 \horbar 1 2 4 3 } &
\shortstack{ 1 1 1 1 \\ 1 2 2 2 \\ 1 2 2 2 \\ 1 2 2 3 \horbar 1 2 3 4 } &
\shortstack{ 1 1 1 1 \\ 1 2 2 2 \\ 1 2 2 2 \\ 1 2 2 4 \horbar 1 2 3 4 } &
\shortstack{ 1 1 1 1 \\ 1 2 2 2 \\ 1 2 2 3 \\ 1 2 3 4 \horbar 1 2 3 4 } &
\shortstack{ 1 1 1 1 \\ 1 2 2 2 \\ 1 2 3 2 \\ 1 2 2 4 \horbar 1 2 3 4 } &
\shortstack{ 1 1 1 1 \\ 1 2 2 2 \\ 1 2 3 2 \\ 1 2 2 4 \horbar 1 2 4 3 } &
\shortstack{ 1 1 1 1 \\ 1 2 2 2 \\ 1 2 3 3 \\ 1 2 3 3 \horbar 1 2 3 4 } &
\shortstack{ 1 1 1 1 \\ 1 2 2 2 \\ 1 2 3 3 \\ 1 2 3 4 \horbar 1 2 3 4 }
\end{array} \\[0.03in]
&
\begin{array}[c]{cccccccccc}
\shortstack{ 1 1 1 1 \\ 1 2 2 2 \\ 1 2 3 4 \\ 1 2 4 3 \horbar 1 2 3 4 } &
\shortstack{ 1 1 1 1 \\ 1 2 2 4 \\ 1 2 2 4 \\ 1 4 4 2 \horbar 1 2 3 4 } &
\shortstack{ 1 1 1 1 \\ 1 2 2 4 \\ 1 2 3 4 \\ 1 4 4 2 \horbar 1 2 3 4 } &
\shortstack{ 1 1 1 1 \\ 1 2 3 3 \\ 1 3 2 2 \\ 1 3 2 2 \horbar 1 2 3 4 } &
\shortstack{ 1 1 1 1 \\ 1 2 3 4 \\ 1 3 4 2 \\ 1 4 2 3 \horbar 1 2 3 4 } &
\shortstack{ 1 1 1 1 \\ 1 2 3 4 \\ 1 3 4 2 \\ 1 4 2 3 \horbar 1 2 4 3 } &
\shortstack{ 1 1 1 4 \\ 1 1 1 4 \\ 1 1 1 4 \\ 4 4 4 1 \horbar 1 2 3 4 } &
\shortstack{ 1 1 1 4 \\ 1 1 1 4 \\ 1 1 1 4 \\ 4 4 4 1 \horbar 1 3 2 4 } &
\shortstack{ 1 1 1 4 \\ 1 1 1 4 \\ 1 1 2 4 \\ 4 4 4 1 \horbar 1 2 3 4 } &
\shortstack{ 1 1 1 4 \\ 1 1 1 4 \\ 1 1 3 4 \\ 4 4 4 1 \horbar 1 2 3 4 }
\end{array} \\[0.03in]
&
\begin{array}[c]{cccccccccc}
\shortstack{ 1 1 1 4 \\ 1 1 2 4 \\ 1 2 3 4 \\ 4 4 4 1 \horbar 1 2 3 4 } &
\shortstack{ 1 1 1 4 \\ 1 2 1 4 \\ 1 1 3 4 \\ 4 4 4 1 \horbar 1 2 3 4 } &
\shortstack{ 1 1 1 4 \\ 1 2 1 4 \\ 1 1 3 4 \\ 4 4 4 1 \horbar 1 3 2 4 } &
\shortstack{ 1 1 1 4 \\ 1 2 2 4 \\ 1 2 2 4 \\ 4 4 4 1 \horbar 1 2 3 4 } &
\shortstack{ 1 1 1 4 \\ 1 2 2 4 \\ 1 2 3 4 \\ 4 4 4 1 \horbar 1 2 3 4 } &
\shortstack{ 1 1 1 4 \\ 1 2 3 4 \\ 1 3 2 4 \\ 4 4 4 1 \horbar 1 2 3 4 } &
\shortstack{ 1 1 3 3 \\ 1 1 3 3 \\ 3 3 1 1 \\ 3 3 1 1 \horbar 1 2 3 4 } &
\shortstack{ 1 1 3 3 \\ 1 1 3 3 \\ 3 3 1 1 \\ 3 3 1 2 \horbar 1 2 3 4 } &
\shortstack{ 1 1 3 3 \\ 1 2 3 3 \\ 3 3 1 1 \\ 3 3 1 1 \horbar 1 2 3 4 } &
\shortstack{ 1 1 3 3 \\ 1 2 3 4 \\ 3 3 1 1 \\ 3 4 1 1 \horbar 1 2 3 4 }
\end{array} \\[0.03in]
&
\begin{array}[c]{cccccccccc}
\shortstack{ 1 1 3 3 \\ 1 2 3 4 \\ 3 3 1 1 \\ 3 4 1 2 \horbar 1 2 3 4 } &
\shortstack{ 1 1 3 3 \\ 2 2 4 4 \\ 1 1 3 3 \\ 2 2 4 4 \horbar 1 3 2 4 } &
\shortstack{ 1 1 3 4 \\ 1 1 3 4 \\ 3 3 4 1 \\ 4 4 1 3 \horbar 1 2 3 4 } &
\shortstack{ 1 1 3 4 \\ 1 1 3 4 \\ 3 3 4 1 \\ 4 4 1 3 \horbar 1 2 4 3 } &
\shortstack{ 1 1 3 4 \\ 1 2 3 4 \\ 3 3 4 1 \\ 4 4 1 3 \horbar 1 2 3 4 } &
\shortstack{ 1 1 3 4 \\ 1 2 3 4 \\ 3 3 4 1 \\ 4 4 1 3 \horbar 1 2 4 3 } &
\shortstack{ 1 2 2 2 \\ 2 1 1 1 \\ 2 1 1 1 \\ 2 1 1 1 \horbar 1 2 3 4 } &
\shortstack{ 1 2 2 2 \\ 2 1 1 1 \\ 2 1 1 1 \\ 2 1 1 1 \horbar 1 2 4 3 } &
\shortstack{ 1 2 2 4 \\ 2 4 4 1 \\ 2 4 4 1 \\ 4 1 1 2 \horbar 1 2 3 4 } &
\shortstack{ 1 2 3 4 \\ 2 1 4 3 \\ 3 4 1 2 \\ 4 3 2 1 \horbar 1 2 3 4 }
\end{array} \\[0.03in]
&
\begin{array}[c]{cccccccccc}
\shortstack{ 1 2 3 4 \\ 2 1 4 3 \\ 3 4 1 2 \\ 4 3 2 1 \horbar 1 2 4 3 } &
\shortstack{ 1 2 3 4 \\ 2 1 4 3 \\ 3 4 2 1 \\ 4 3 1 2 \horbar 1 2 3 4 } &
\shortstack{ 1 2 3 4 \\ 2 1 4 3 \\ 3 4 2 1 \\ 4 3 1 2 \horbar 1 2 4 3 } &
& & & & & &
\end{array}
\end{align*}
}

\section{Acknowledgements}

The authors are grateful to the reviewers for their insightful reports and suggestions.

\end{document}